\documentclass[12pt,twoside]{article}
\usepackage[margin=1in]{geometry}
\usepackage{fancyhdr}
\usepackage{graphicx}
\usepackage{subcaption}
\usepackage{amssymb}
\usepackage{amsmath}
\usepackage{amsthm}
\usepackage{amsfonts}
\usepackage{mathrsfs}
\usepackage{mathtools}
\usepackage{bm}
\usepackage{color}
\usepackage{setspace}
\usepackage{exscale}
\usepackage{relsize}
\usepackage{float}
\usepackage{cite}
\usepackage{hyperref}
\usepackage{cleveref}
\usepackage{nicefrac}
\usepackage{ulem}
\usepackage{caption}
\usepackage{soul}
\usepackage{dsfont}
\usepackage{euscript,mathrsfs}
\usepackage{comment}
\usepackage{diagbox}

\usepackage{booktabs,multirow} 
\usepackage{array} 
\usepackage{paralist} 
\usepackage{fancyhdr} 
\theoremstyle{plain}                    
\newtheorem{thm}{Theorem}[section]

\newtheorem{definition}[thm]{Definition}

\newenvironment{acknowledgment}{{\flushleft \bf Acknowledgment:}}{}

\usepackage{tikz}
\usetikzlibrary{positioning}
\usepackage{xcolor}
\allowdisplaybreaks[1]
\numberwithin{equation}{section}
\numberwithin{figure}{section}
\numberwithin{table}{section}
\newcommand\eref[1]{(\ref{#1})}
\newcommand\fref[1]{Figure~\ref{#1}}
\newcommand\tref[1]{Table~\ref{#1}}
\newcommand*\xbar[1]{%
  \hbox{%
    \vbox{%
      \hrule height 0.5pt 
      \kern0.4ex
      \hbox{%
        \kern-0.05em
        \ensuremath{#1}%
        \kern-0.00em
      }%
    }%
  }%
}

\newcommand{\dx}{\Delta x}

\newcommand{\dy}{\Delta y}
\newcommand{\hf}{{\frac{1}{2}}}

\newcommand{\softd}{d\hspace{-0.2mm}'}

\newcommand{\Td}{\mathbb{T}^d}

\newcommand{\R}{\Bbb R}

\newcommand{\vc}[1]{{\bm #1}}

\newcommand{\vm}{\vc{m}}

\newcommand{\Grad}{\nabla}
\newcommand{\bnabla}{{\bm\nabla}}

\newcommand{\bfphi}{\bm{\varphi}}

\graphicspath{{Figures/}}

\hypersetup{colorlinks=true,
            breaklinks=true,
            urlcolor=blue,
            linkcolor=black,
            bookmarksopen=false,
            filecolor=black,
            citecolor=black,
            linkbordercolor=blue
}

\title{Numerical Study of Random Kelvin-Helmholtz Instability}
\author{Alina Chertock\thanks{Department of Mathematics, North Carolina State University, Raleigh, NC, USA;
{\href{mailto:chertock@math.ncsu.edu}{chertock@math.ncsu.edu}}}, Michael Herty\thanks{Department of Mathematics, RWTH Aachen University,
Aachen, Germany; {\href{mailto:herty@igpm.rwth-aachen.de}{herty@igpm.rwth-aachen.de}}}, Arsen S. Iskhakov\thanks{Department of Mechanical
and Nuclear Engineering, Kansas State University, Manhattan, KS, USA; {\href{mailto:aiskhak@ksu.edu}{aiskhak@ksu.edu}}}, Anna
Iskhakova\thanks{Department of Mechanical and Nuclear Engineering, Kansas State University, Manhattan, KS, USA; 
{\href{mailto:aiskh@ksu.edu}{aiskh@ksu.edu}}},\\ Alexander Kurganov\thanks{Department of Mathematics and Shenzhen International Center for 
Mathematics, Southern University of Science and Technology, Shenzhen 518055, China; {\href{mailto:alexander@sustech.edu.cn}
{alexander@sustech.edu.cn}}}, and M\'{a}ria Luk\'a\v{c}ov\'a-Medvi{\softd}ov\'a\thanks{Institute of Mathematics, Johannes Gutenberg 
University Mainz, Mainz, Germany; {\href{mailto:lukacova@uni-mainz.de}{lukacova@uni-mainz.de}}}}
\date{}

\begin{document}
\maketitle
\begin{abstract}
In this paper, we study random dissipative weak solutions of the compressible Euler equations in the Kelvin–Helmholtz (KH) instability.
Motivated by the fact that weak entropy solutions are not unique and can be viewed as inviscid limits of Navier–Stokes flows, we take a
statistical approach following ideas from turbulence theory. Our aim is to identify solution features that remain consistent across
different realizations and mesh resolutions. For this purpose, we compute stable numerical solutions using a stochastic collocation method
implemented with the help of a fifth-order alternative weighted essentially non-oscillatory (A-WENO) scheme and seventh-order central
weighted essentially non-oscillatory (CWENO) interpolation in the random space. The obtained solutions are averaged over several embedded
uniform grids, resulting in Ces\`aro averages, which are studied using stochastic tools. The analysis includes Reynolds stress and energy
defects, probability density functions of averaged quantities, and reduced-order representations using proper orthogonal decomposition. The
presented numerical experiments illustrate that random KH instabilities can be systematically described using statistical methods,
averaging, and reduced-order modeling, providing a robust methodology for capturing the complex and chaotic dynamics of inviscid
compressible flows.
\end{abstract}



\section{Introduction}
We consider the stochastic compressible Euler equations
\begin{equation}
\begin{aligned}
&\rho_t+\bnabla\!\cdot\!\bm m=0,\\
&\bm m_t+\bnabla\!\cdot\!(\bm m\otimes\bm u)+\nabla p=0,\\
&E_t+\bnabla\!\cdot\!((E+p)\bm u)=0,
\end{aligned}
\label{1.1}
\end{equation}
where $\rho(\bm x,t;\bm\xi), \bm m(\bm x,t;\bm\xi)$, and $E(\bm x,t;\bm\xi)$ are the conservative random variables, representing the
density, momentum, and the total energy, respectively. Here, $t$ is time, $\bm x\in\mathbb R^d$ are spatial variables, and
$\bm\xi\in\Omega\subset\mathbb R^s$ are real-valued random variables. We denote a complete probability space
$(\Omega,\mathfrak M,\mathbb P)$, where $\Omega$ is a set of events, $\mathfrak M$ is the $\sigma$-algebra of Borel measurable sets,
$\mathbb P$ is a probability measure, and $\mu(\bm\xi):\mathbb R^s\to\mathbb R^+$ denotes the probability density function (PDF) of
$\bm\xi$. Further, $p$ and $\bm u=\bm m/\rho$ stand for the pressure and velocity. The system \eref{1.1} is closed using the following
equation of state:
\begin{equation}
E=\frac{|\bm m|^2}{2\rho}+\rho e,\quad e=c_V\rho^{\gamma-1}\exp\left(\frac{S}{c_V\rho}\right),\quad
S=c_V\rho\ln\left(\frac{p}{\rho^\gamma}\right),
\label{1.2}
\end{equation}
where $e$ is the internal energy, $S$ is the total entropy, $1<\gamma\le5/3$ is the adiabatic coefficient, and $c_V=\frac{1}{\gamma-1}$ is
the specific heat at constant volume.

In addition to the conservation laws \eref{1.1}, we impose the second law of thermodynamics, expressed by the entropy inequality requiring
that entropy is nondecreasing in time:
\begin{equation}
S_t+\bnabla\!\cdot\!(S\bm u)\ge0.
\label{1.3}
\end{equation}
This condition serves as an admissibility criterion for weak solutions of the Euler system, ruling out nonphysical states. Nevertheless,
even under the entropy inequality, the multidimensional Euler equations may admit infinitely many weak entropy solutions, which is the
source of ill-posedness discussed below. 

Even in the deterministic case, that is, when $\rho=\rho(\bm x,t)$, $\bm m=\bm m(\bm x,t)$, and $E=E(\bm x,t)$, solutions of
\eref{1.1}--\eref{1.3} are known to develop discontinuities in finite time even for infinitely smooth initial data. Since a classical
solution may not exist, \eref{1.1}--\eref{1.3} are considered in the weak (distributional) sense. However, it was shown in \cite{dlsz2}
that one can construct infinitely many weak entropy solutions of the multidimensional compressible Euler equations; see also
\cite{kreml,feireisl}. Because of the ill-posedness of multidimensional Euler equations in the class of weak entropy solutions, there
is a need to propose new selection criteria to obtain a physically reasonable solution concept. We note that these questions are still 
open and pose challenges for numerical computations. Namely, different numerical methods may potentially produce different results for 
the same specific initial data. Moreover, numerical solutions computed by the same numerical method do not necessarily exhibit strong 
convergence as the mesh is refined; see, e.g., \cite{book,Feireisl2021,FMT_Acta,CCHKL}. Therefore, one may consider approximating suitable 
observable quantities obtained by an averaging procedure. For instance, it was shown in \cite{Feireisl2021}, that the so-called Ces\`aro 
averages computed over several mesh resolutions converge strongly to a generalized solution of the compressible Euler equations.  

The compressible Euler system \eref{1.1}--\eref{1.3} can be regarded as the inviscid limit of the Navier-Stokes equations. In this limit,
the absence of viscous dissipation leads to increasingly fine-scale oscillations, which are characteristic of turbulent flows. Even from
smooth initial data, solutions may evolve into complex structures whose deterministic description is ill-posed, while statistical 
quantities still remain meaningful. This connects the issue of non-uniqueness of weak entropy solutions with the broader question of turbulence modeling.

A prototypical mechanism for the onset of turbulence is the Kelvin-Helmholtz (KH) instability, where a shear layer rolls up into vortical
structures, which undergo secondary instabilities and cascade into progressively smaller scales. This process highlights the difficulty of
predicting a unique deterministic solution, while simultaneously motivating statistical approaches. Indeed, turbulence is commonly 
described not by single realizations but by ensemble or averaged quantities, such as mean fields, variances, and energy spectra, which 
exhibit reproducible behavior.

In this work, we therefore adopt a statistical viewpoint inspired by the turbulence theory. Building on this motivation, we propose a
statistical framework to study the non-uniqueness of the compressible Euler equations by considering the random system
\eref{1.1}--\eref{1.3}. We focus on the KH instability as a representative case and investigate the statistical properties of the resulting
solutions. Our goal is to identify robust features that persist across different realizations, thereby gaining insight into the complex and
potentially chaotic behavior of the system. To this end, we compute numerical solutions using a stochastic collocation method, which belongs
to a class of non-intrusive algorithms, in which one seeks to satisfy the governing equations at a discrete set of nodes in the random space
employing the same numerical solver as for the deterministic problem, and then using interpolation and quadrature rules to evaluate
statistical moments numerically; see, e.g., \cite{Xiu09,Xiu10}. At each collocation point, the deterministic compressible Euler equations
are numerically solved by the fifth-order alternative weighted essentially non-oscillatory (A-WENO) scheme from \cite{Shaoshuai2022} on a
sequence of embedded uniform spatial meshes, and the obtained solutions are used to compute the Ces\`aro averages. The generated data are
then interpolated in the random space using the seventh-order central weighted essentially non-oscillatory (CWENO7) interpolation
\cite{CIK,CPSV,DZP} resulting in a piecewise polynomial approximation, which is, in turn, integrated to compute the statistical moments.

Equipped with the constructed collocation method, we compute Ces\`aro averages, Reynolds stresses and energy defects, and perform
reduced-order analysis, such as proper orthogonal decompositions (POD), to characterize the stochastic solution space of the KH instability.
This perspective is in line with the classical statistical approach to turbulence, where universal features emerge at the level of averaged
or distributional quantities rather than individual flow realizations.

The paper is organized as follows. In \S\ref{sec2}, we introduce the concept of dissipative weak (DW) solutions for both deterministic and
random compressible Euler equations. \S\ref{sec3} describes the numerical methodology, including the computation of quantities of interest.
In \S\ref{sec4}, we present a detailed numerical study of the KH instability, including analysis of Reynolds stress and energy defect,
statistical properties, and reduced-order modeling via POD. Finally, \S\ref{sec5} summarizes our main findings and outlines possible
directions for future research.

\section{Dissipative Weak (DW) Solutions}\label{sec2}
In this section, we describe a concept of dissipative weak solutions for both deterministic (\S\ref{sec21}) and random (\S\ref{sec22})
compressible Euler equations.

\subsection{Deterministic Solutions of the Compressible Euler Equations}\label{sec21}
In view of the ill–posedness of the compressible Euler equations \eref{1.1}--\eref{1.2} in the class of weak entropy solutions, the
relevance of the system \eref{1.1}--\eref{1.2} to describe the behavior of fluids in higher space dimensions may be questionable. In fact,
\eref{1.1}--\eref{1.3} should be seen as an inviscid (vanishing viscosity) limit of a more realistic viscous fluid model. The low viscosity
regime is typical for turbulent flows, where the solutions may develop oscillatory behavior. As it was shown in \cite{book,FLSS}, a weak
limit of weak solutions of the compressible Navier-Stokes equations may not be a weak solution of \eref{1.1}--\eref{1.3}. Instead, it is a
generalized DW solution, which is defined as follows.

Let us consider \eref{1.1}--\eref{1.3} on a space-time cylinder $\Td\times[0,T],$ where $T>0$ and $\Td:=[0,1]^d$, $d=2$, $3$ is a flat
torus, subject to the initial data:
\begin{equation}
\rho(\bm x,0)=\rho_0(\bm x),\quad\bm m(\bm x,0)=\bm m_0(\bm x),\quad E(\bm x,0)=E_0(\bm x),\qquad\bm x\in \Td,
\label{2.1}
\end{equation}
and the periodic boundary conditions. DW solution satisfies the Euler equations in the weak sense modulo the Reynolds stress $\mathfrak R$
and energy $\mathfrak E$ defects, which are positive Radon measures, that is,
$$
\begin{aligned}
&\rho_t+\bnabla\!\cdot\!\vm=0,\\
&\bm m_t+\bnabla\!\cdot\!\left(\frac{\bm m\otimes\bm m}{\rho}\right)+\nabla p(\rho,S)+\bnabla\!\cdot\!\mathfrak R=0,\\
&S_t+\bnabla\!\cdot\!\left(S\frac{\bm m}{\rho}\right)\ge0,\\
&\int\limits_{\Td}E_0(\bm x)\,{\rm d}\bm x\ge\int\limits_{\Td}E\big(\rho(\bm x,t),\bm m(\bm x,t),S(\bm x,t)\big){\rm d}\bm x+
\int\limits_{\Td}{\rm d}\big(\mathfrak E(\bm x,t)\big)\mbox{~a.a. }t\in(0,T),
\end{aligned}
$$
where the trace of the Reynolds stress and energy defects satisfy the following inequality:
\begin{equation}
d_1\mathfrak E\le{\rm tr}\,\mathfrak R\le d_2\mathfrak E,\quad d_1=\min\big\{2,d(\gamma-1)\big\},~~d_2=\max\big\{2,d(\gamma-1)\big\}.
\label{2.2f}
\end{equation}
A rigorous definition of DW solutions is provided in Appendix \ref{appA}.

\subsubsection{Basic Properties of DW Solutions}
Unlike the weak entropy solutions, the DW solutions exist globally in time; see \cite{book}. Moreover, let us consider a sequence of
approximate solutions $\{(\rho_m,\bm m_m,S_m)\}_{m=1}^\infty$, which is consistent in the sense that each $(\rho_m,\bm m_m,S_m)$ satisfies
the weak formulation of \eref{1.1}--\eref{1.3}, \eref{2.1} with local consistency errors as $m\to\infty$, and stable, that is, uniformly
bounded with respect to $m$. It was shown in \cite{book}, such a sequence converges weakly to a DW solution.

In addition, in \cite{Feireisl2021}, the following theorem establishing a strong convergence of the so-called Ces\`aro averages,
\begin{equation}
\langle\rho\rangle_M:=\frac{1}{M}\sum_{m=1}^M\rho_m,\quad\langle\bm m\rangle_M:=\frac{1}{M}\sum_{m=1}^M\bm m_m,\quad
\langle S\rangle_M:=\frac{1}{M}\sum_{m=1}^MS_m,
\label{2.8f}
\end{equation}
has been proved.

\begin{thm}[{\bf ${\cal K}$-convergence}]\label{thm21}
Let the initial data $\{(\rho_{0,m},\bm m_{0,m},E_{0,m})\}_{m=1}^\infty$ satisfy
$$
\rho_{0,m}\ge\underline{\rho}>0,\quad E_{0,m}-\frac{ |\bm m_{0,m}|^2}{2\rho_{0,m}}>0,\quad m=1,2,\dots,
$$
where $\underline{\rho}$ is a constant independent of $m$, and let $\{(\rho_m,\bm m_m,S_m)\}_{m=1}^\infty$ be a consistent approximate
solution of \eref{1.1}--\eref{1.3}. Further, let
$$
\rho_m(\bm x,t)\ge\underline{\rho}>0,\quad\mbox{and}\quad E_m(\bm x,t)\le\,\xbar{E},
$$
where $\xbar{E}$ is another constant independent of $m$.

Then, the sequence $\{(\rho_m,\bm m_m,S_m)\}_{m=1}^\infty$ is uniformly bounded and there exists its subsequence
$(\rho_{m_n},\bm m_{m_n},S_{m_n})$ that converges strongly to a DW solution $(\rho,\bm m,S)$ in the following sense.

\smallskip
\noindent
{\bf (i) Strong convergences of Ces\` aro averages:}
\begin{equation*}
\frac{1}{M}\sum_{n=1}^M(\rho_{m_n},\bm m_{m_n},S_{m_n})\to(\rho,\bm m,S),\quad
\frac{1}{M}\sum_{n=1}^ME(\rho_{m_n},\bm m_{m_n},S_{m_n})\to\big<{\cal V}_{\bm x,t},E(\tilde\rho,\tilde{\bm m},\tilde S)\big>
\end{equation*}
as $M\to\infty$ in $L^q\big(\Td\times(0,T);\R^{d+2}\big)$ for any $0\le q<\infty$. Here ${\cal V}_{\bm x,t}$ is a space-time parametrized
probability measure on $\R^{d+2}$ and $(\rho,\bm m,S)$ are the mean values with respect to the Young measure ${\cal V}$.

\smallskip
\noindent
{\bf (ii) Strong convergence to the Young measure in the Wasserstein metric:}
\begin{equation*}
W_r\bigg[\frac{1}{M}\sum_{n=1}^M\delta_{[\rho_{m_n},\bm m_{m_n},S_{m_n}]};{\cal V}_{\bm x,t}\bigg]\to0
\end{equation*}
as $M\to\infty$ in $L^{\tilde r}(\Td\times(0,T))$ for any $1\le\tilde r<r<\infty$. Here $W_r$ denotes the Wasserstein metric of order $r$
\footnote{The Wasserstein distance of $q$-th order of probability measures ${\cal N}$ and ${\cal V}$ is defined as
$W_q({\cal N},{\cal V}):=\Big\{\inf_{\pi\in\Pi({\cal N},{\cal V})}
\int_{\mathbb R^{d+3}\times\mathbb R^{d+3}}|\bm{\zeta}_1-\bm{\zeta}_2|^q\,{\rm d}\pi(\bm{\zeta}_1,\bm{\zeta}_2)\Big\}^{1/q}$,
$q\in[1,\infty)$, where $\Pi({\cal N},{\cal V})$ is the set of probability measures on $\mathbb R^{d+3}\times\mathbb R^{d+3}$ with marginals
${\cal N}$ and ${\cal V}$.}.
\end{thm}

Moreover, DW solutions satisfy the following properties:

\smallskip
\noindent
$\bullet$ {\bf Weak-strong uniqueness}

\noindent
If a strong solution to the compressible Euler equations \eref{1.1}--\eref{1.3} exists, then any DW solution emanated from the same initial
data coincides with the strong solution on its lifespan; see \cite{brezina};

\smallskip
\noindent
$\bullet$ {\bf Compatibility}

\noindent
If a DW solution $(\rho,\bm u,S)\in C^1\big(\Td\times[0,T];\mathbb R^{d+2}\big)$, $\inf_{\Td\times(0,T)}\rho>0$, $\bm u=\bm m/\rho$, then
$(\rho,\bm m,S)$ is a classical solution of \eref{1.1}--\eref{1.3}; see \cite{book}. Specifically,
\begin{equation*}
\mathfrak R=0,\quad{\cal V}_{\bm x,t}=\delta_{[\rho(\bm x,t),\bm m(\bm x,t),S(\bm x,t)]}~\mbox{ for }(\bm x,t)\in\Td\times(0,T).
\end{equation*}

\subsection{Random Solutions of the Compressible Euler Equations}\label{sec22}
Random DW solutions can be defined as in Definition \ref{def21} with \eref{2.2}--\eref{2.7} hold $\mathbb P$-{\it a.s.} in $\Omega$, and 
the following theorem provides the convergence results for the numerical solutions obtained by the Monte Carlo method. Its proof can be
established analogously to the proof of \cite[Theorem 5.6]{FLMY}, where the barotropic Euler system was considered.
\begin{thm}[{\bf Convergence of the Monte Carlo method}]\label{thm22}
Suppose the initial data,
\begin{equation}
\rho(\bm x,0;\bm\xi)=\rho_0(\bm x;\bm\xi),~\bm m(\bm x,0;\bm\xi)=\bm m_0(\bm x;\bm\xi),~S(\bm x,0;\bm\xi)=S_0(\bm x;\bm\xi),~
\bm x\in\Td,\,\bm\xi\in\Omega,
\label{2.8}
\end{equation}
are measurable for each $\bm\xi$ and satisfy the following bounds: $\int_{\Td}E(\rho_0,\bm m_0,S_0)\,{\rm d}\bm x<\infty$ and
\begin{equation*}
\frac{1}{C}\le\rho_0(\bm x;\bm\xi)\le C,\quad|\bm m_0(\bm x;\bm\xi)|\le C,\quad\frac{1}{C}\le S_0(\bm x;\bm\xi)\le C~
\mbox{~for a.a.}~\bm x\in\Td,~\mathbb P\mbox{-{\it a.s.}}
\end{equation*}
for some constant $C>0$, and $\{(\rho_{0,\ell},\bm m_{0,\ell},S_{0,\ell})\}_{\ell=1}^\infty$ with
$\rho_{0,\ell}:=\rho_0(\bm x;\bm\xi_\ell)$, $\bm m_{0,\ell}:=\bm m_0(\bm x;\bm\xi_\ell)$, and $S_{0,\ell}:=S_0(\bm x;\bm\xi_\ell)$ are 
their pairwise independent identically distributed representations. Let $\{(\rho_{m,\ell},\bm m_{m,\ell},S_{m,\ell})\}_{m=1}^\infty$ with
$\rho_{m,\ell}:=\rho_m(\bm x,t;\bm\xi_\ell)$, $\bm m_{m,\ell}:=\bm m_m(\bm x,t;\bm\xi_\ell)$, and $S_{m,\ell}:=S_m(\bm x,t;\bm\xi_\ell)$ be
a consistent and stable approximation of \eref{1.1}--\eref{1.3} for each $\ell=1,\dots,\infty$.

Then there is a subsequence such that
\begin{equation*}
\mathbb E\bigg[\bigg\|\frac{1}{ML}\sum_{n=1}^M\sum_{\ell=1}^L\big(\rho_{m_n,\ell},\bm m_{m_n,\ell},S_{m_n,\ell}\big)
-\mathbb E\big[(\rho,\bm m,S)\big]
\bigg\|_{L^q(\Td\times(0,T);\mathbb R^{d+1})}\bigg]\to0\mbox{ as } L,M\to\infty
\end{equation*}
for any $1<q\le\frac{2\gamma}{\gamma+1}$, where $(\rho,\bm m,S)$ is a DW solution of the initial value problem \eref{1.1}--\eref{1.3},
\eref{2.8}, and $\mathbb E$ is the expected values with respect to the PDF $\mu(\bm\xi)$.
\end{thm}

\section{Methodology}\label{sec3}
In this section, we construct consistent and stable numerical solutions of the random initial value problem \eref{1.1}--\eref{1.2},
\eref{2.8} and discuss the analysis of their properties using stochastic tools.

Without loss of generality, we consider the case of two space dimensions in $\bm x=(x,y)$ ($d=2$) and begin by introducing embedded uniform
grids $(x_{j_m},y_{k_m})$ with $x_{j_m}=j_m\dx^m$ and $y_k^m=k_m\dy^m$, $j_m,k_m=0,\dots,N_m$, where $\dx^m=\dy^m=1/N_m$,
$N_m=2^{m-1}(2^{m_0+1}-1)$, $m_0\ge0$ is a fixed integer number, and $m=1,\dots,M$.

Next, we choose the uniformly distributed collocation points $\bm\xi_\ell$, $\ell=1,\dots,L$, and numerically solve \eref{1.1}--\eref{1.2},
\eref{2.8} on the aforementioned sequence of embedded meshes using the fifth-order A-WENO scheme from \cite{Shaoshuai2022}.

We then denote the discrete solutions obtained at time level $t$ for the sequence of embedded meshes for $m=1,\dots,M$ by
\begin{equation}
\bm U_{m,\ell}(t):\approx\big\{\bm U(x_{j_m},y_{k_m},t;\bm\xi_\ell)\big\}_{j_m,k_m=0}^{N_m},
\label{3.1f}
\end{equation}
where $\bm U:=(\rho,\bm m,S)^\top$ (we stress
that we evolve in time the conservative variables $\rho$, $\bm m$, and $E$, and then recalculate $S$), and evaluate the Ces\`aro averages at
the final time $T$ using \eref{2.8f}. To this end, we recall that in \eref{2.8f}, discrete solutions computed on several embedded meshes are
to be averaged. We therefore first project all of the solutions corresponding to $m=1,\dots,M-1$ onto the finest mesh, which corresponds to
$m=M$. This is done using the one-dimensional (1-D) uniformly seventh-order accurate CWENO7 interpolation \cite{CIK} applied in a
``dimension-by-dimension'' manner (first in the $x$-direction and then in the $y$-direction). As a result, we obtain the Ces\`aro averages
$\langle\bm U\rangle_M(x_{j_M},y_{k_M},T;\bm\xi_\ell)$, for which we compute statistical quantities with respect to $\xi$, which, from now
on, will be assumed to be 1-D ($s=1$). Specifically, we compute the mean, variance, and standard deviation,
\begin{equation}
\mathbb E[\psi]:=\int\limits_\Omega\psi(\xi)\mu(\xi)\,{\rm d}\xi,\quad{\rm Var}[\psi]:=\mathbb E[\psi^2]-(\mathbb E[\psi])^2,\quad
\sigma[\psi]:=\sqrt{{\rm Var}[\psi]},
\label{3.1}
\end{equation}
for all $j_M,k_M=0,\dots,N_M$ and each of the components of $\langle\bm U\rangle_M(x_{j_M},y_{k_M},T;\xi)$, which are denoted by
$\psi(\xi)$ in \eref{3.1}.

Notice that we only have the discrete values $\psi(\xi_\ell)$ available, where the collocation points $\xi_\ell$ are uniformly distributed 
over the interval $\Omega=[a,b]$ so that $\xi_\ell=(\ell-1)\Delta\xi$, $\Delta\xi=(b-a)/(L-1)$, $\ell=1,\dots,L$. Hence, we need to use a 
proper quadrature in the integrals in \eref{3.1}. To this end, we use the CWENO7 interpolation from \cite{CIK} to obtain a piecewise 
polynomial approximation of $\psi$:
\begin{equation}
\sum_{\ell=1}^L\psi_\ell(\xi)\mbox{{\Large$\chi$}}_{[\xi_{\ell-\hf},\xi_{\ell+\hf}]}(\xi),
\label{3.2}
\end{equation}
where $\psi_\ell$ are the CWENO7 polynomial pieces described in \cite{CIK} and $\mbox{{\Large$\chi$}}_{[\xi_{\ell-\hf},\xi_{\ell+\hf}]}$ is
a characteristic function of the interval $[\xi_{\ell-\hf},\xi_{\ell+\hf}]$ with $\xi_{\ell\pm\hf}=(\xi_{\ell\pm1}+\xi_\ell)/2$. We then
substitute \eref{3.2} into \eref{3.1} to end up with the following approximations of $\mathbb E[\psi]$ and $\sigma[\psi]$:
\begin{equation}
\begin{aligned}
\xbar\psi&=\int\limits_{\xi_1}^{\xi_{\frac{3}{2}}}\psi_1(\xi)\mu(\xi)\,{\rm d}\xi+
\sum_{\ell=2}^{L-1}\int\limits_{\xi_{\ell-\hf}}^{\xi_{\ell+\hf}}\psi_\ell(\xi)\mu(\xi)\,{\rm d}\xi+
\int\limits_{\xi_{L-\hf}}^{\xi_L}\psi_L(\xi)\mu(\xi)\,{\rm d}\xi,\\
\xbar\sigma&=\Bigg(\int\limits_{\xi_1}^{\xi_{\frac{3}{2}}}(\psi_1(\xi)-\,\xbar\psi)^2\mu(\xi)\,{\rm d}\xi+
\sum_{\ell=2}^{L-1}\int\limits_{\xi_{\ell-\hf}}^{\xi_{\ell+\hf}}(\psi_\ell(\xi)-\,\xbar\psi)^2\mu(\xi)\,{\rm d}\xi\\
&\hspace*{4.9cm}+\int\limits_{\xi_{L-\hf}}^{\xi_L}(\psi_L(\xi)-\,\xbar\psi)^2\mu(\xi)\,{\rm d}\xi\Bigg)^\hf,
\end{aligned}
\label{3.3}
\end{equation}
which can be evaluated either exactly or with high accuracy using a proper Gaussian quadrature.

As mentioned in the Introduction, the solution $\bm U$ is not expected to be unique but can be characterized by a family of $(x,y,t;\xi)$
parameterized Young measure ${\cal V}_{x,y,t}$, which we approximate using the obtained mean of the Ces\`aro averages
$\xbar{\langle\bm U\rangle}_M(x_{j_M},y_{k_M},T)$. More precisely, we fix a small spatial window $D\subset\mathbb R^2$, where histograms for
$\rho$, $\rho u$, $\rho v$, and $E$ (in the 2-D case, $\bm m=\rho\bm u$, $\bm u=(u,v)$), as well as for other quantities of
interest such as the total entropy $S$ are computed using the data of $\,\xbar{\langle\bm U\rangle}_M(x_{j_M},y_{k_M},T)$ for
$(x_{j_M},y_{k_M})\in D$.

Other quantities of interest for performing the analysis of turbulent statistics are the mean of the trace of the Reynolds stress defect
${\rm tr}(\,\xbar{\mathfrak R})(x,y)$ and the mean of the energy defect $\,\xbar{\mathfrak E}(x,y)$, which are approximated as follows:
\begin{equation}
\begin{aligned}
&\mathfrak R\approx\mathfrak R_M:=\Bigl\langle\frac{\bm m\otimes\bm m}{\rho}\Bigr\rangle_M+\langle p(\rho,S)\rangle_M\,I-
\frac{\langle\bm m\rangle_M\otimes\langle\bm m\rangle_M}{\langle\rho\rangle_M}-p\bigl(\langle\rho\rangle_M,\langle S\rangle_M\bigr)\,I,\\
&\mathfrak E\approx\mathfrak E_M=\hf\Bigl\langle\frac{|\bm m|^2}{\rho}\Bigr\rangle_M+\langle\rho\,e(\rho,S)\rangle_M-
\frac{|\langle\bm m\rangle_M|^2}{2\langle\rho\rangle_M}-\langle\rho\rangle_M\,e\bigl(\langle\rho\rangle_M,\langle S\rangle_M\bigr),
\end{aligned}
\label{3.4}
\end{equation}
where $I$ is the identity matrix. In the context of turbulence modeling, Reynolds stresses quantify the transport of momentum by 
unresolved fluctuations, while energy defects capture the mismatch between averaged and instantaneous energy balances. The computed 
quantities $\mathfrak R_M$ and $\mathfrak E_M$ thus serve as turbulence-style diagnostics, measuring the degree of
fluctuation-induced transport in our random KH flows.

\section{Numerical Study of KH Instabilities}\label{sec4}
We consider the following initial conditions, which correspond to the KH instability problem studied in \cite{Feireisl2021}:
\begin{equation}
(\rho,u,v,p)(x,y,0)=\left\{\begin{aligned}&(2,-0.5,0,2.5)&&\mbox{if~}I_1(x,y)<y<I_2(x,y),\\&(1,0.5,0,2.5)&&\text{otherwise},
\end{aligned}\right.
\label{4.1}
\end{equation}
subject to the periodic boundary conditions in the computational domain $[0,1]\times[0,1]$. The interface profiles in \eref{4.1} are given
by
\begin{equation*}
I_i(x,y)=J_i+0.05Y_i(x,y,\xi),\quad i=1,2,
\end{equation*}
where $J_1=0.25$, $J_2=0.75$, and small perturbations of the interfaces are introduced using the terms
\begin{equation}
	Y_i(x,y;\xi)=(1+\tau\tanh\xi)\sum_{k=1}^{10}a_i^k\cos(b_i^k+10k\pi x),\quad i=1,2,
\label{4.2}
\end{equation}
where $a_i^k\in[0,1]$ and $b_i^k\in[-\pi,\pi]$ ($k=1,\dots,10$) are uniformly distributed random variables. To ensure
$|I_i(x,y)-J_i|\le0.05$, the coefficients $a_i^k$ are normalized such that $\sum_{k=1}^{10}a_i^k=1$. The random numbers $a_i^k$ and $b_i^k$
are generated once for repeatability and consistency.

In space, we use the embedded grids specified in \S\ref{sec3} with $M=5$. Along the $\xi$-direction, we set $\Omega=[-1,1]$ and take
$L=101$. Note that according to \eref{4.2}, larger values of $\xi$ introduce larger initial instability amplitude, while the uncertainty 
parameter $\tau$ linearly magnifies this effect (controls the spread in $\xi$). Below, we take $\tau=1.1$ unless specified differently.

\begin{figure}[ht!]
\centerline{\small $\xi=-1$}
\centerline{\includegraphics[trim=0.0cm 0cm 2.8cm 0cm, clip, width=0.165\textwidth]{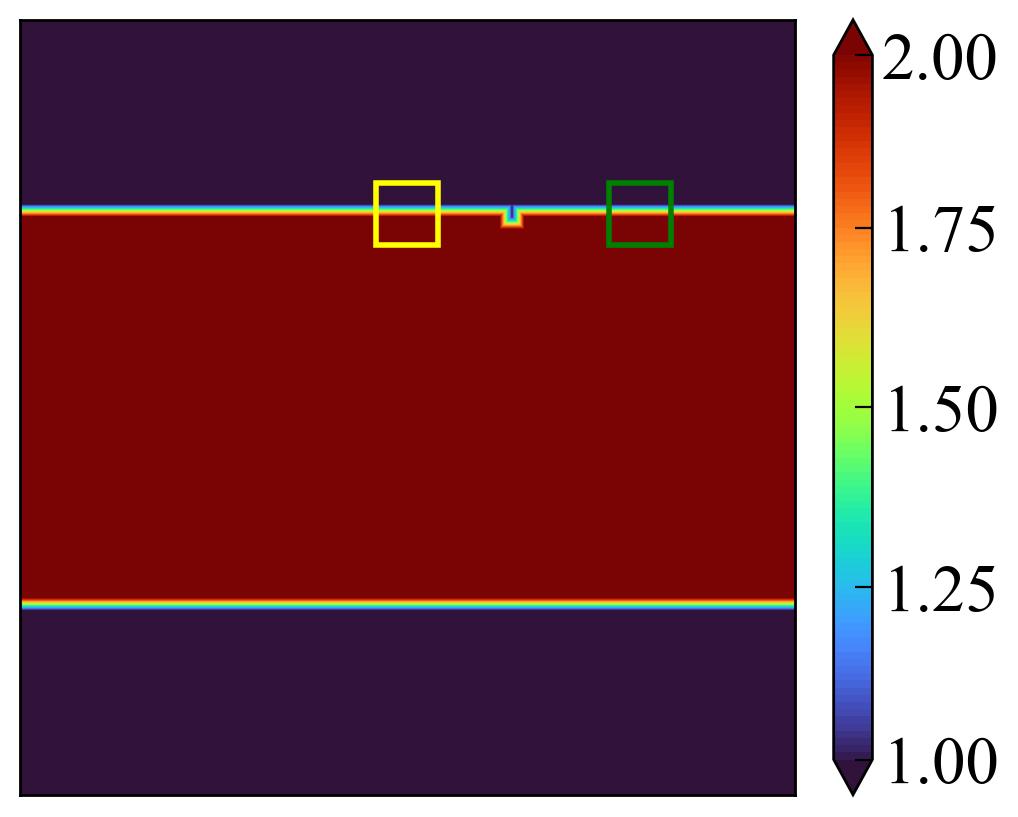}\hspace*{0.45cm}
            \includegraphics[trim=0.0cm 0cm 2.8cm 0cm, clip, width=0.165\textwidth]{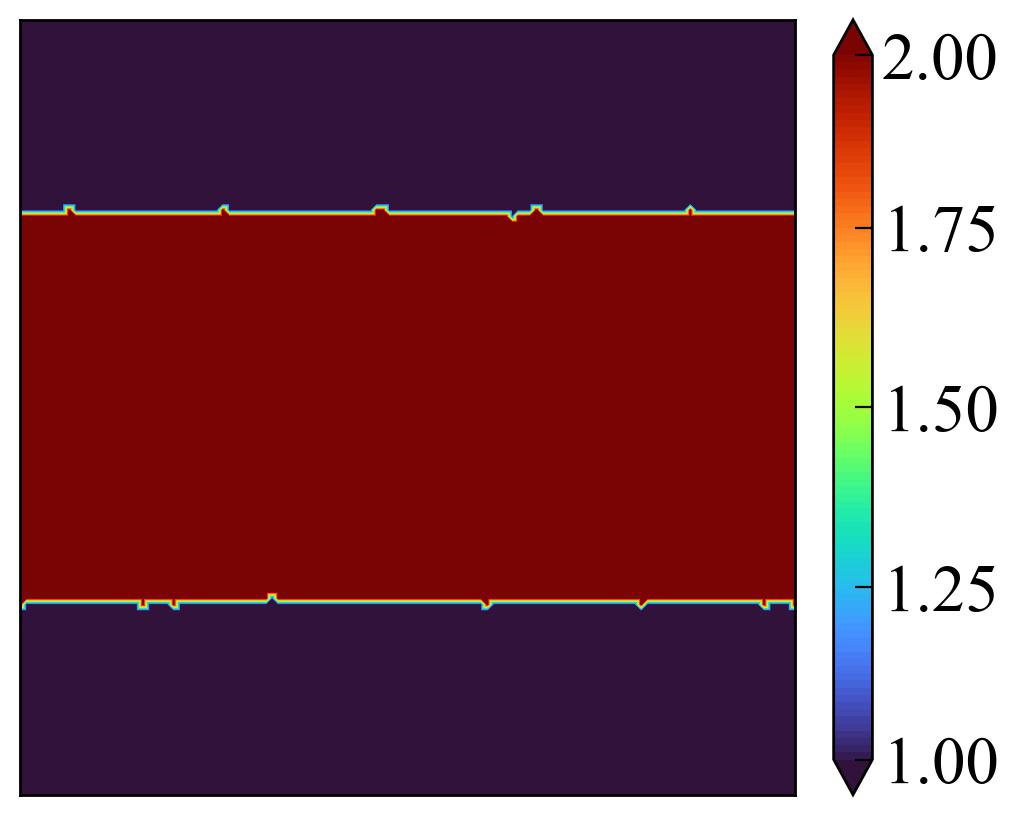}\hspace*{0.45cm}
            \includegraphics[trim=0.0cm 0cm 2.8cm 0cm, clip, width=0.165\textwidth]{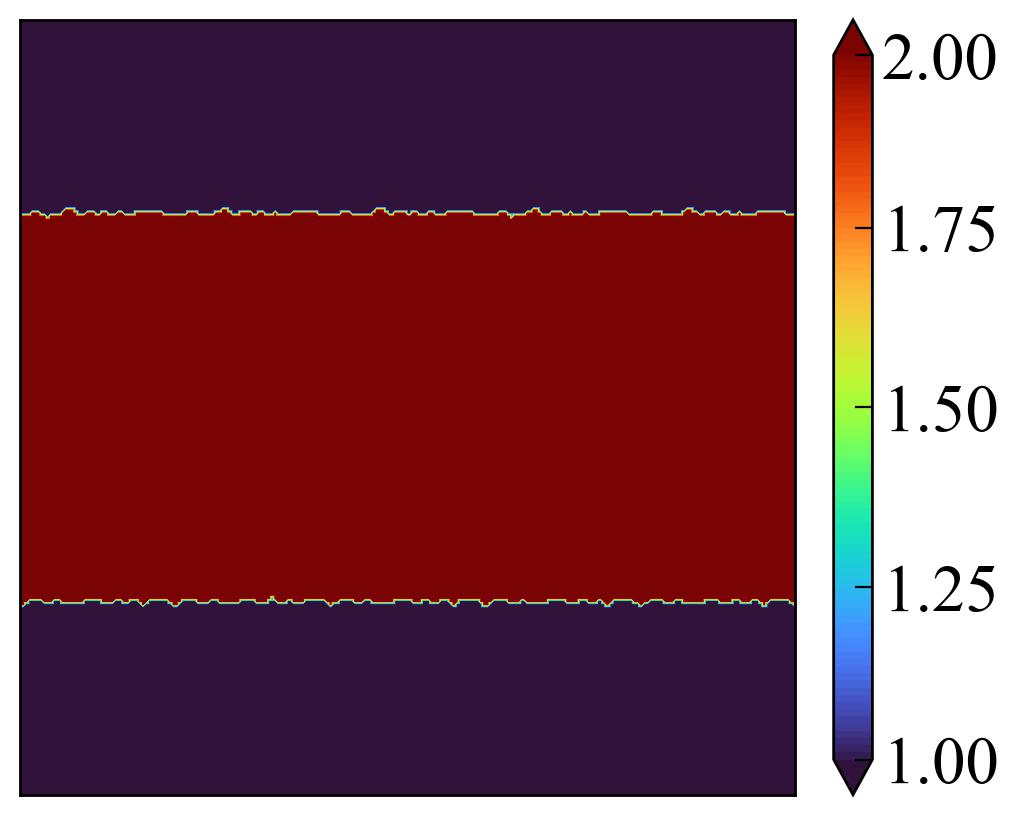}\hspace*{0.45cm}
            \includegraphics[trim=0.0cm 0cm 2.8cm 0cm, clip, width=0.165\textwidth]{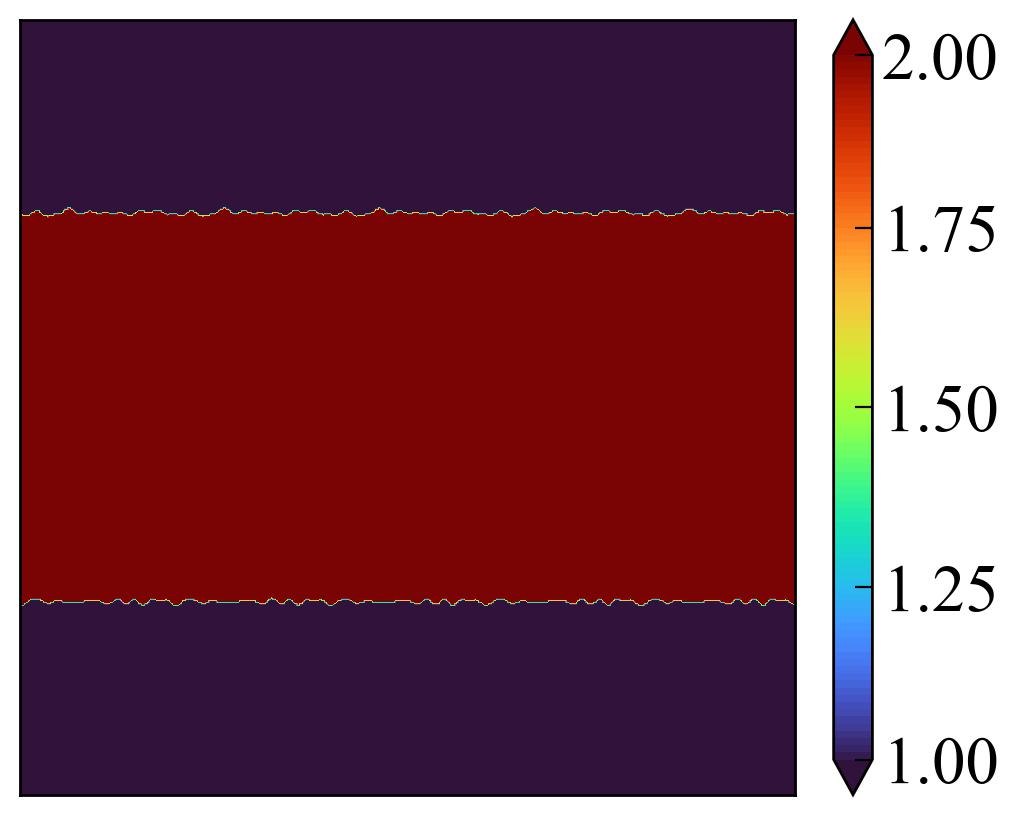}\hspace*{0.45cm}
            \includegraphics[trim=0.0cm 0cm 0.0cm 0cm, clip, width=0.21\textwidth]{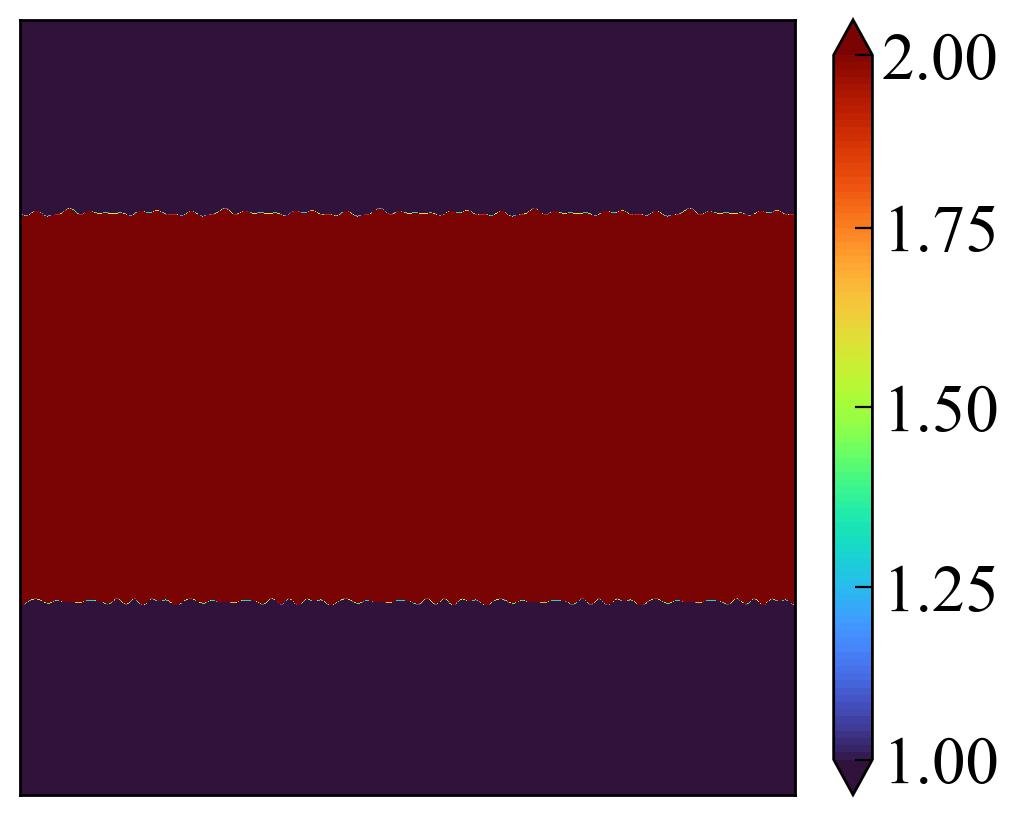}}
\vskip5pt
\centerline{\small $\xi=0$}
\centerline{\includegraphics[trim=0.0cm 0cm 2.8cm 0cm, clip, width=0.165\textwidth]{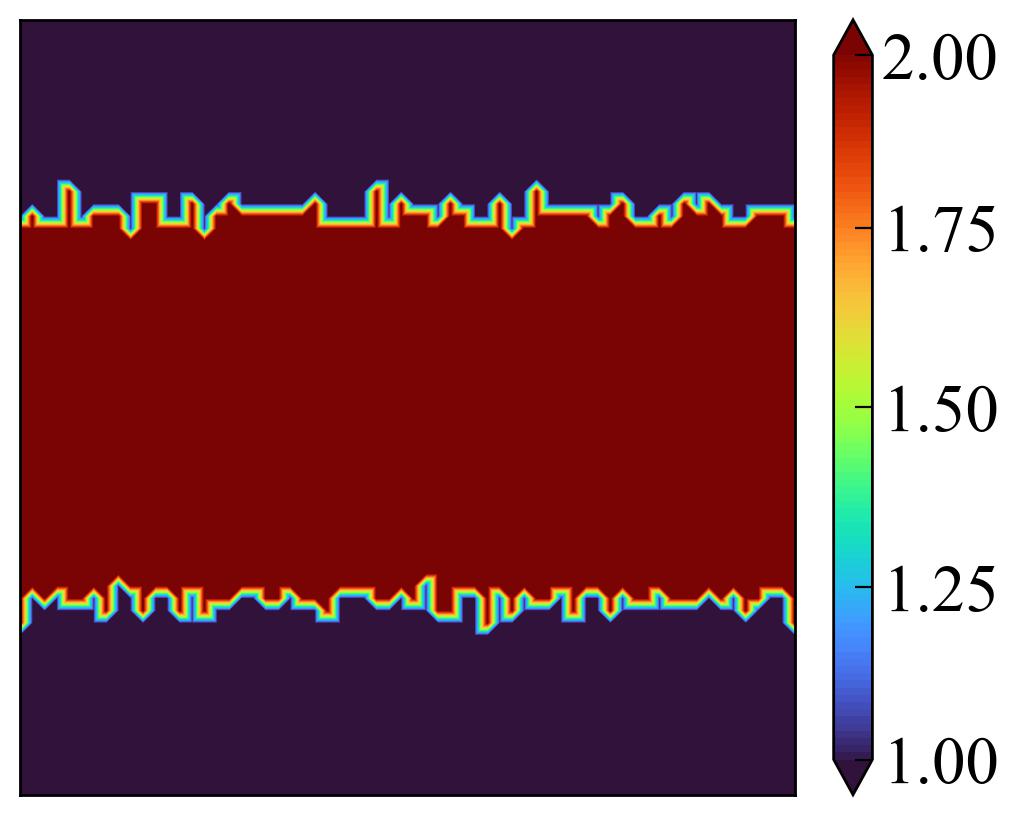}\hspace*{0.45cm}
            \includegraphics[trim=0.0cm 0cm 2.8cm 0cm, clip, width=0.165\textwidth]{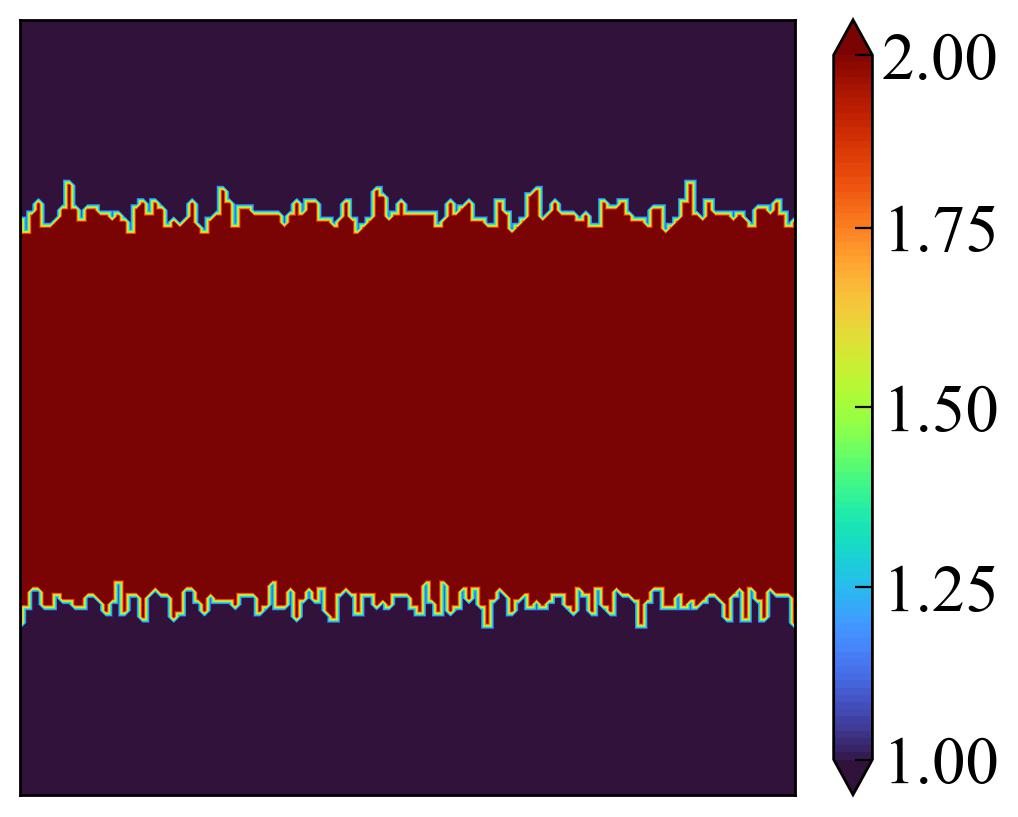}\hspace*{0.45cm}
            \includegraphics[trim=0.0cm 0cm 2.8cm 0cm, clip, width=0.165\textwidth]{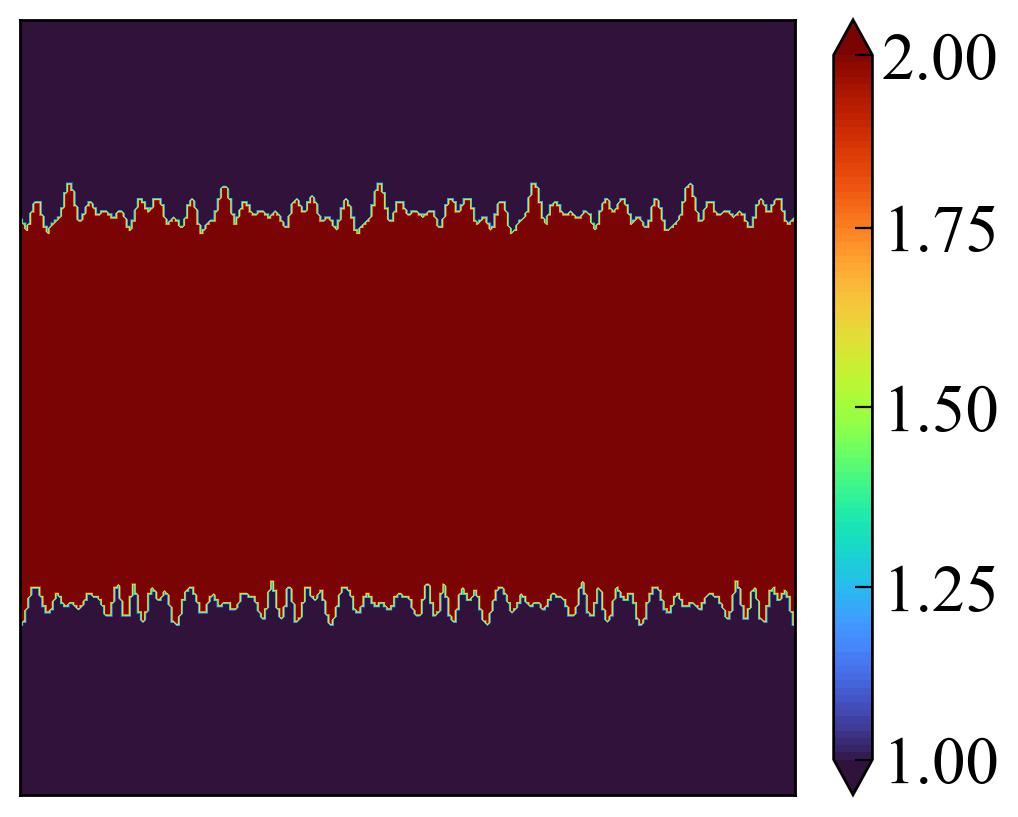}\hspace*{0.45cm}
            \includegraphics[trim=0.0cm 0cm 2.8cm 0cm, clip, width=0.165\textwidth]{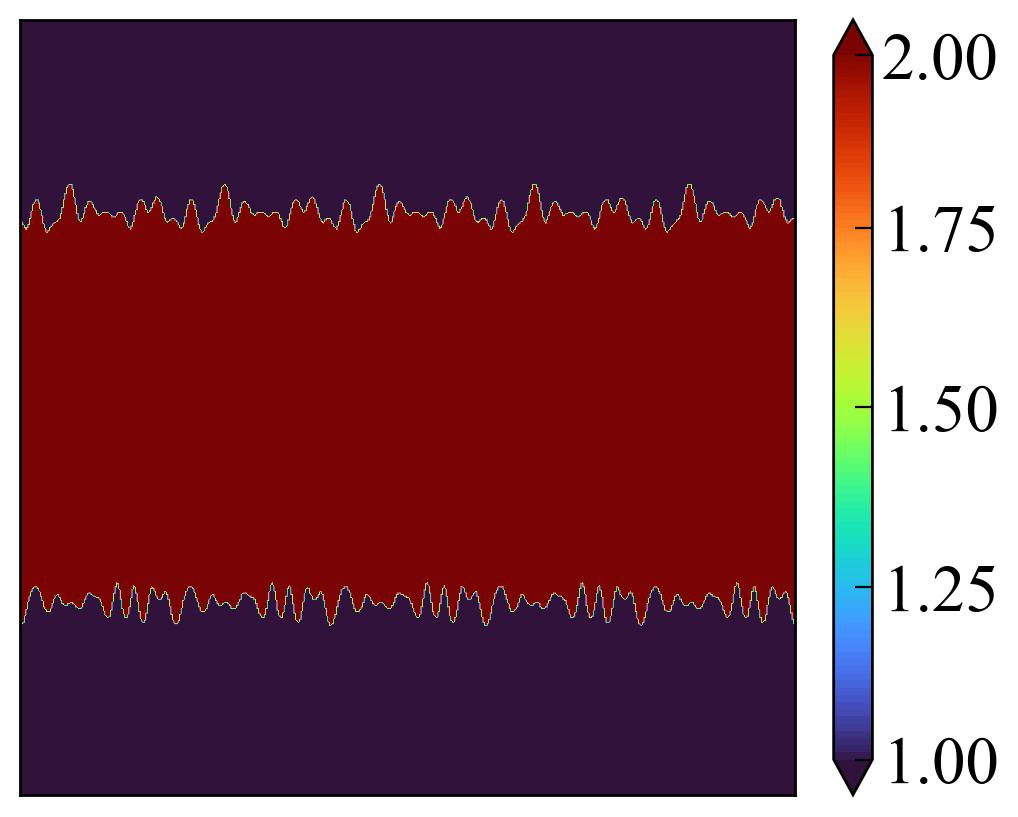}\hspace*{0.45cm}
            \includegraphics[trim=0.0cm 0cm 0.0cm 0cm, clip, width=0.21\textwidth]{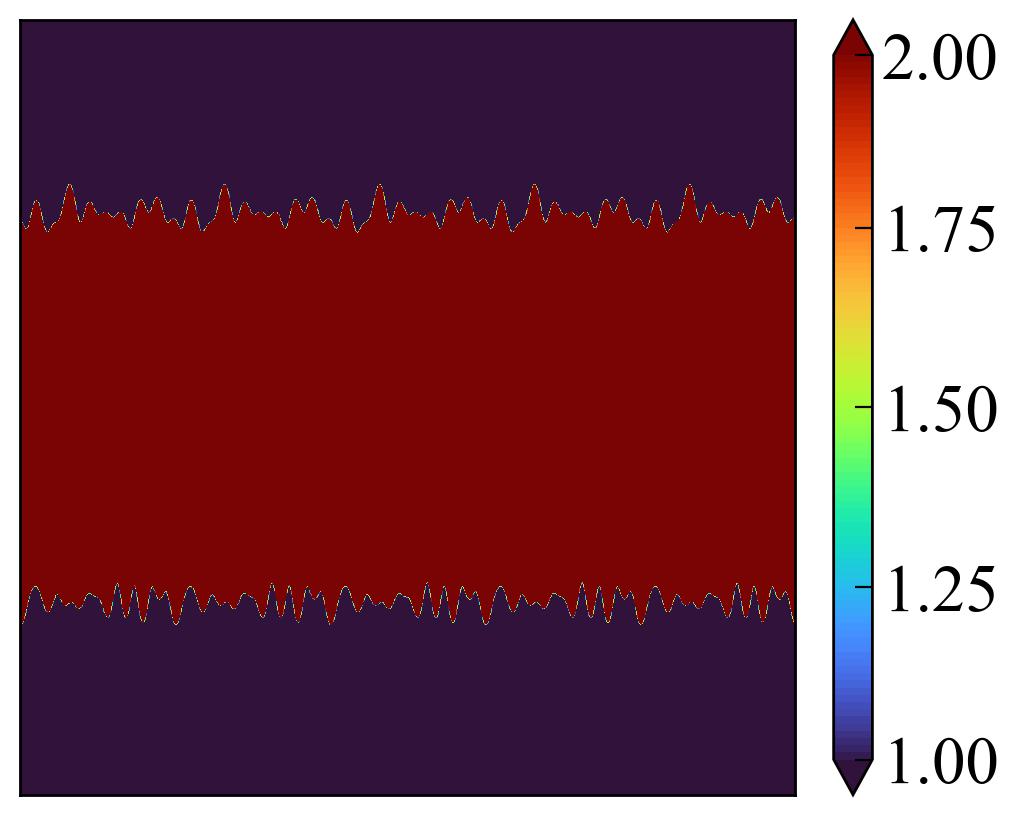}}
\vskip5pt
\centerline{\small $\xi=1$}
\centerline{\includegraphics[trim=0.0cm 0cm 2.8cm 0cm, clip, width=0.165\textwidth]{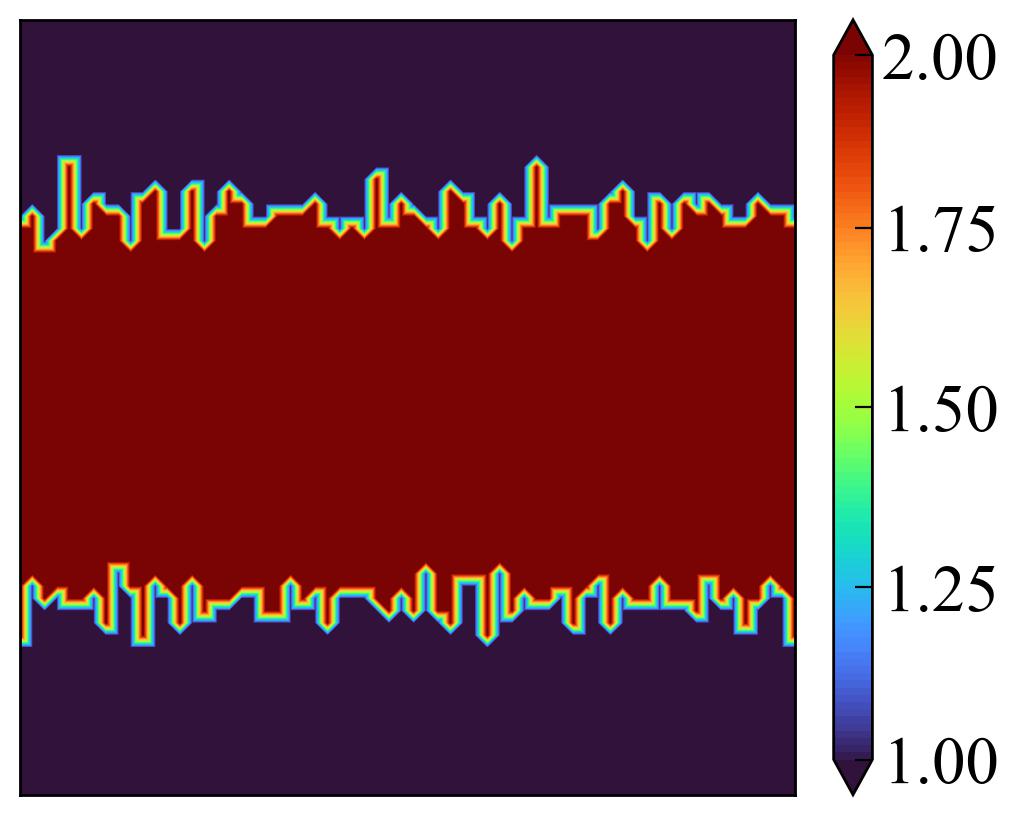}\hspace*{0.45cm}
            \includegraphics[trim=0.0cm 0cm 2.8cm 0cm, clip, width=0.165\textwidth]{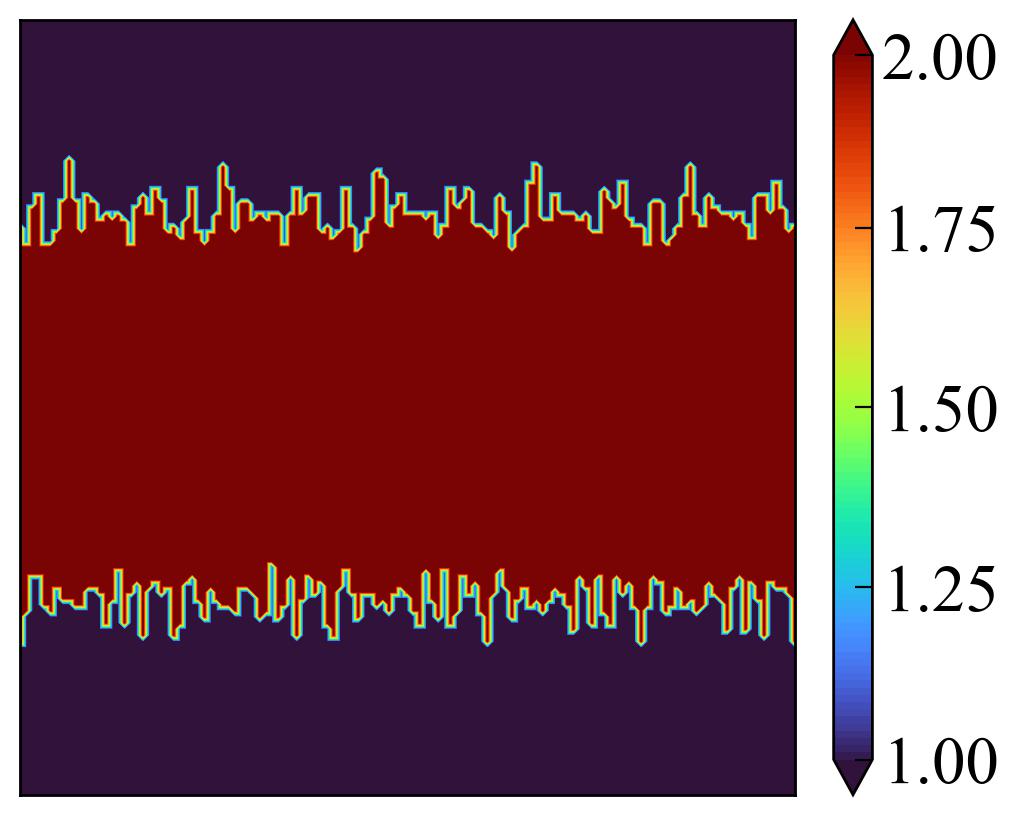}\hspace*{0.45cm}
            \includegraphics[trim=0.0cm 0cm 2.8cm 0cm, clip, width=0.165\textwidth]{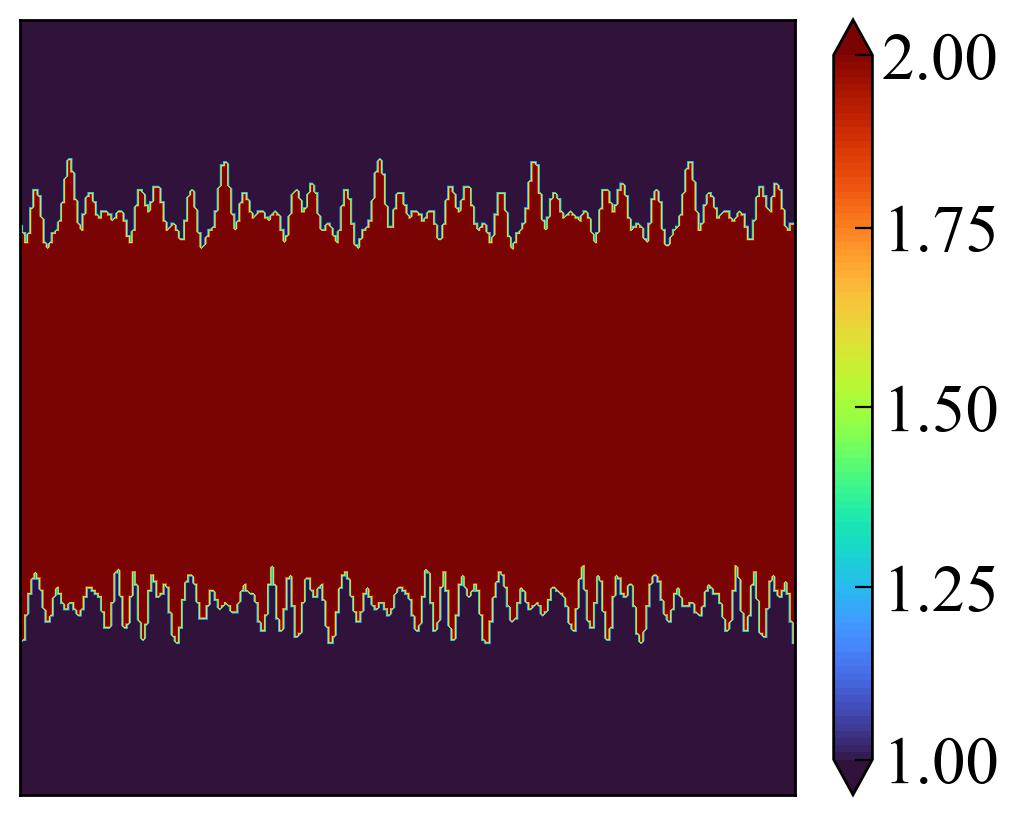}\hspace*{0.45cm}
            \includegraphics[trim=0.0cm 0cm 2.8cm 0cm, clip, width=0.165\textwidth]{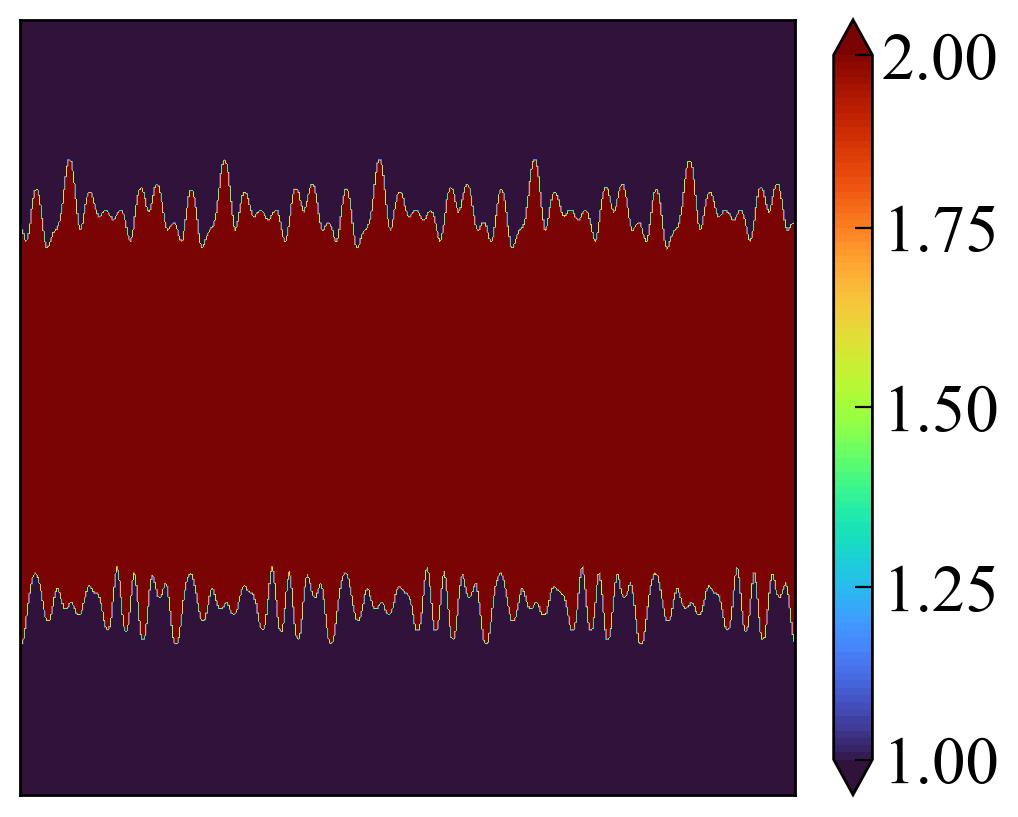}\hspace*{0.45cm}
            \includegraphics[trim=0.0cm 0cm 0.0cm 0cm, clip, width=0.21\textwidth]{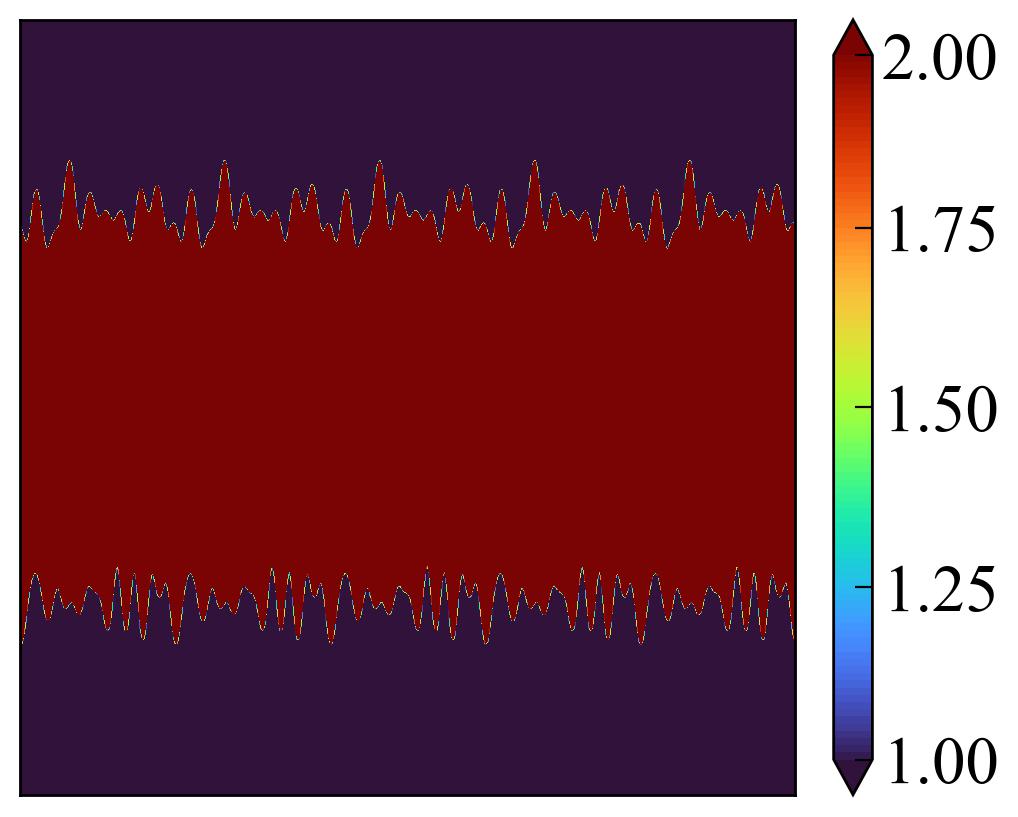}}
\caption{\sf $\rho(x,y,0;\xi)$ for $\xi=-1$ (top row), $0$ (middle row), and $1$ (bottom row), and five embedded uniform meshes with 
$m=1,\dots,5$ (from left to right). Yellow and green squares in the top left panel show the regions $D_1=[0.46,0.54]\times[0.71,0.79]$
and $D_2=[0.76,0.84]\times[0.71,0.79]$, which will be used below.\label{fig41}}
\end{figure}

\begin{figure}[ht!]
\centerline{\small $\xi=-1$}
\centerline{\includegraphics[trim=0.0cm 0cm 2.8cm 0cm, clip, width=0.165\textwidth]{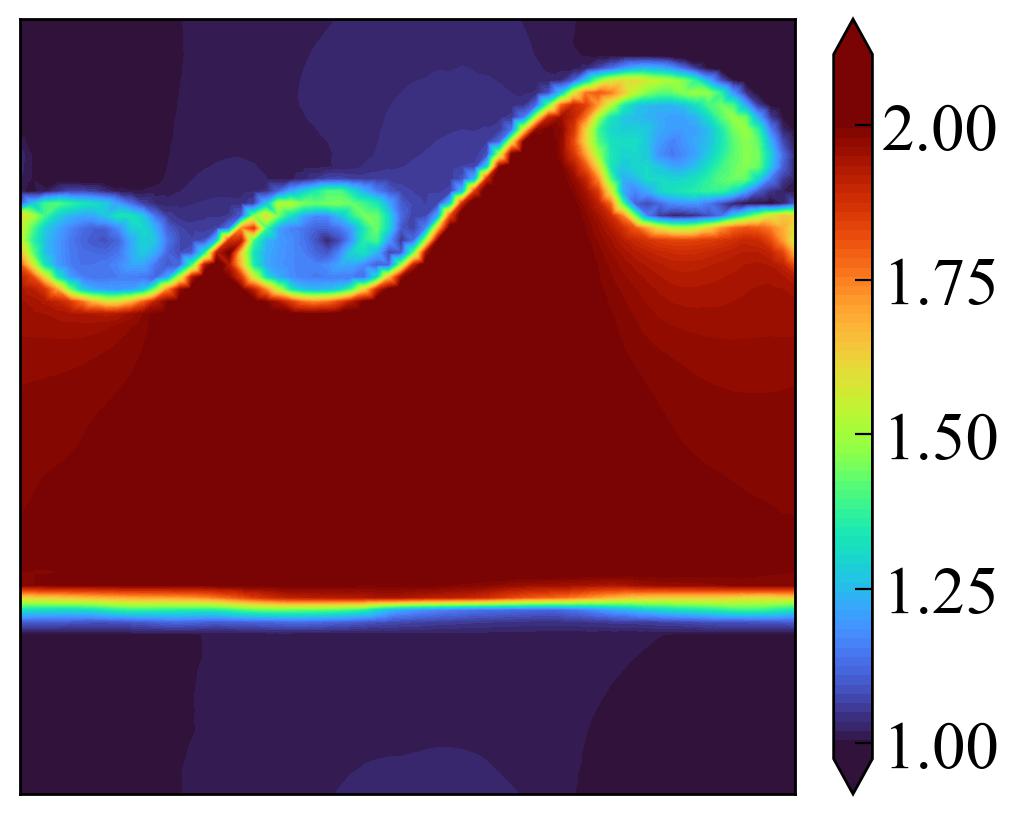}\hspace*{0.45cm}
            \includegraphics[trim=0.0cm 0cm 2.8cm 0cm, clip, width=0.165\textwidth]{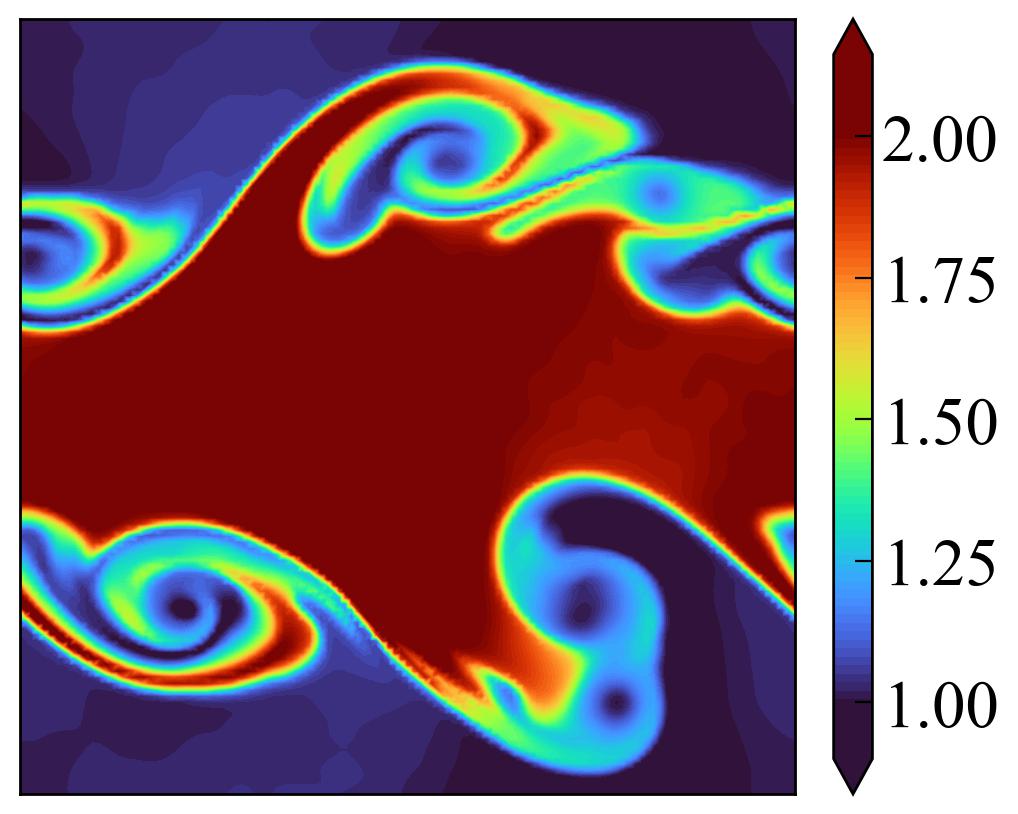}\hspace*{0.45cm}
            \includegraphics[trim=0.0cm 0cm 2.8cm 0cm, clip, width=0.165\textwidth]{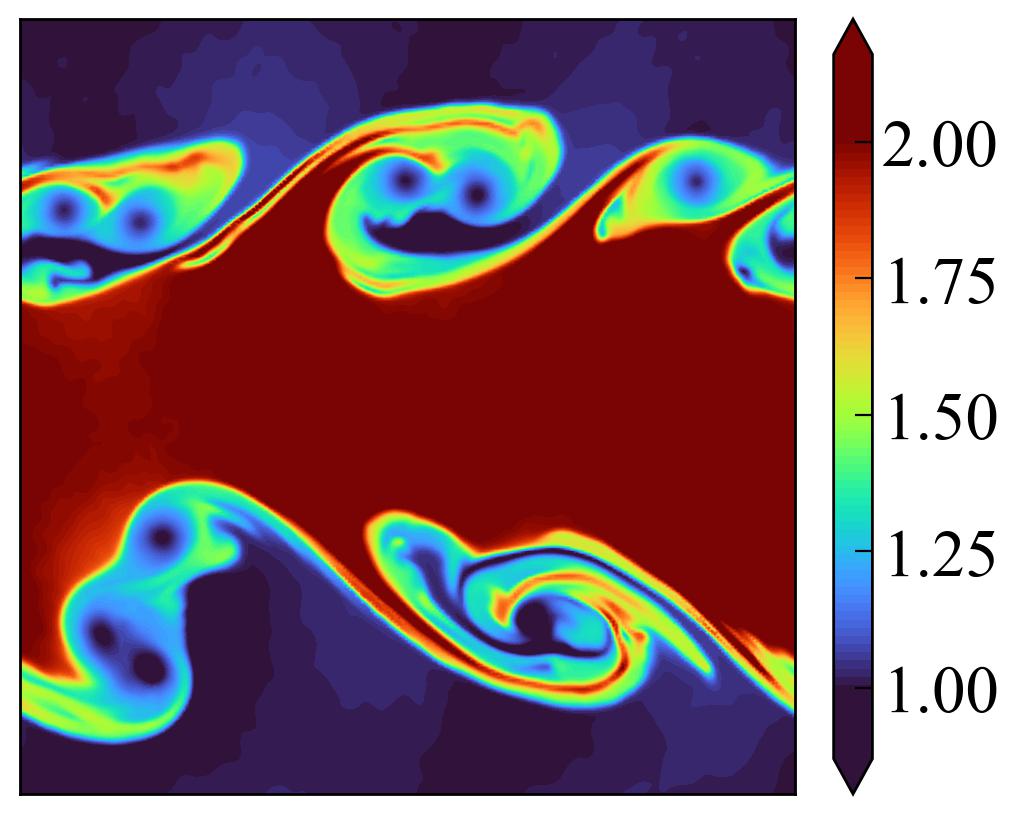}\hspace*{0.45cm}
            \includegraphics[trim=0.0cm 0cm 2.8cm 0cm, clip, width=0.165\textwidth]{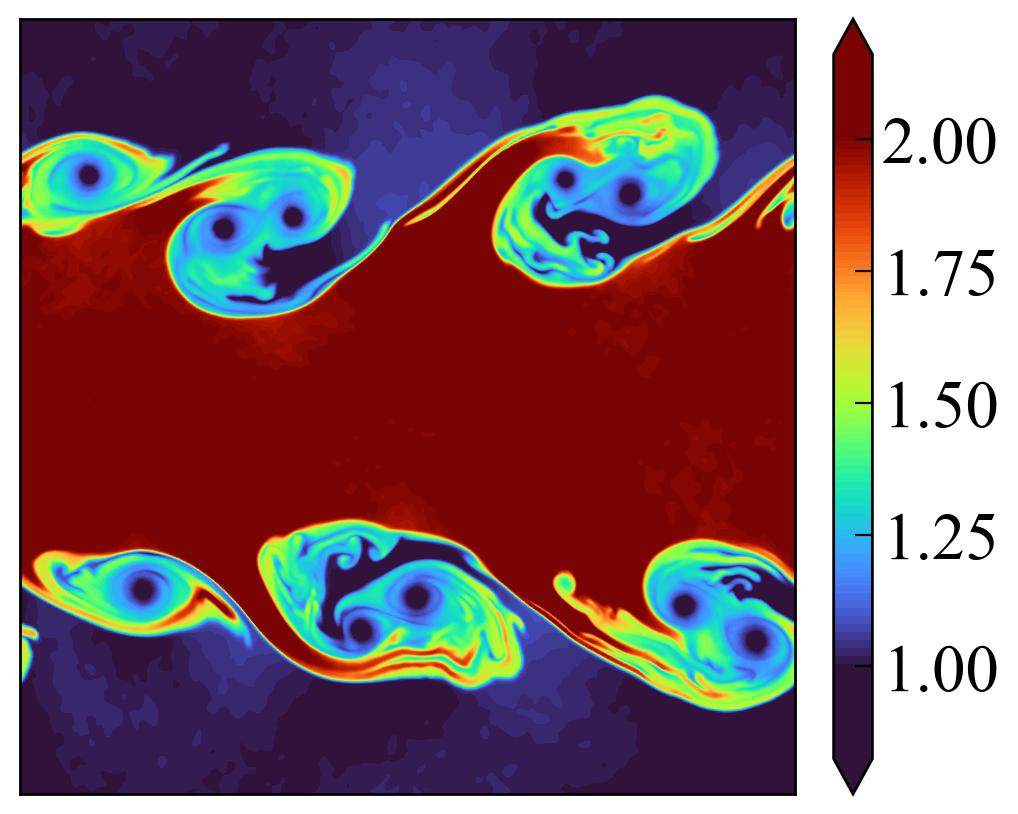}\hspace*{0.45cm}
            \includegraphics[trim=0.0cm 0cm 0.0cm 0cm, clip, width=0.21\textwidth]{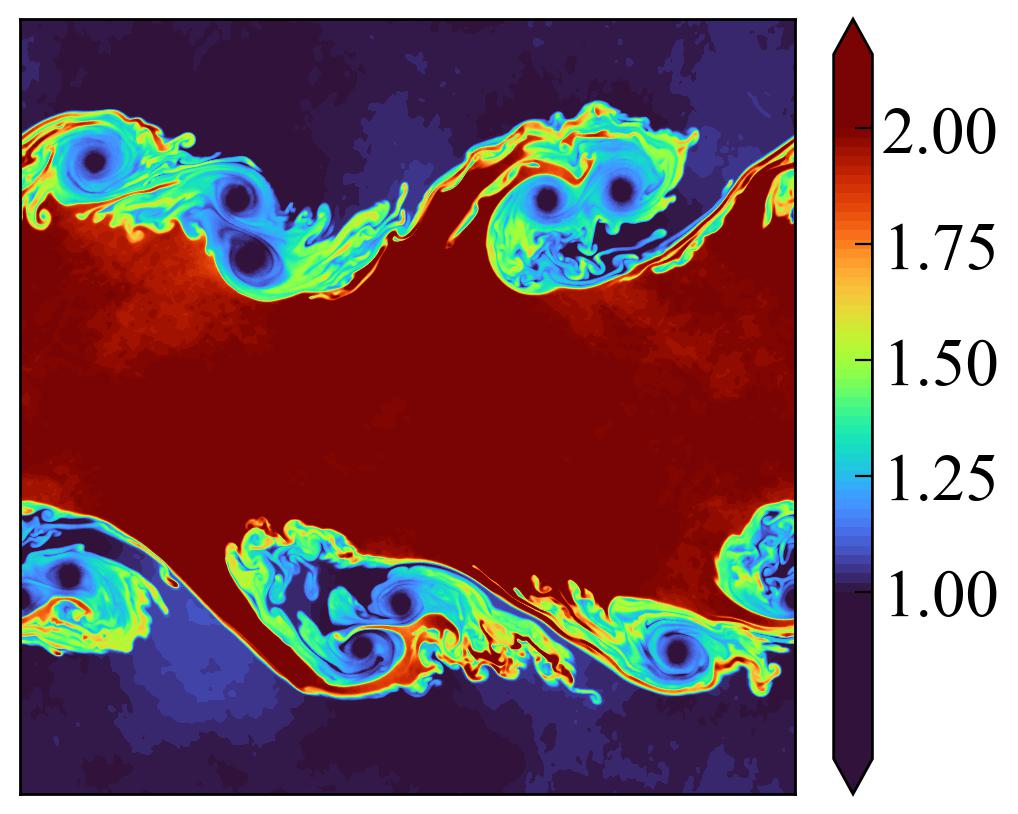}}
\vskip5pt
\centerline{\small $\xi=0$}
\centerline{\includegraphics[trim=0.0cm 0cm 2.8cm 0cm, clip, width=0.165\textwidth]{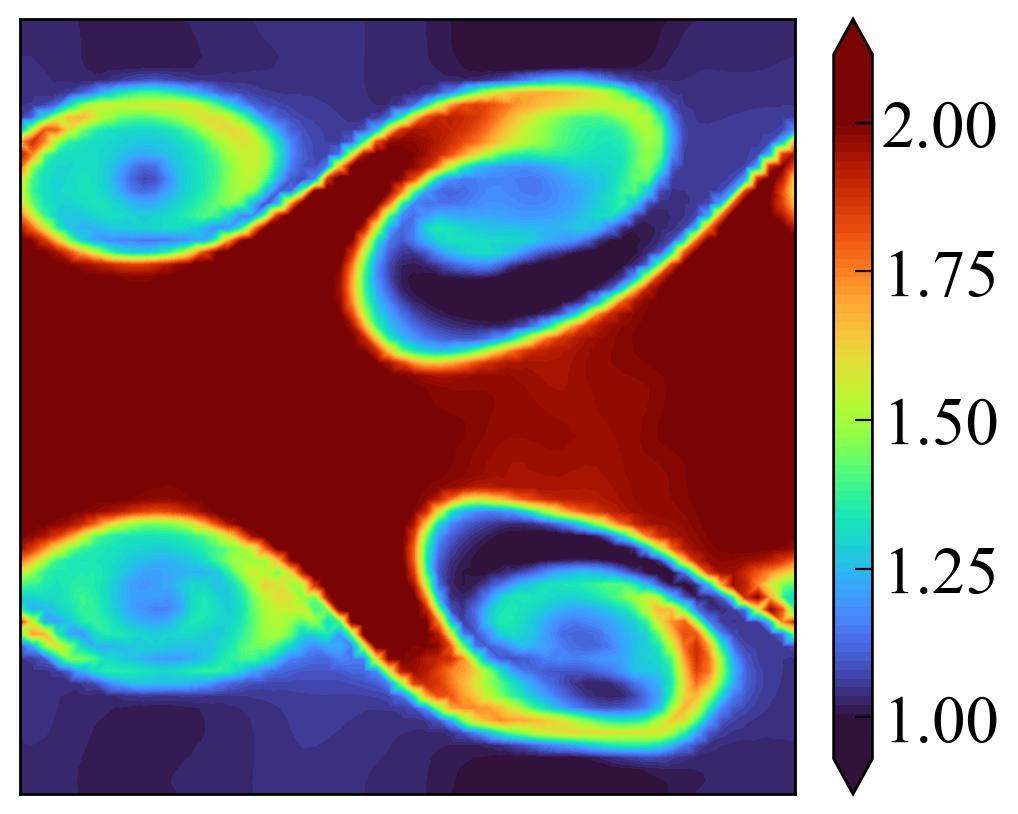}\hspace*{0.45cm}
            \includegraphics[trim=0.0cm 0cm 2.8cm 0cm, clip, width=0.165\textwidth]{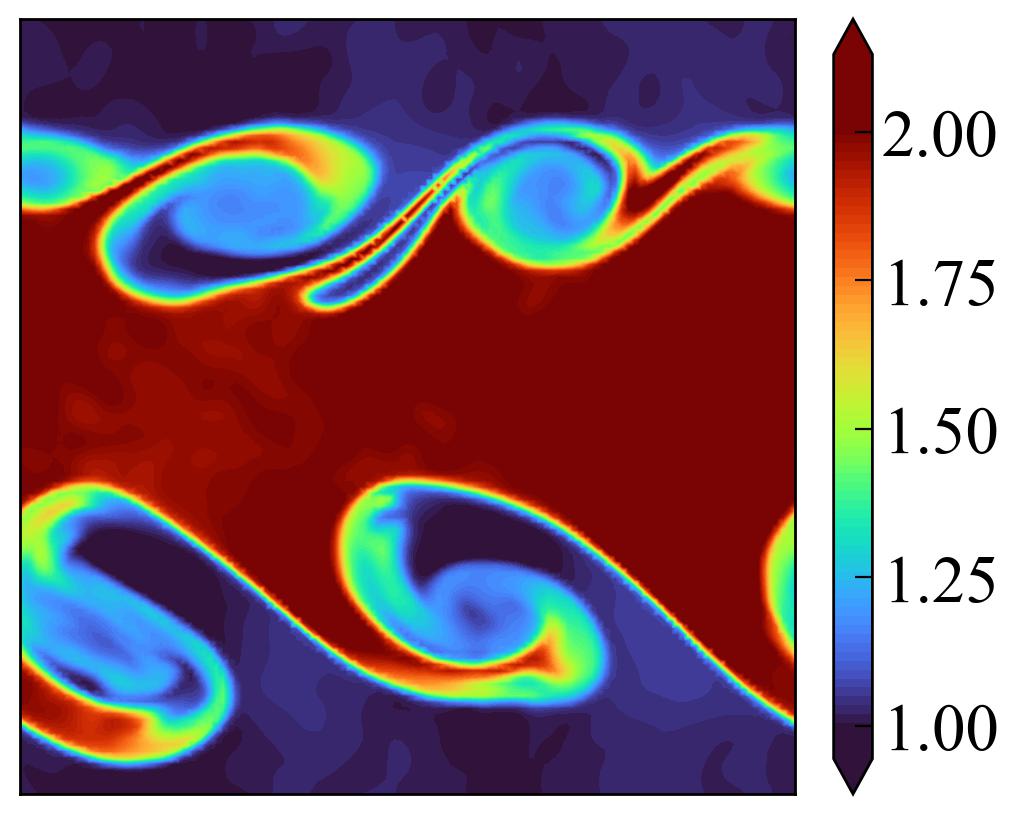}\hspace*{0.45cm}
            \includegraphics[trim=0.0cm 0cm 2.8cm 0cm, clip, width=0.165\textwidth]{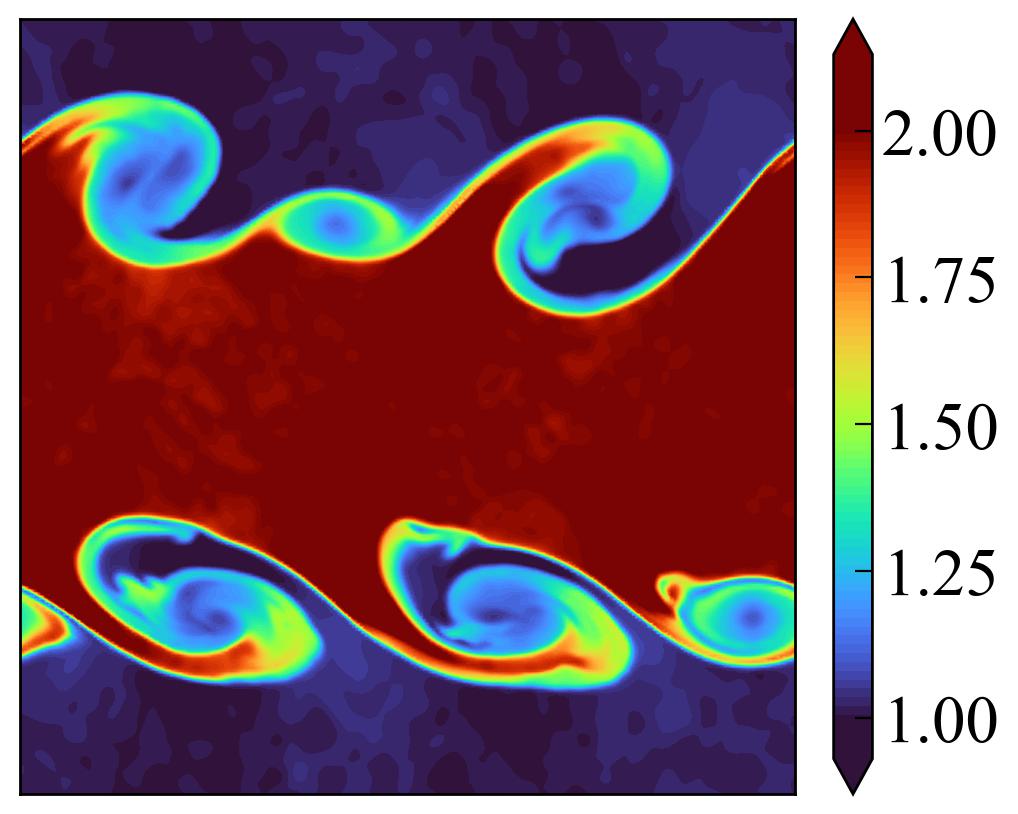}\hspace*{0.45cm}
            \includegraphics[trim=0.0cm 0cm 2.8cm 0cm, clip, width=0.165\textwidth]{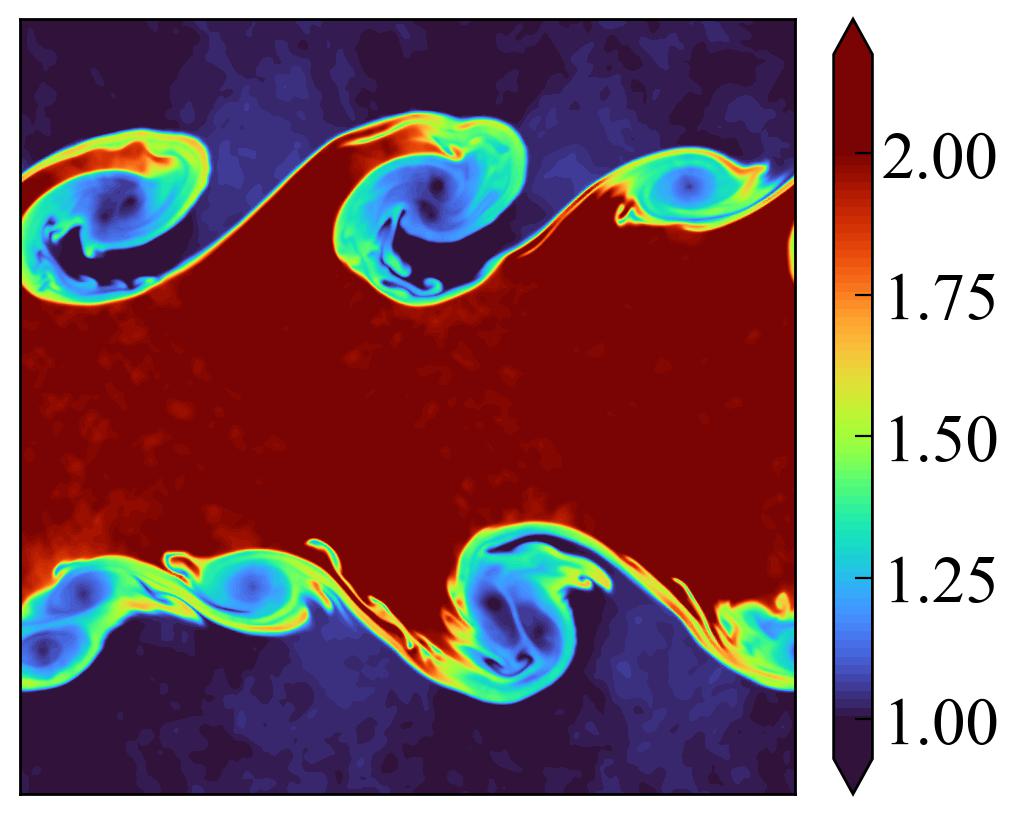}\hspace*{0.45cm}
            \includegraphics[trim=0.0cm 0cm 0.0cm 0cm, clip, width=0.21\textwidth]{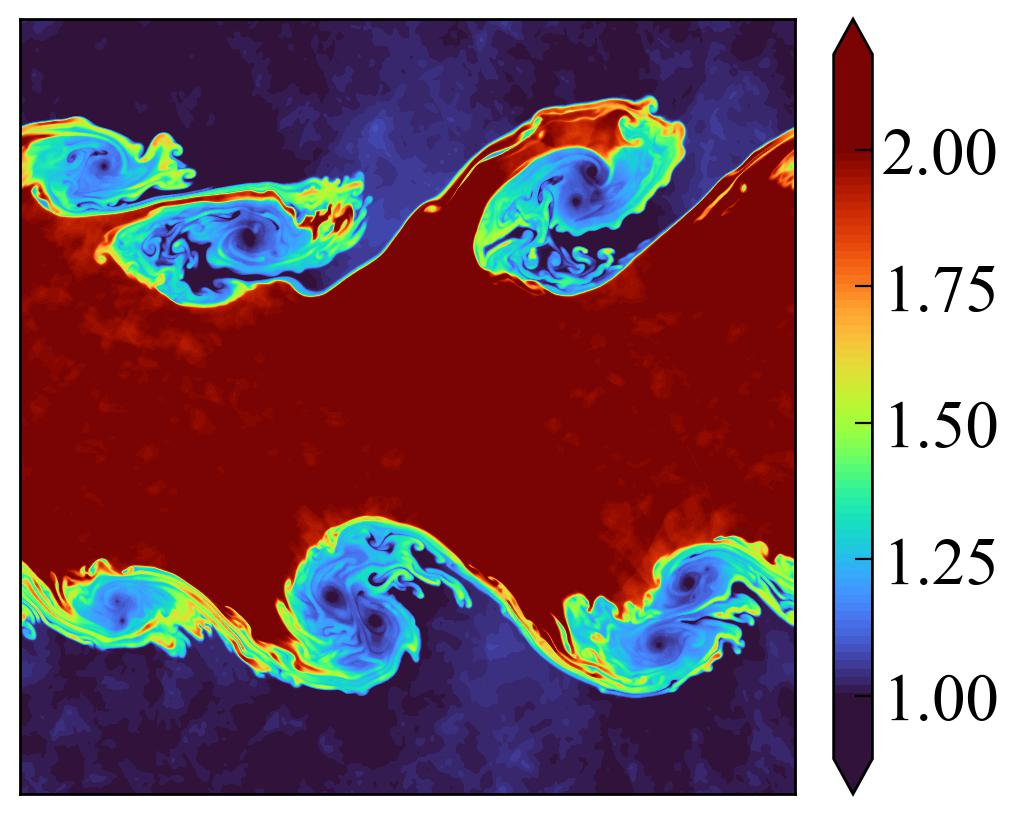}}
\vskip5pt
\centerline{\small $\xi=1$}
\centerline{\includegraphics[trim=0.0cm 0cm 2.8cm 0cm, clip, width=0.165\textwidth]{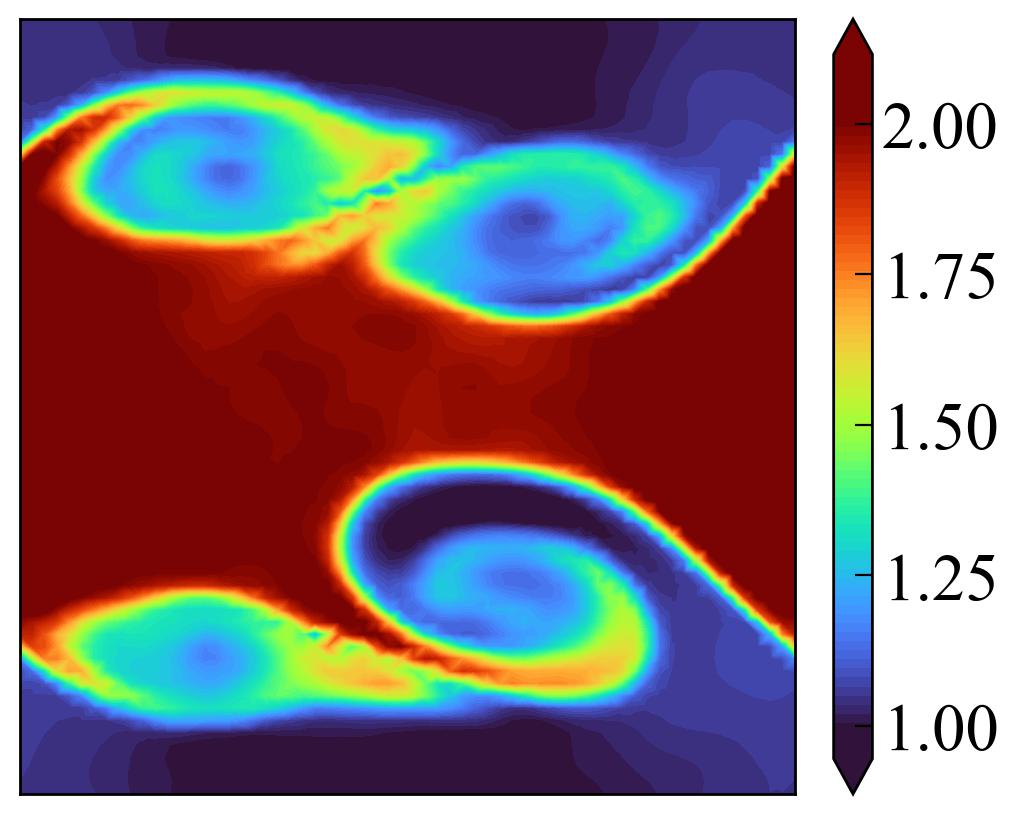}\hspace*{0.45cm}
            \includegraphics[trim=0.0cm 0cm 2.8cm 0cm, clip, width=0.165\textwidth]{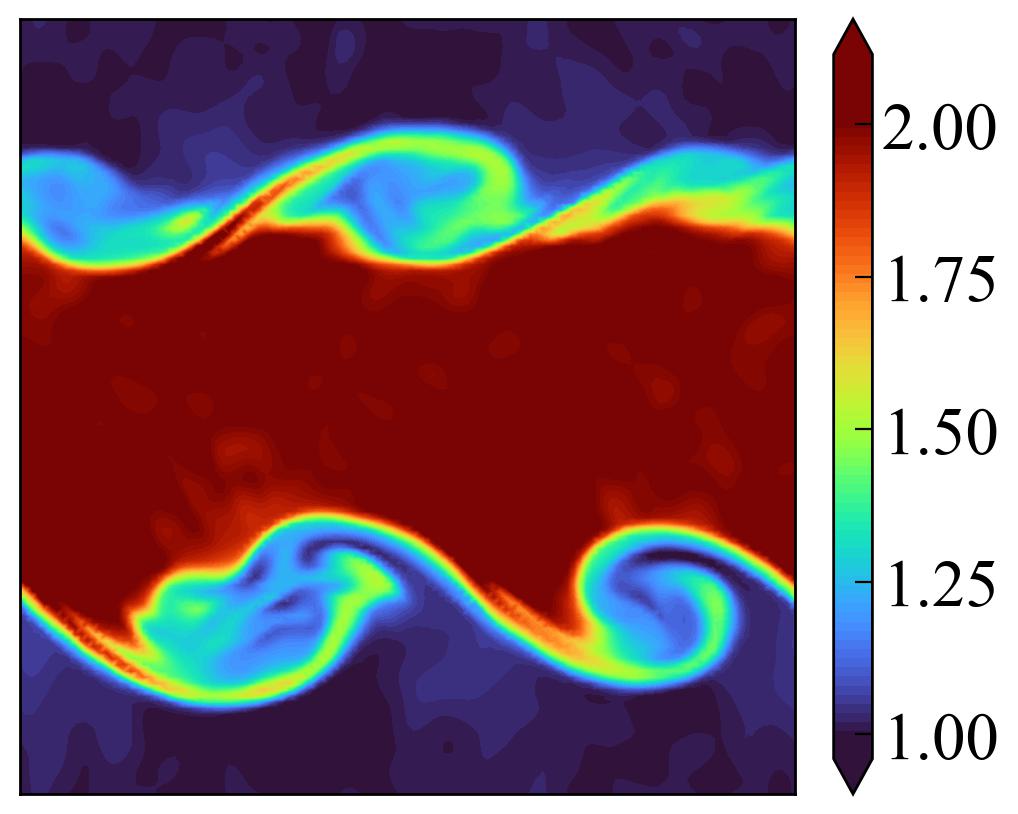}\hspace*{0.45cm}
            \includegraphics[trim=0.0cm 0cm 2.8cm 0cm, clip, width=0.165\textwidth]{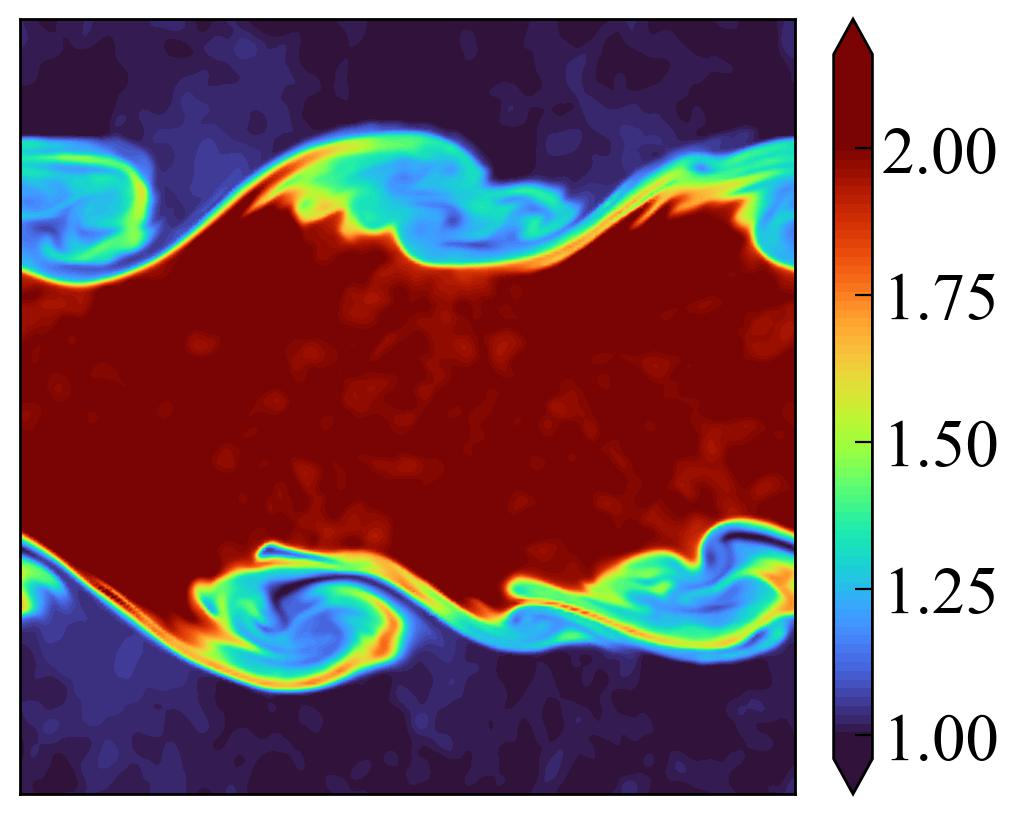}\hspace*{0.45cm}
            \includegraphics[trim=0.0cm 0cm 2.8cm 0cm, clip, width=0.165\textwidth]{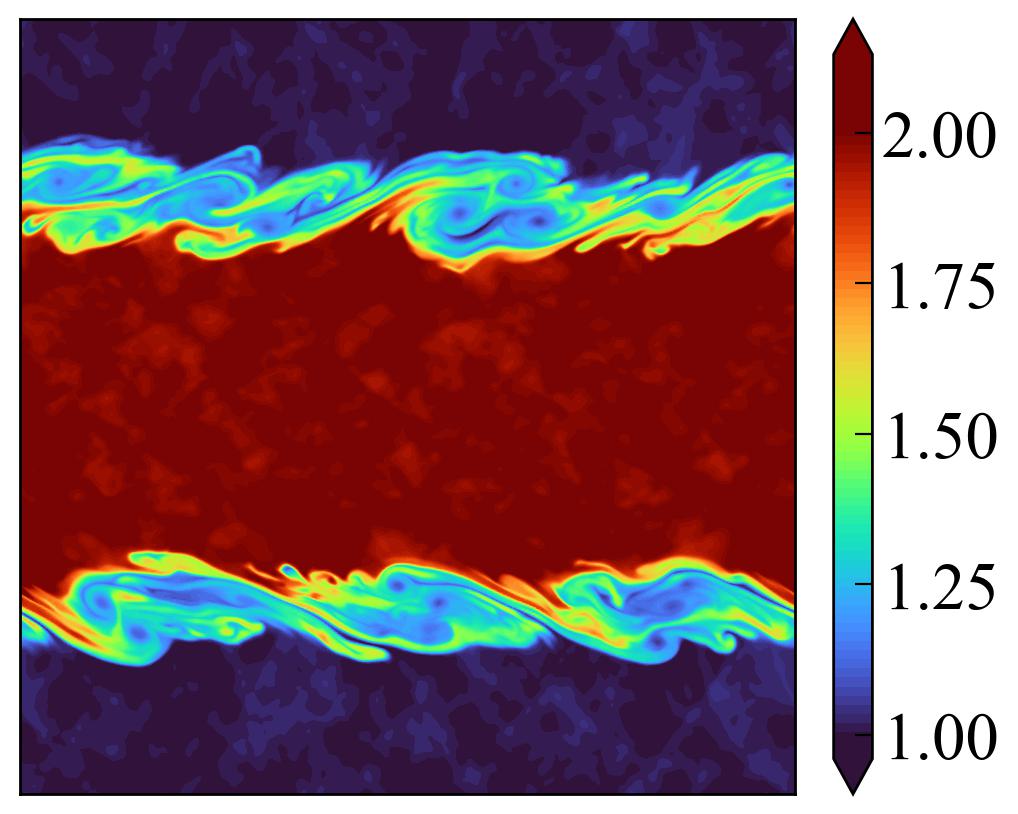}\hspace*{0.45cm}
            \includegraphics[trim=0.0cm 0cm 0.0cm 0cm, clip, width=0.21\textwidth]{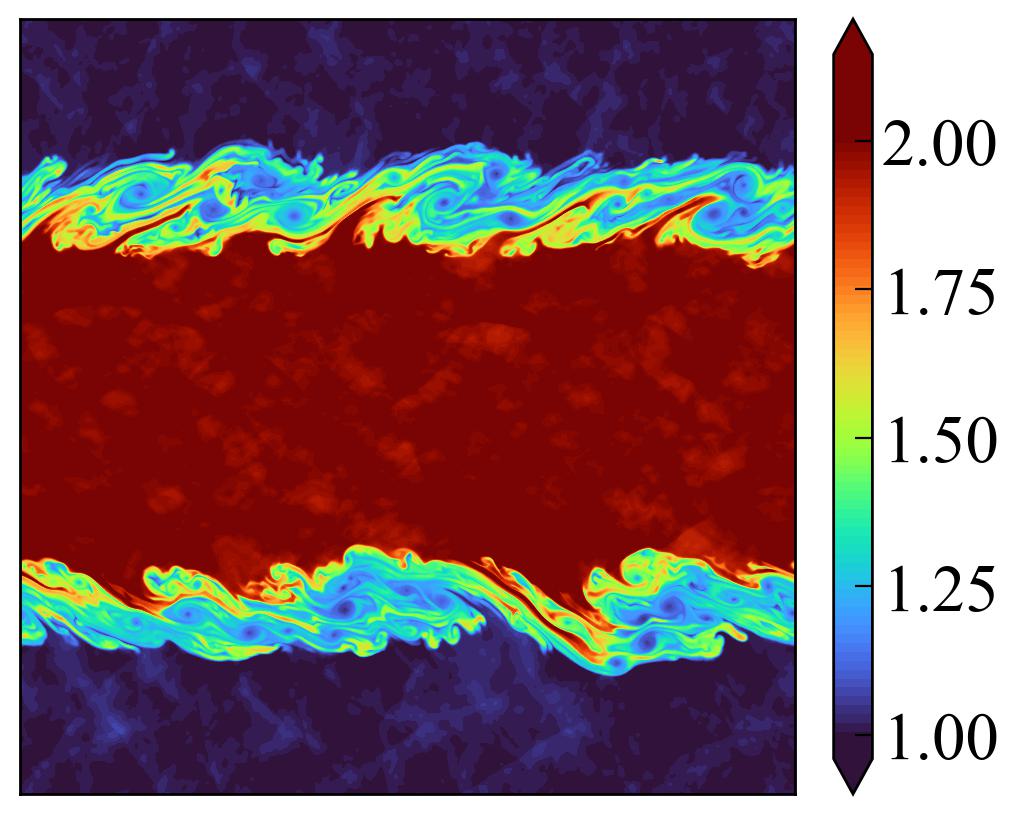}}
\caption{\sf $\rho(x,y,2;\xi)$ for $\xi=-1$ (top row), $0$ (middle row), and $1$ (bottom row) and five embedded uniform meshes
with $m=1,\dots,5$ (from left to right).\label{fig42}}
\end{figure}

\fref{fig41} shows the initial density distribution for the selected values of $\xi=-1$, $0$, and $1$ for five embedded uniform
meshes with $m=1,\dots,5$. We conduct simulations until the final time $T=2$, and plot, in \fref{fig42}, the obtained densities that
correspond to these initial conditions. As one can see, as the resolution increases, finer structures are resolved, indicating a more
accurate representation of the KH instability dynamics. Another observation is that larger values of $\xi$ (larger initial instability
amplitude) tend to produce finer structures that are more localized near the interface region. These increasingly finer roll-up structures,
especially at larger values of $\xi$, are reminiscent of the onset of turbulence, where coherent vortices undergo secondary instabilities
and break down into smaller scales. This behavior highlights the link between random KH instabilities and transitional turbulent mixing.

\begin{figure}[ht!]
\centering
\begin{subfigure}[b]{0.21\textwidth}
\centering
\hspace*{-0.4cm}\small $\langle\rho\rangle_5$\\
\centerline{\includegraphics[width=\textwidth,trim=0.0cm 0cm 0.2cm 0.0cm,clip]{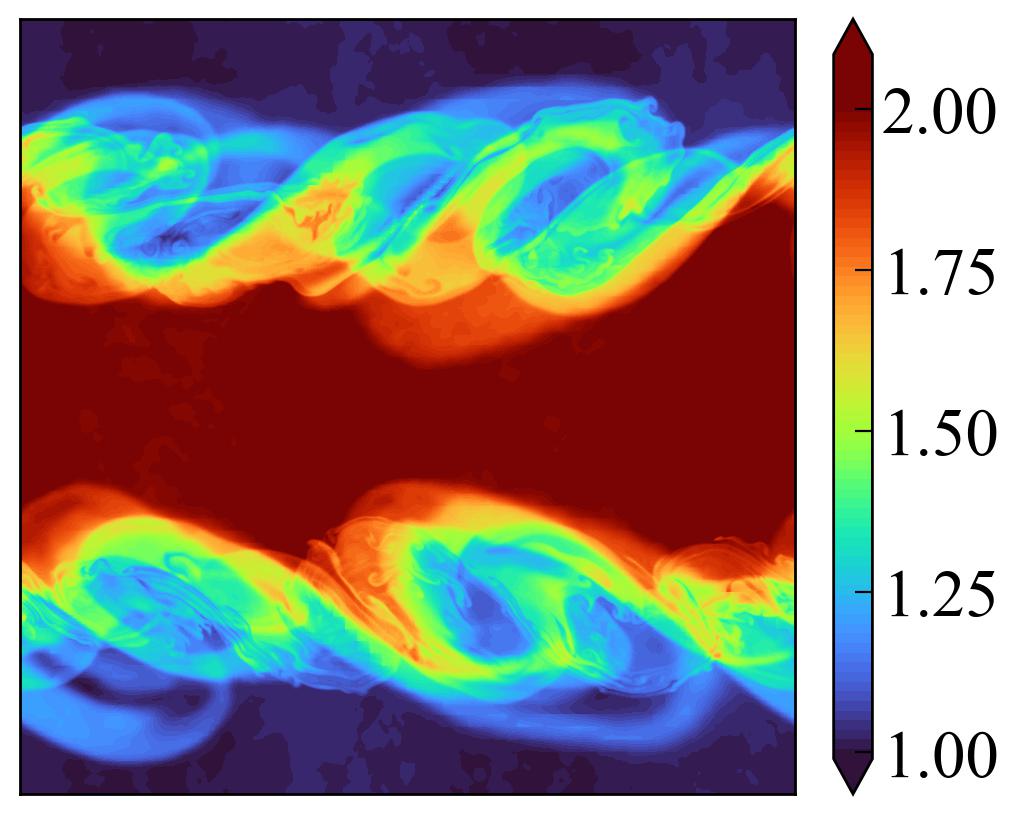}}
\end{subfigure}
\hfill
\begin{subfigure}[b]{0.21\textwidth}
\centering
\hspace*{-0.5cm}\small $\langle\rho u\rangle_5$\\
\includegraphics[width=\textwidth,trim=0cm 0cm 0.2cm 0cm,clip]{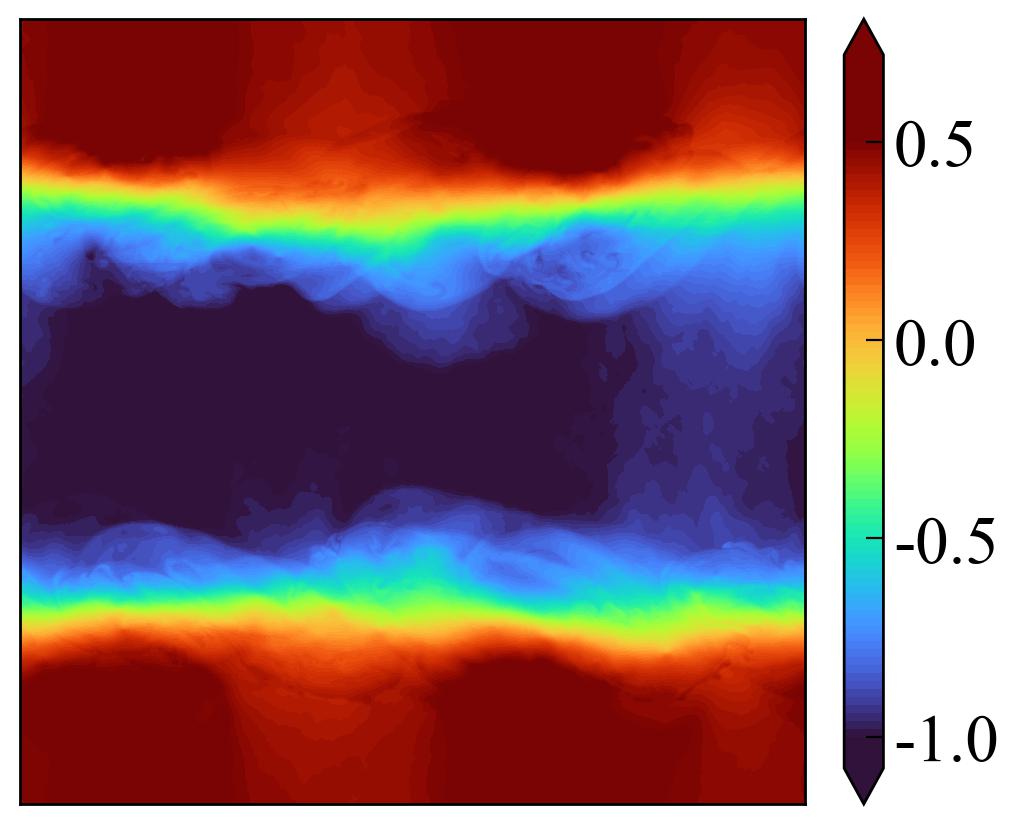}
\end{subfigure}
\hfill
\begin{subfigure}[b]{0.22\textwidth}
\centering
\hspace*{-0.6cm}\small $\langle\rho v\rangle_5$\\
\includegraphics[width=\textwidth,trim=0cm 0cm 0.2cm 0cm,clip]{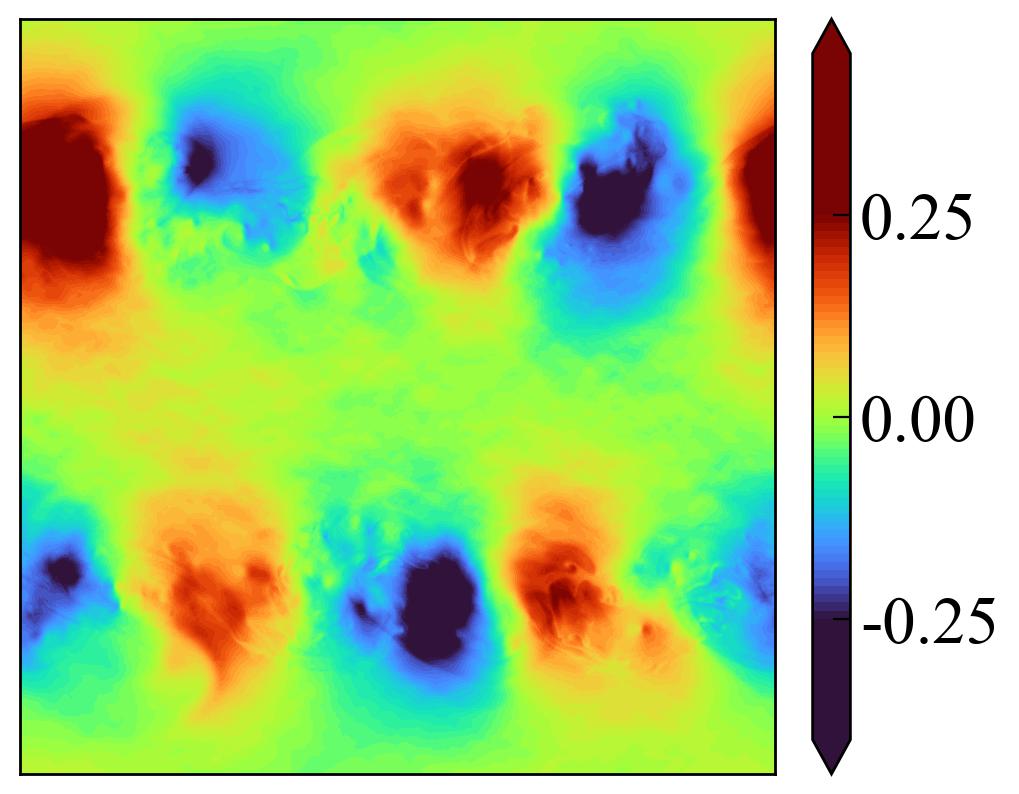}
\end{subfigure}
\hfill
\begin{subfigure}[b]{0.21\textwidth}
\centering
\hspace*{-0.4cm}\small $\langle S\rangle_5$\\
\includegraphics[width=\textwidth,trim=0cm 0cm 0.2cm 0cm,clip]{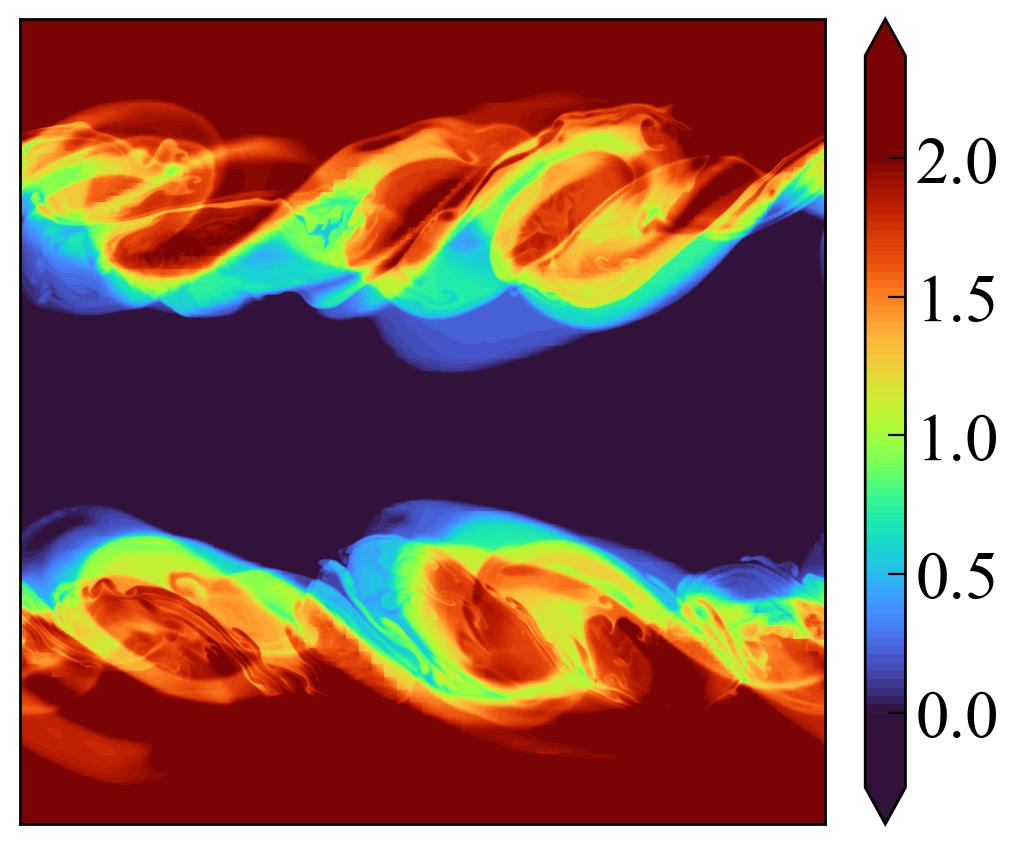}
\end{subfigure}
\caption{\sf Ces\`aro averages for the deterministic problem with $\xi=0$.\label{fig43}}
\end{figure}

\begin{figure}[ht!]
\centering
\begin{subfigure}[b]{0.21\textwidth}
\centering
\hspace*{-0.4cm}\small $\xbar{\langle\rho\rangle}_5$\\
\includegraphics[width=\textwidth,trim=0cm 0cm 0.2cm 0cm,clip]{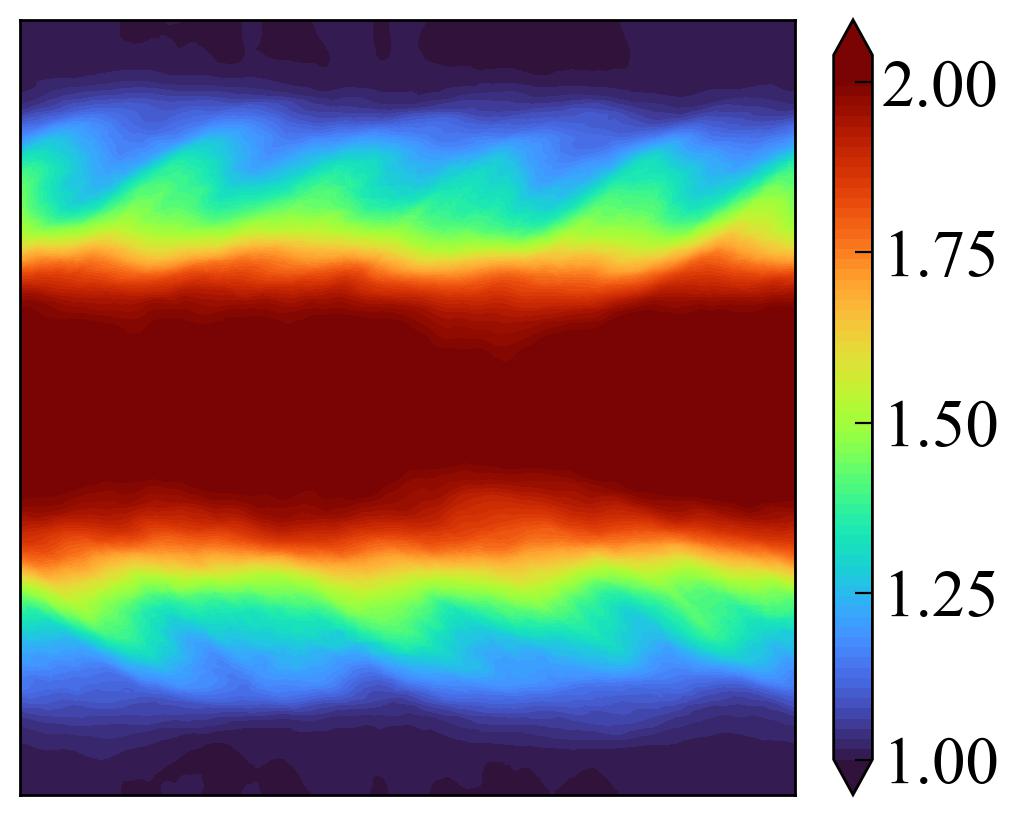}
\end{subfigure}
\hfill
\begin{subfigure}[b]{0.21\textwidth}
\centering
\hspace*{-0.5cm}\small $\xbar{\langle\rho u\rangle}_5$\\
\includegraphics[width=\textwidth,trim=0cm 0cm 0.2cm 0cm,clip]{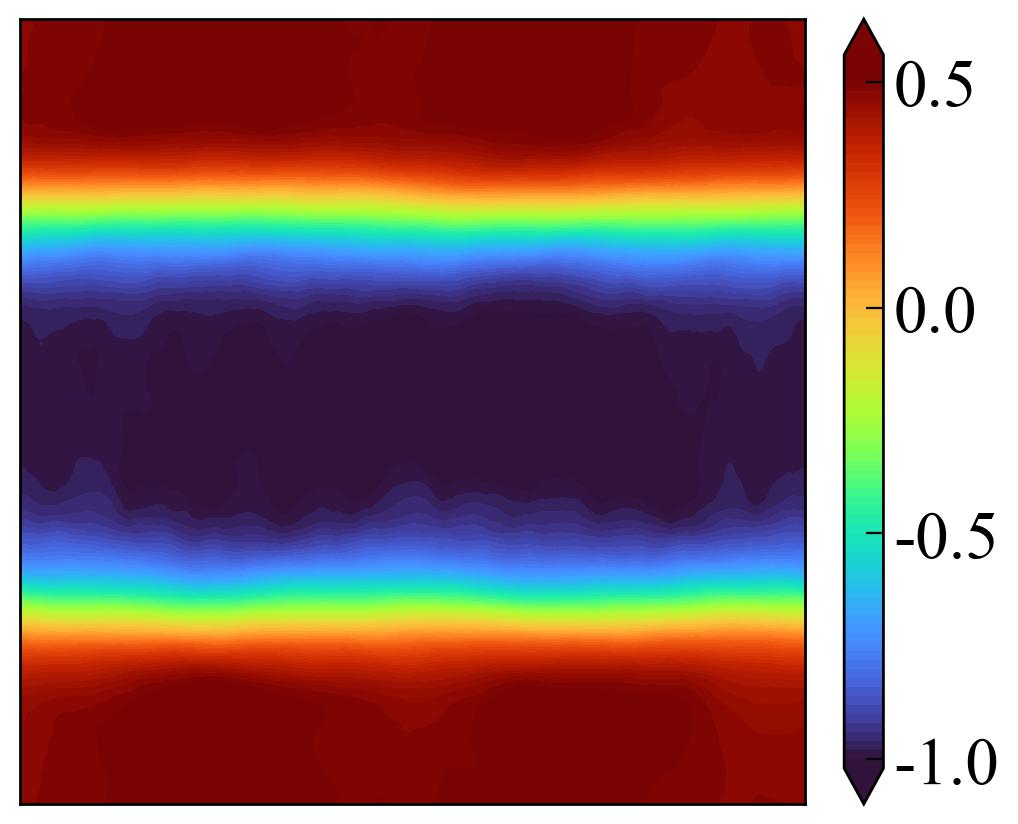}
\end{subfigure}
\hfill
\begin{subfigure}[b]{0.22\textwidth}
\centering
\hspace*{-0.6cm}\small $\xbar{\langle\rho v\rangle}_5$\\
\includegraphics[width=\textwidth,trim=0cm 0cm 0.2cm 0cm,clip]{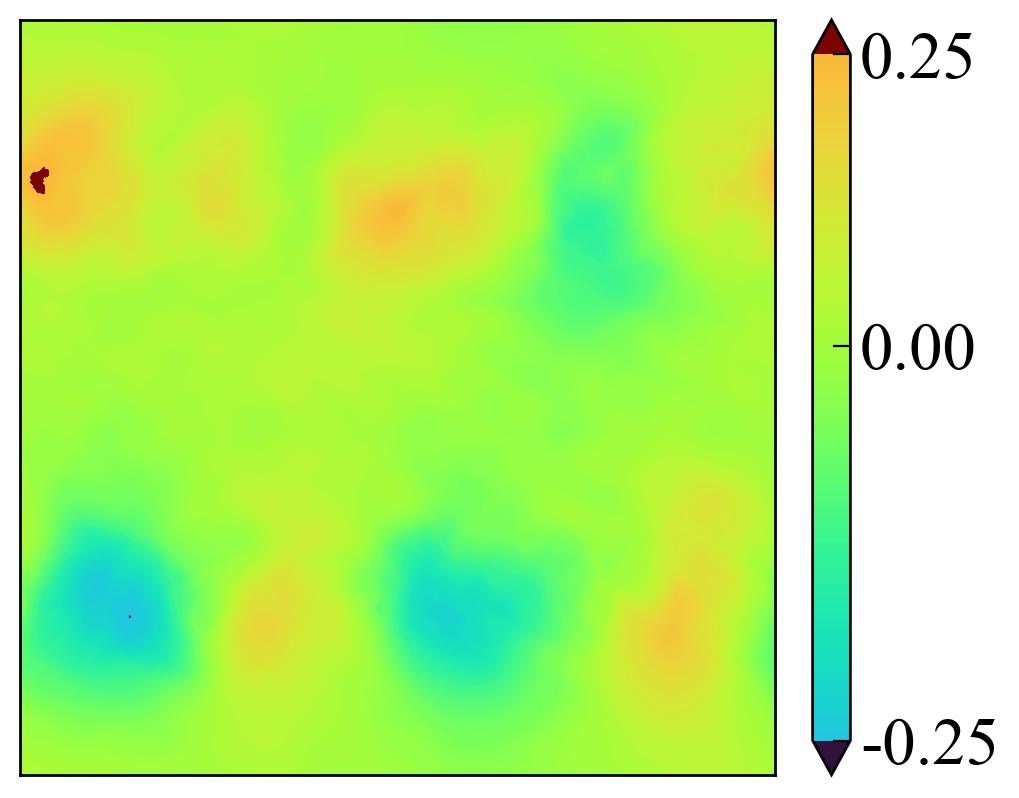}
\end{subfigure}
\hfill
\begin{subfigure}[b]{0.21\textwidth}
\centering
\hspace*{-0.4cm}\small $\xbar{\langle S\rangle}_5$\\
\includegraphics[width=\textwidth,trim=0cm 0cm 0.2cm 0cm,clip]{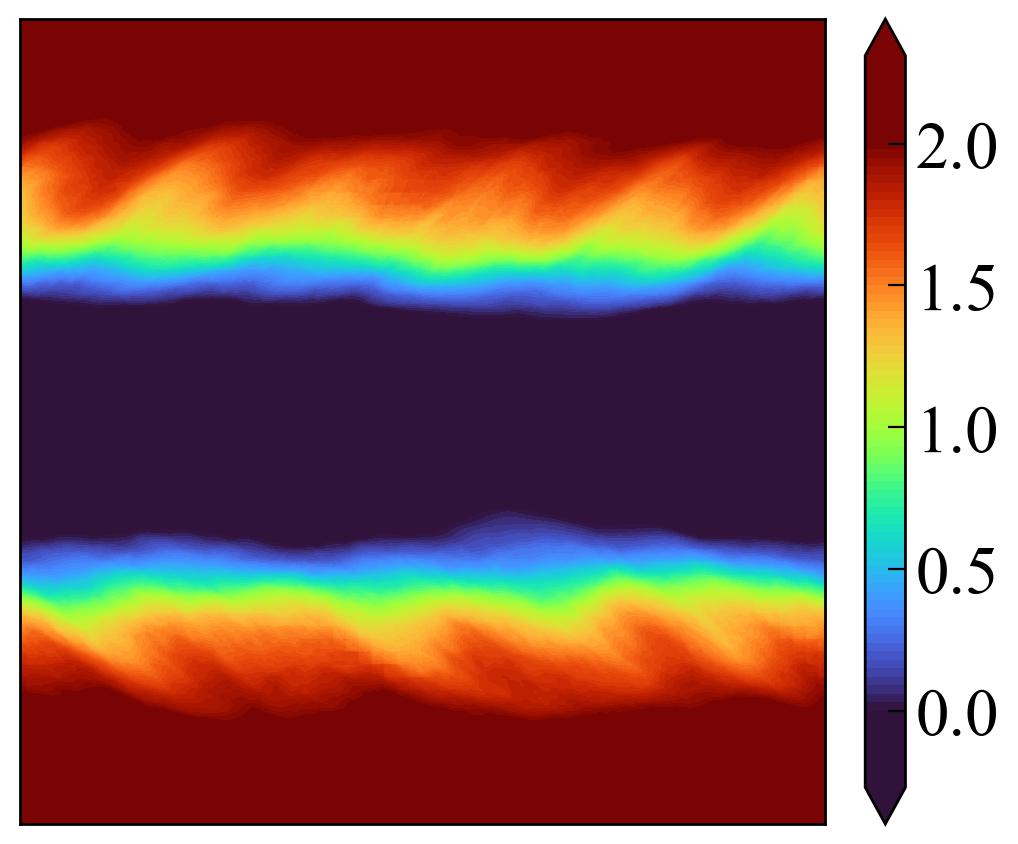}
\end{subfigure}
\caption{\sf Means of the Ces\'aro averages for the stochastic problem.\label{fig44}}
\end{figure}

\paragraph{Ces\`aro averages.} We compute the Ces\`aro averages defined in \eref{2.8f} for $M=5$ for the deterministic problem with a
specific value $\xi=0$, and plot the obtained results in \fref{fig43}. These Ces\`aro averages provide an approximation of the DW solutions
in the strong sense specified in Theorem \ref{thm21}. We then compute the mean values of the Ces\`aro averages 
$\langle\rho\rangle_5, \langle\rho u\rangle_5,\langle\rho v\rangle_5$, and $\langle S\rangle_5$ as described in \S\ref{sec3}: Since $\xi$ is
uniformly distributed, the integrals in \eref{3.3} are computed exactly. The obtained results, presented in \fref{fig44}, are expected to
approximate the mean of the DW solution in the strong sense as indicated by Theorem \ref{thm22}.

\begin{figure}[ht!]
\centering
\begin{subfigure}[t]{0.3\textwidth}
\centering
\text{\small $\text{tr}(\overline{\mathfrak{R}}_5)$}\\
\includegraphics[width=0.75\textwidth,trim=0cm 0cm 0.2cm 0cm,clip]{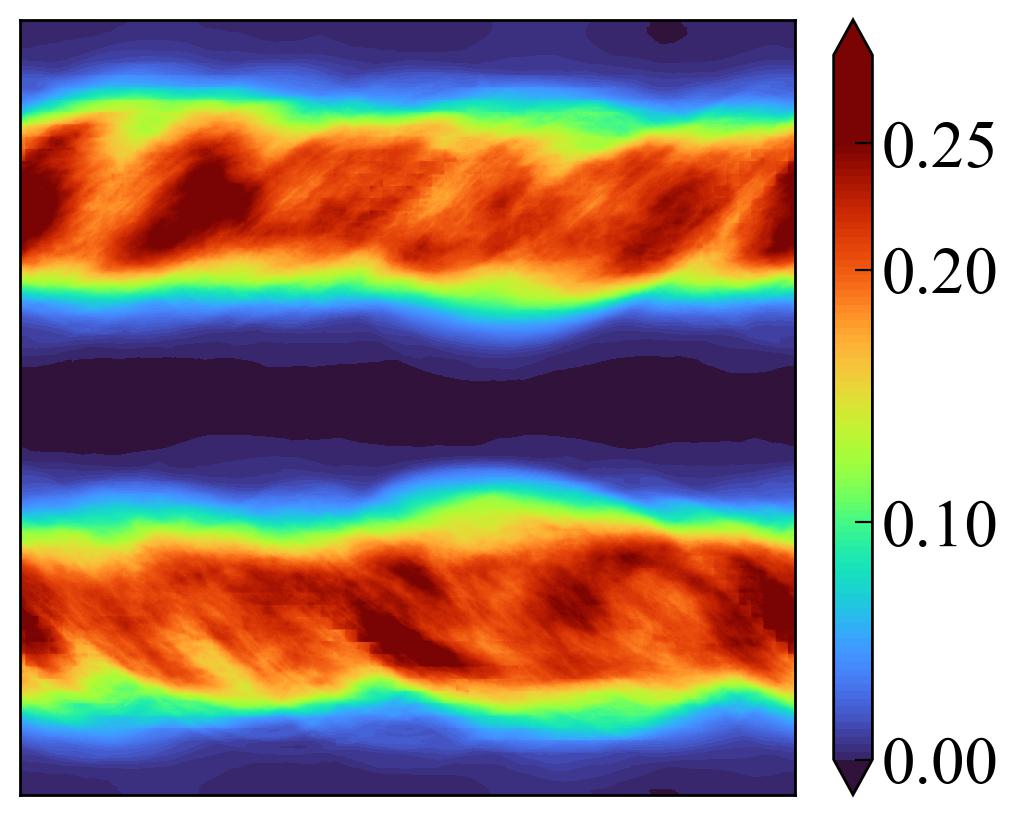}
\end{subfigure}
\begin{subfigure}[t]{0.3\textwidth}
\centering
\text{\small $\overline{\mathfrak{E}}_5$}\\
\includegraphics[width=0.75\textwidth,trim=0cm 0cm 0.2cm 0cm,clip]{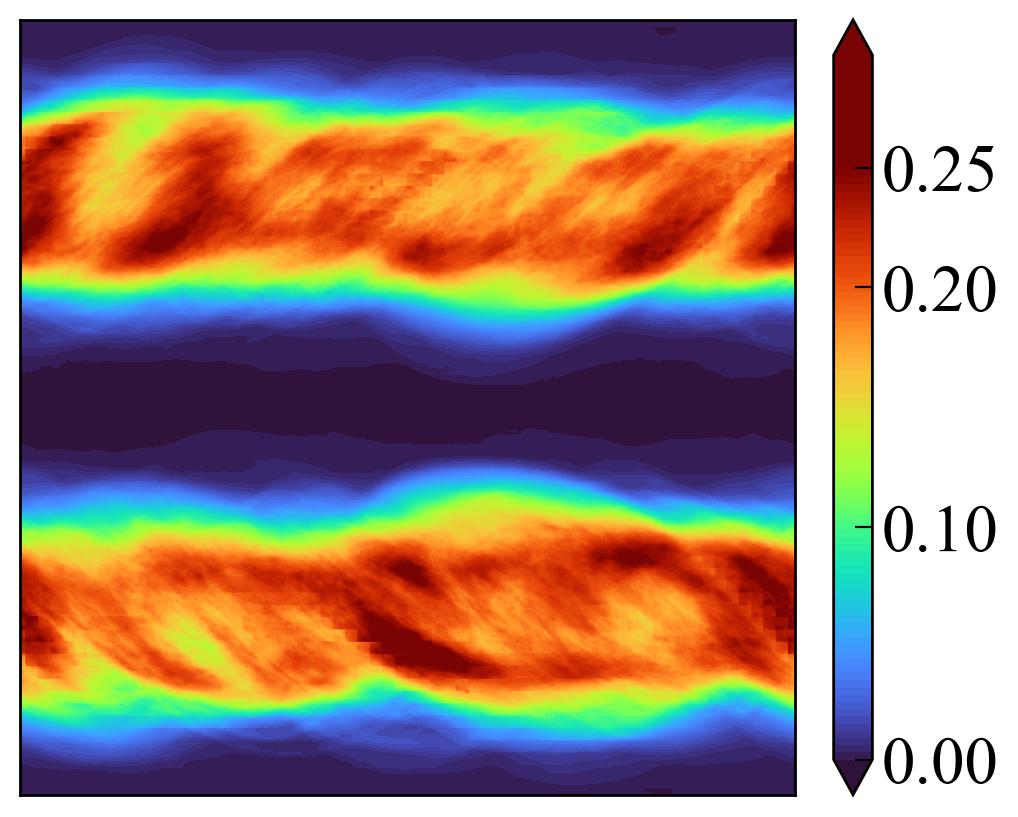}
\end{subfigure}
\begin{subfigure}[t]{0.3\textwidth}
\centering
\text{\small ${\overline{\mathfrak{E}}_5}/{\text{tr}(\overline{\mathfrak{R}}_5)}$}\\
\includegraphics[width=0.75\textwidth,trim=0cm 0cm 0.2cm 0cm,clip]{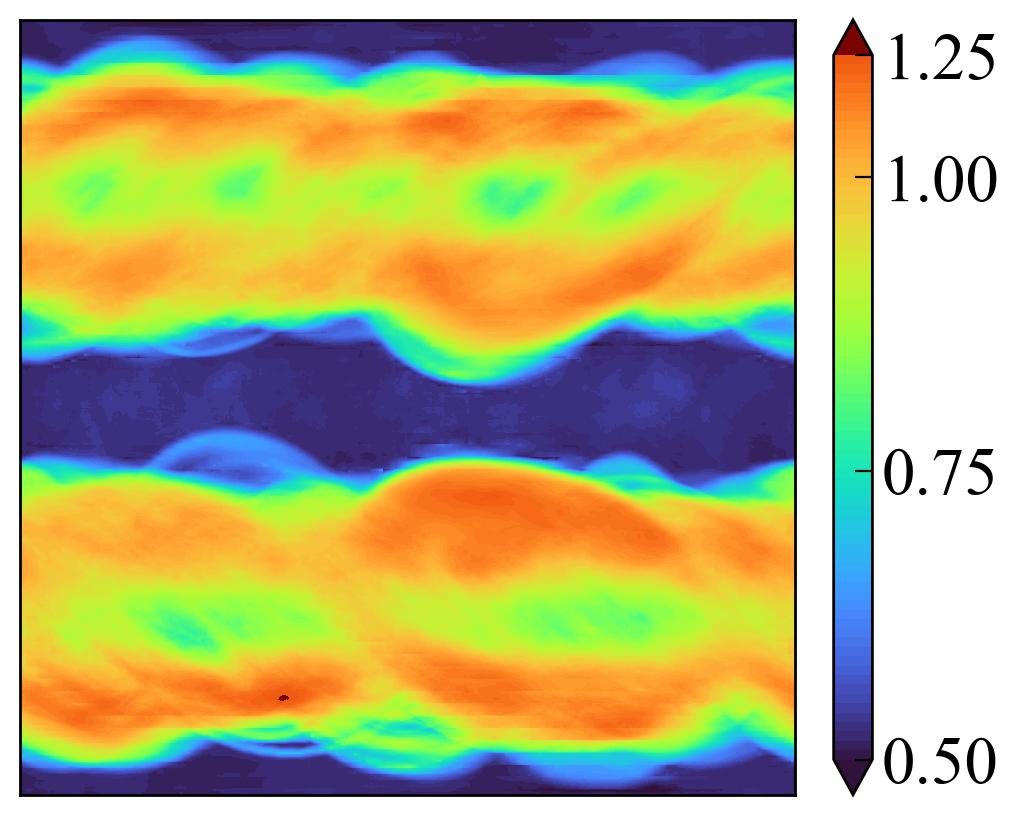}
\end{subfigure}
\caption{\sf Reynolds stress and energy defect, along with their ratio.\label{fig45}}
\end{figure}

\paragraph{Reynolds stress and energy defects.} \fref{fig45} shows the mean of the trace of the Reynolds stress defect
${\rm tr}(\,\xbar{\mathfrak R}_5)$ and the mean of the energy defect $\,\xbar{\mathfrak E}_5$ computed according to \eref{3.4}, as well as
their ratio $\,\xbar{\mathfrak E}_5/{\rm tr}(\,\xbar{\mathfrak R}_5)$. As expected, the Reynolds stress and energy defect have similar
structures and their ratio stays within the theoretically bounds $0.5$ and $1.25$, specified in \eref{2.2f}.

\begin{figure}[ht!]
\centerline{\includegraphics[width=0.37\textwidth,trim=0.9cm 0.8cm 0.6cm 0.7cm,clip]{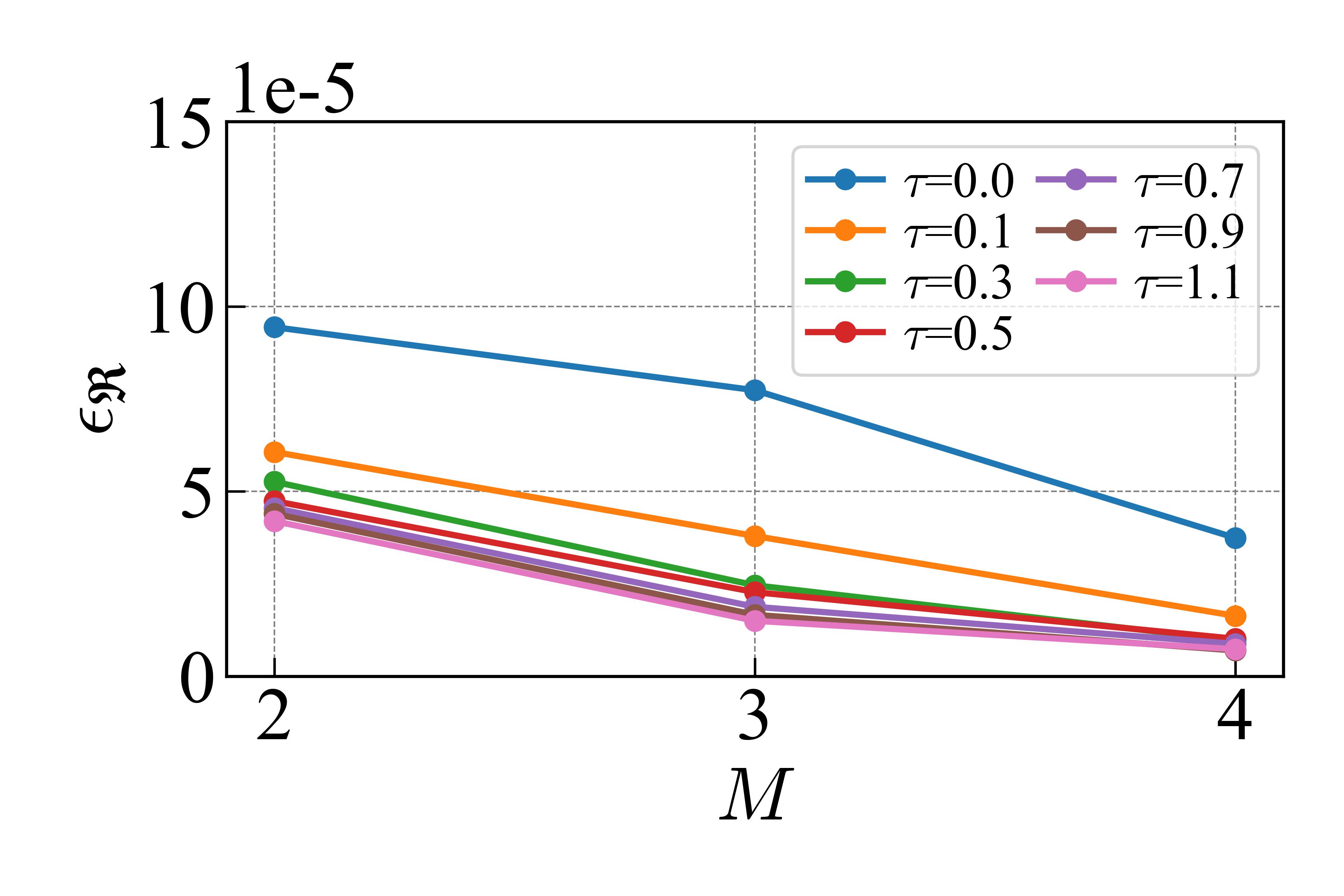}\hspace*{0.8cm}
            \includegraphics[width=0.37\textwidth,trim=0.9cm 0.8cm 0.6cm 0.7cm,clip]{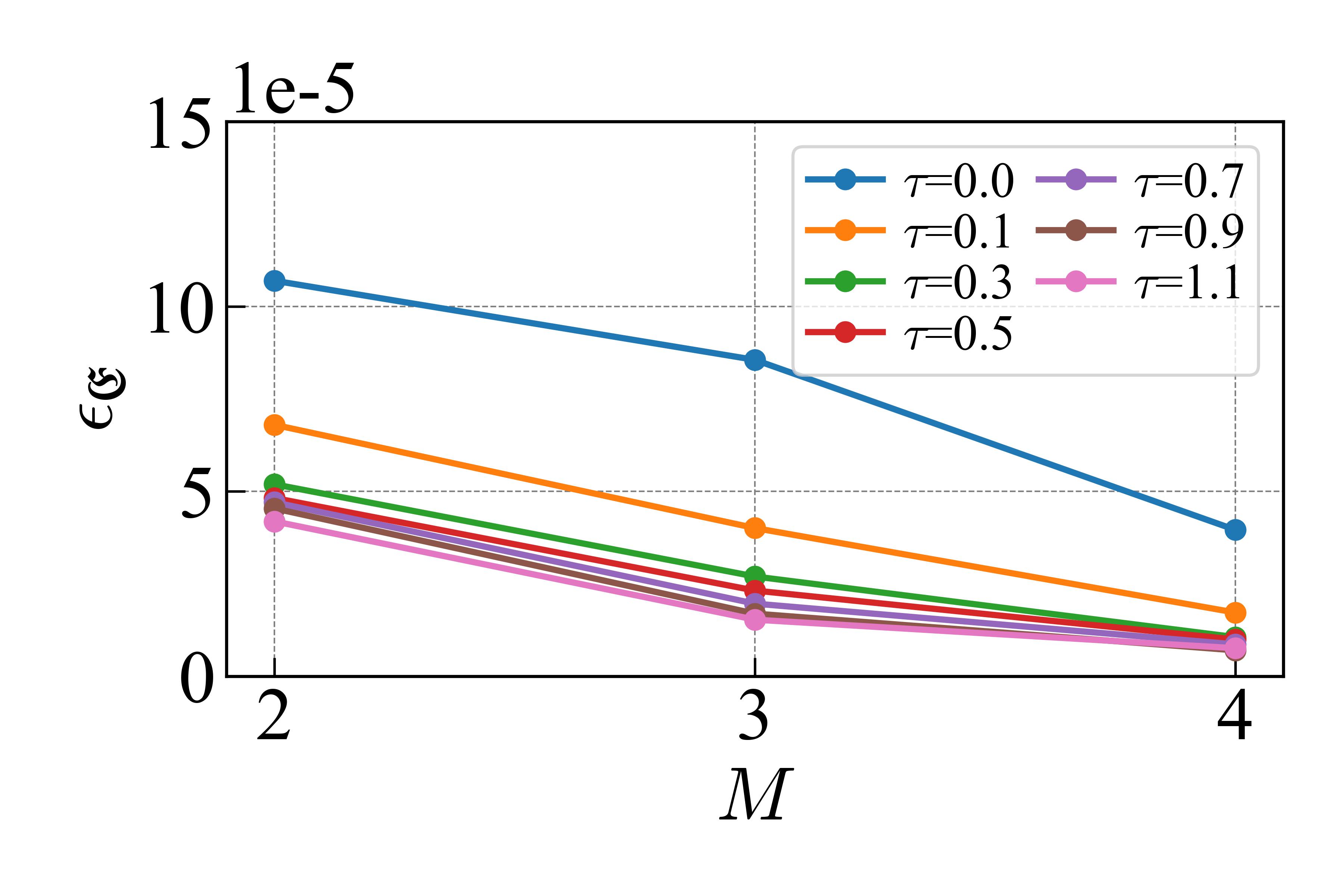}}
\centerline{\includegraphics[width=0.30\textwidth,trim=0.9cm 0.8cm 0.6cm 0.75cm,clip]{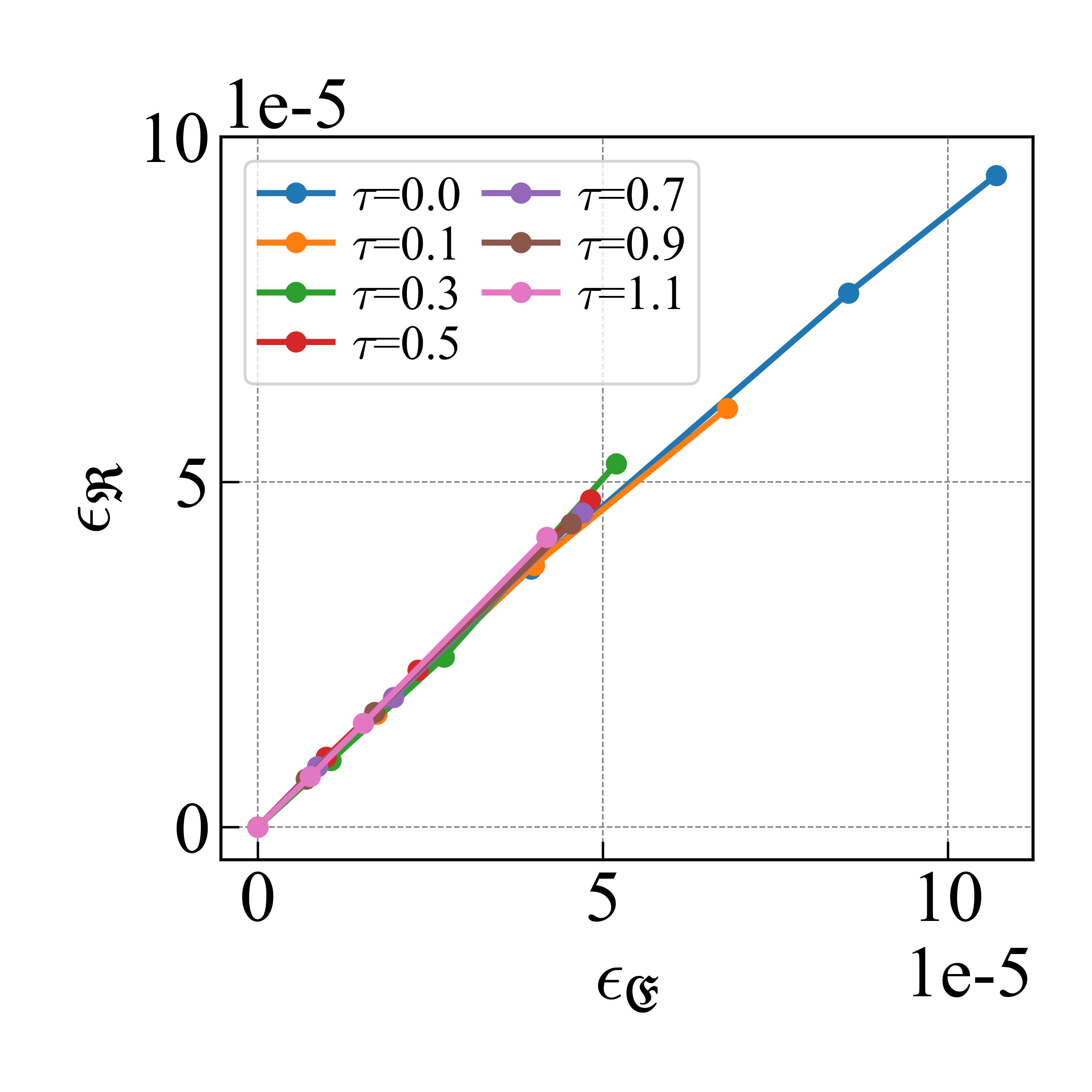}}
\caption{\sf $\epsilon_{\mathfrak R}$ (top left) and $\epsilon_{\mathfrak E}$ (top right) as functions of $M$; $\epsilon_{\mathfrak E}$ as a
function of $\epsilon_{\mathfrak R}$ (bottom).\label{fig46}}
\end{figure}

We also experimentally study the convergence of these defects by measuring the following quantities:
\begin{equation*}
\epsilon_{\mathfrak R_M}=\|{\rm tr}(\,\xbar{\mathfrak R}_M)-{\rm tr}(\,\xbar{\mathfrak R}_5)||_1,\quad
\epsilon_{\mathfrak E_M}=||\,\xbar{\mathfrak E}_M-\,\xbar{\mathfrak E}_5||_1,
\end{equation*}
for $M=2,3,4$ and for different initial perturbations $\tau$ ranging from $0$ to $1.1$; see \fref{fig46}. One can observe that both
$\epsilon_{\mathfrak R_M}$ and $\epsilon_{\mathfrak E_M}$ decrease steadily as $M$ increases, and that for all $M$, larger values of the
instability parameter $\tau$ yields smaller residual norms. In addition, we plot $\epsilon_{\mathfrak R_M}$ as a functions of
$\epsilon_{\mathfrak E_M}$ demonstrating the near-unity slope, which indicates the same convergence rate for both the Reynolds stresses and
energy defect. The latter confirms that these turbulence-like quantities stabilize under Ces\`aro averaging, supporting the statistical
framework as a robust description of turbulent variability.

\begin{figure}[ht!]
\centering
\begin{subfigure}[b]{0.305\textwidth}
\centering
\text{\small $\xbar{\langle\rho\rangle}_5$}\\
\includegraphics[trim=0.2cm 0.8cm 0.6cm 0.6cm,clip,width=\textwidth]{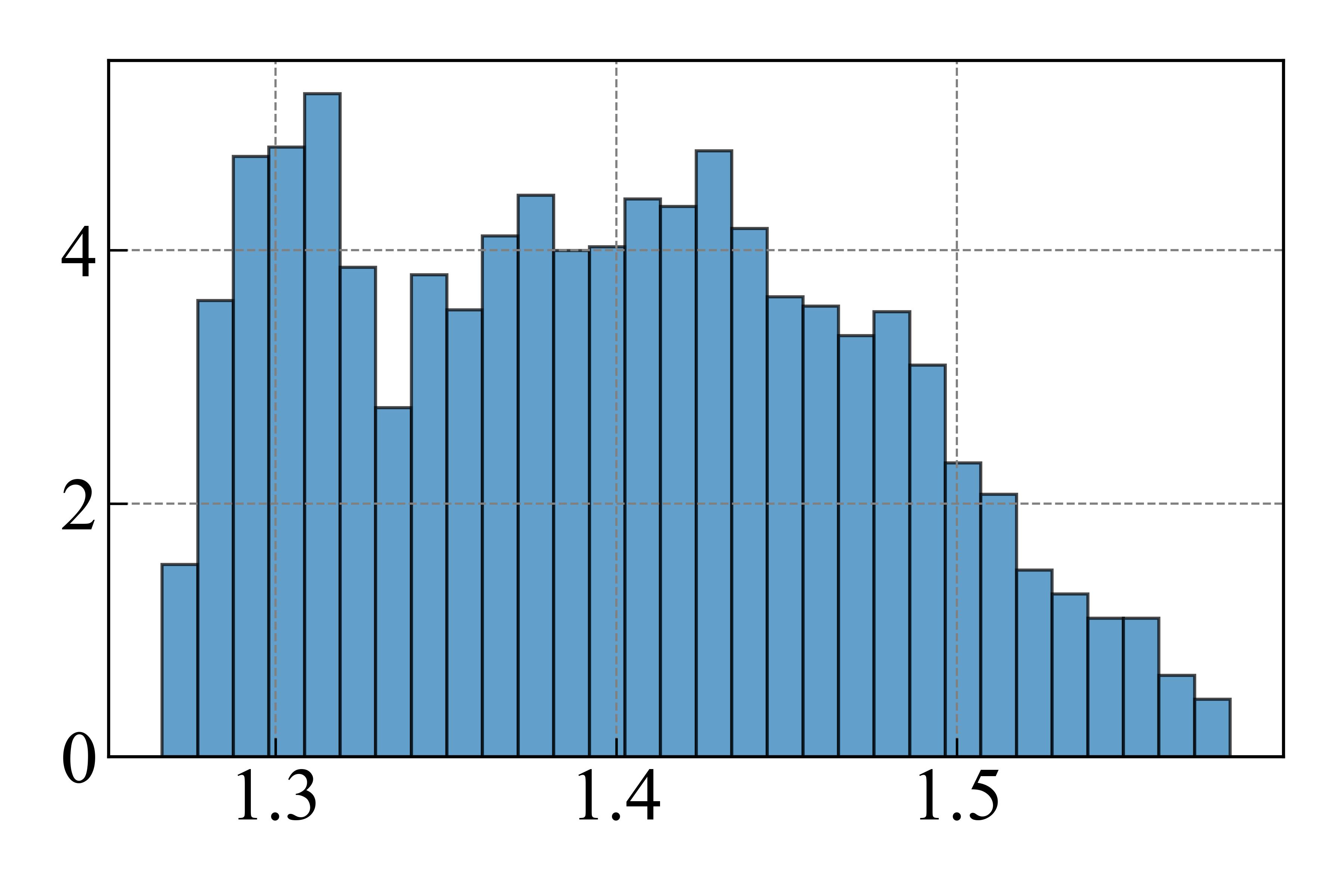}
\end{subfigure}
\hfill
\begin{subfigure}[b]{0.305\textwidth}
\centering
\text{\small $\xbar{\langle S\rangle}_5$}\\
\includegraphics[trim=0.2cm 0.8cm 0.6cm 0.6cm,clip,width=\textwidth]{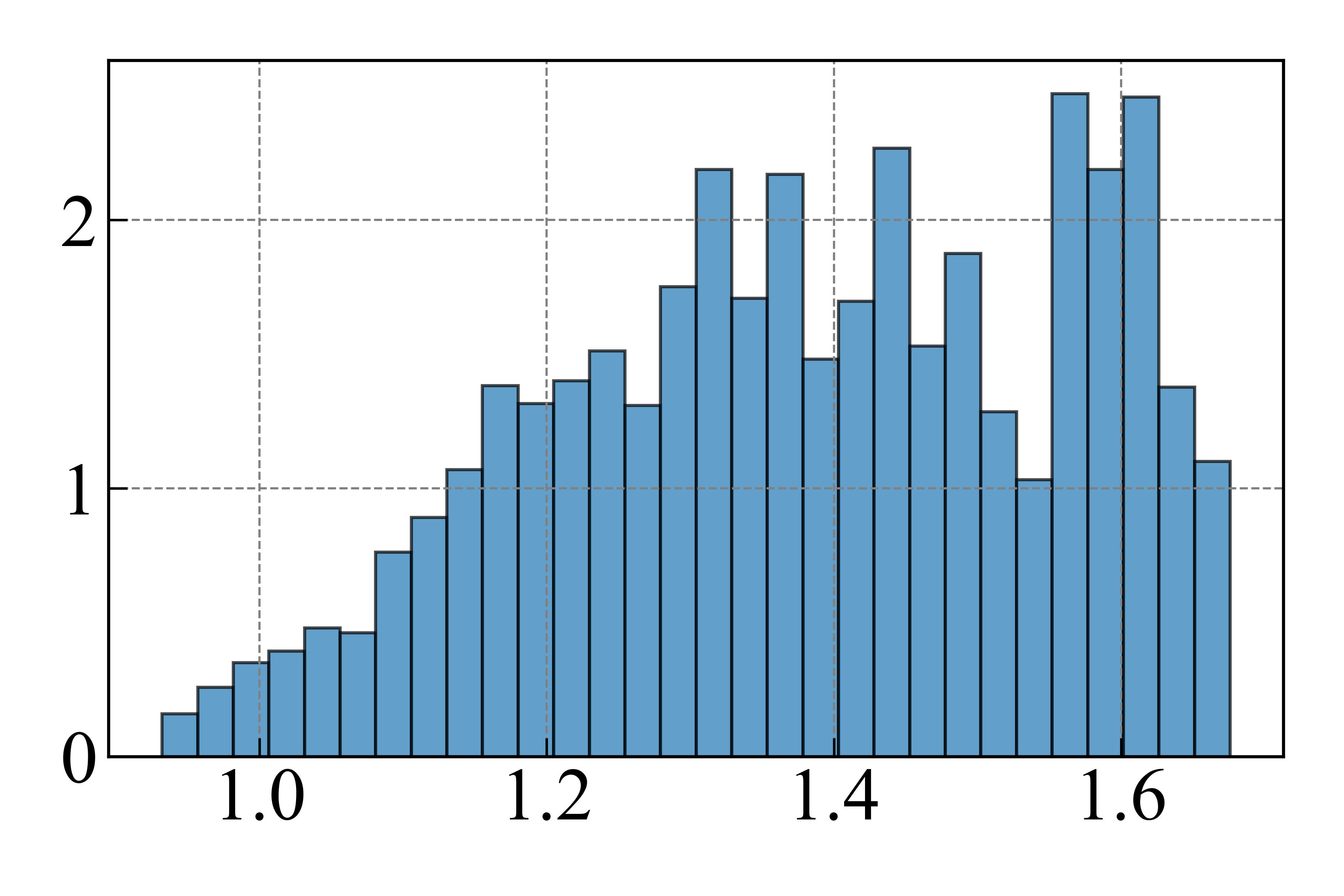}
\end{subfigure} 
\hfill
\begin{subfigure}[b]{0.305\textwidth}
\centering
\text{\small ${\rm tr}(\,\xbar{\mathfrak R}_5)$}\\
\includegraphics[trim=0.2cm 0.8cm 0.6cm 0.6cm,clip,width=\textwidth]{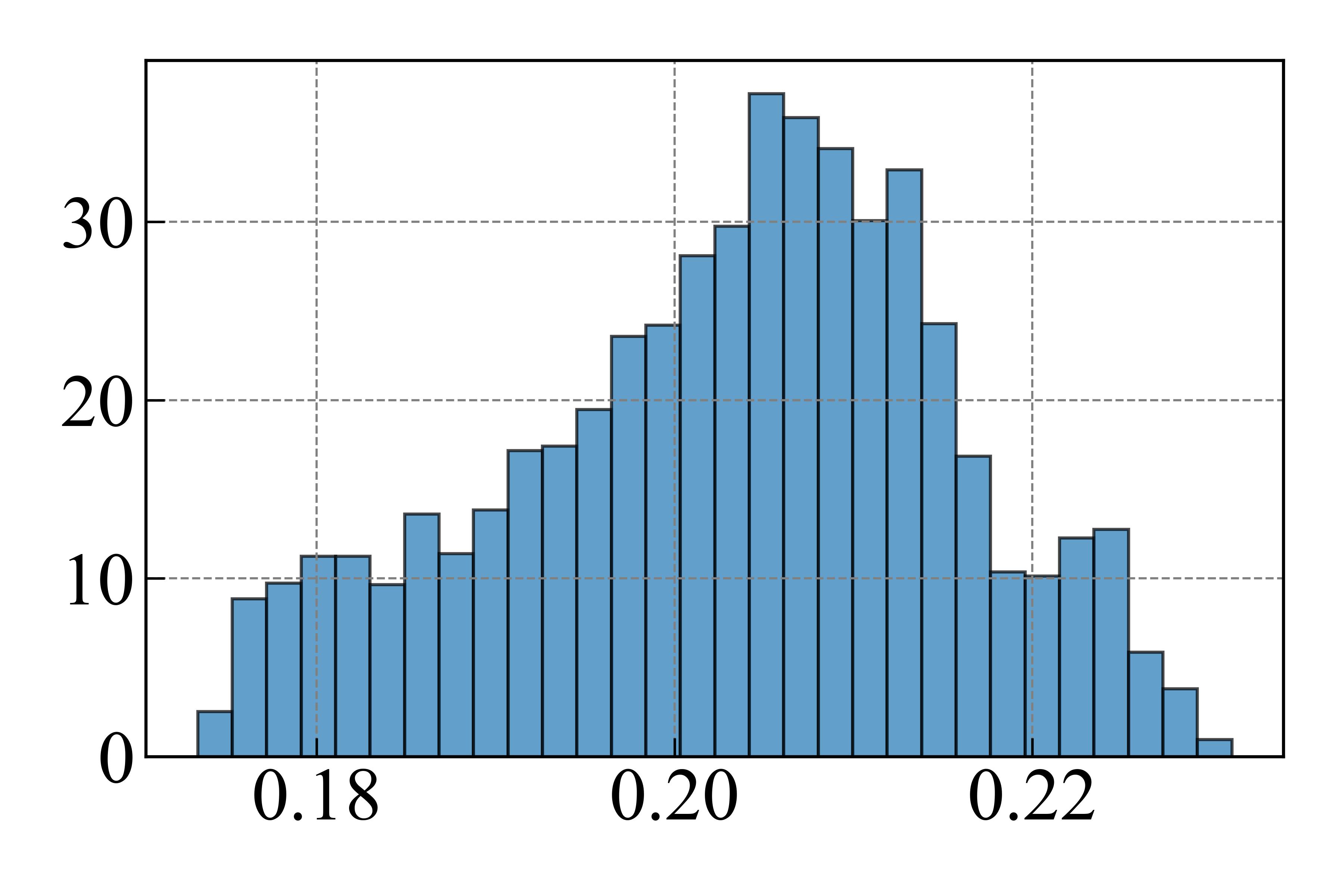}
\end{subfigure}
\hfill
\vskip8pt
\centering
\begin{subfigure}[b]{0.305\textwidth}
\centering
\text{\small $\xbar{\langle\rho\rangle}_5$}\\
\includegraphics[trim=0.2cm 0.8cm 0.6cm 0.6cm,clip,width=\textwidth]{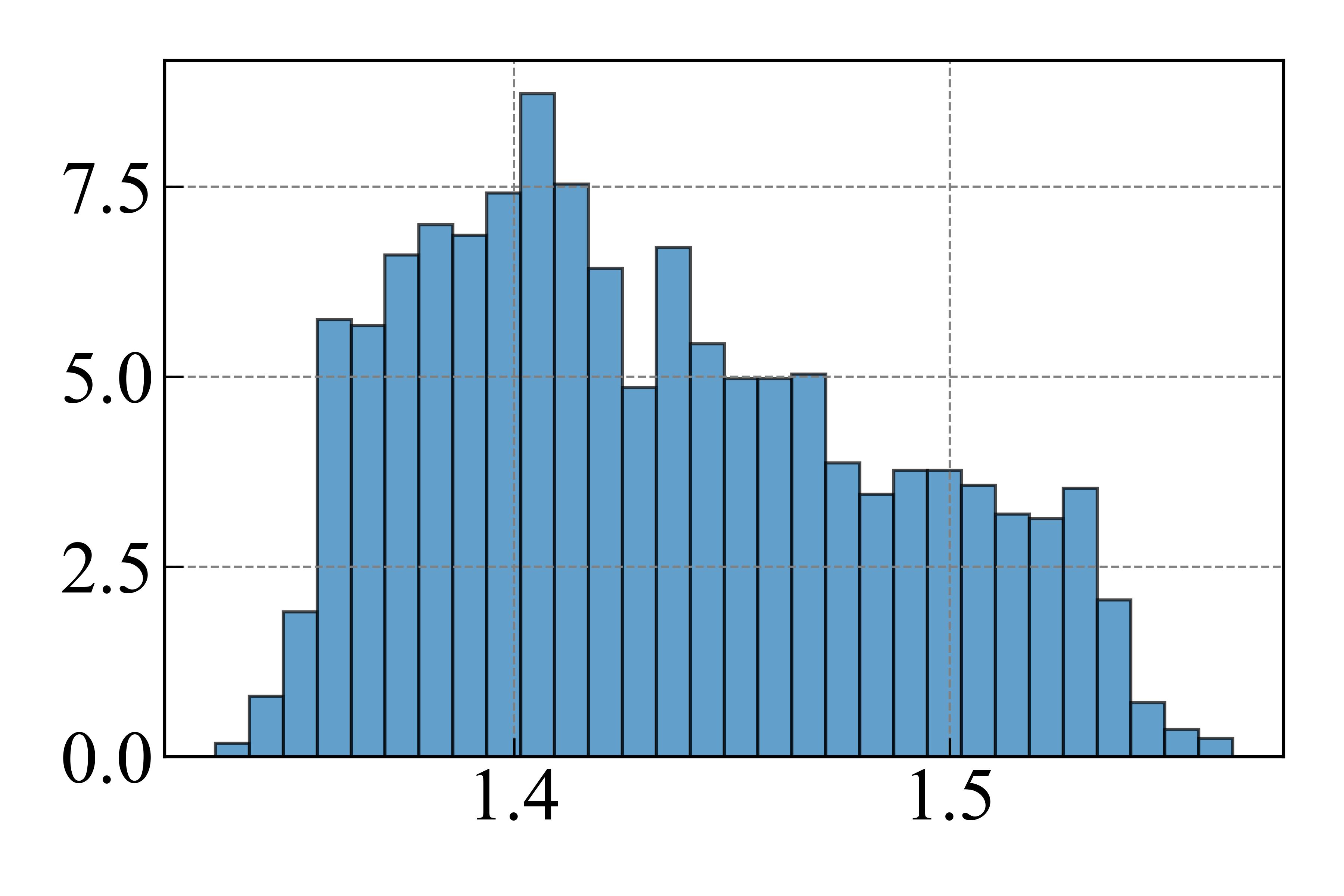}
\end{subfigure}
\hfill
\begin{subfigure}[b]{0.305\textwidth}
\centering
\text{\small $\xbar{\langle S\rangle}_5$}\\
\includegraphics[trim=0.2cm 0.8cm 0.6cm 0.6cm,clip,width=\textwidth]{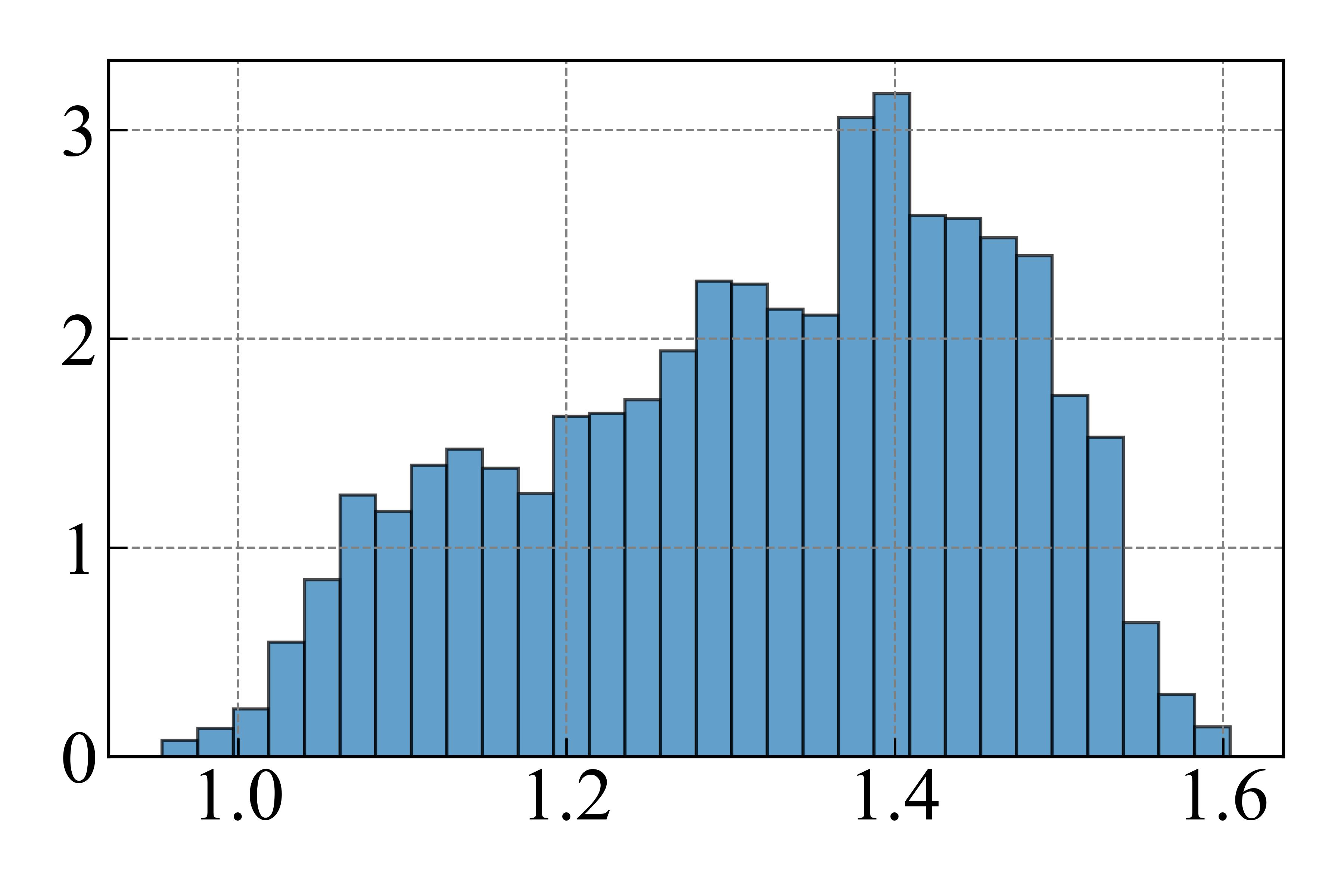}
\end{subfigure} 
\hfill
\begin{subfigure}[b]{0.305\textwidth}
\centering
\text{\small ${\rm tr}(\,\xbar{\mathfrak R}_5)$}\\
\includegraphics[trim=0.2cm 0.8cm 0.6cm 0.6cm,clip,width=\textwidth]{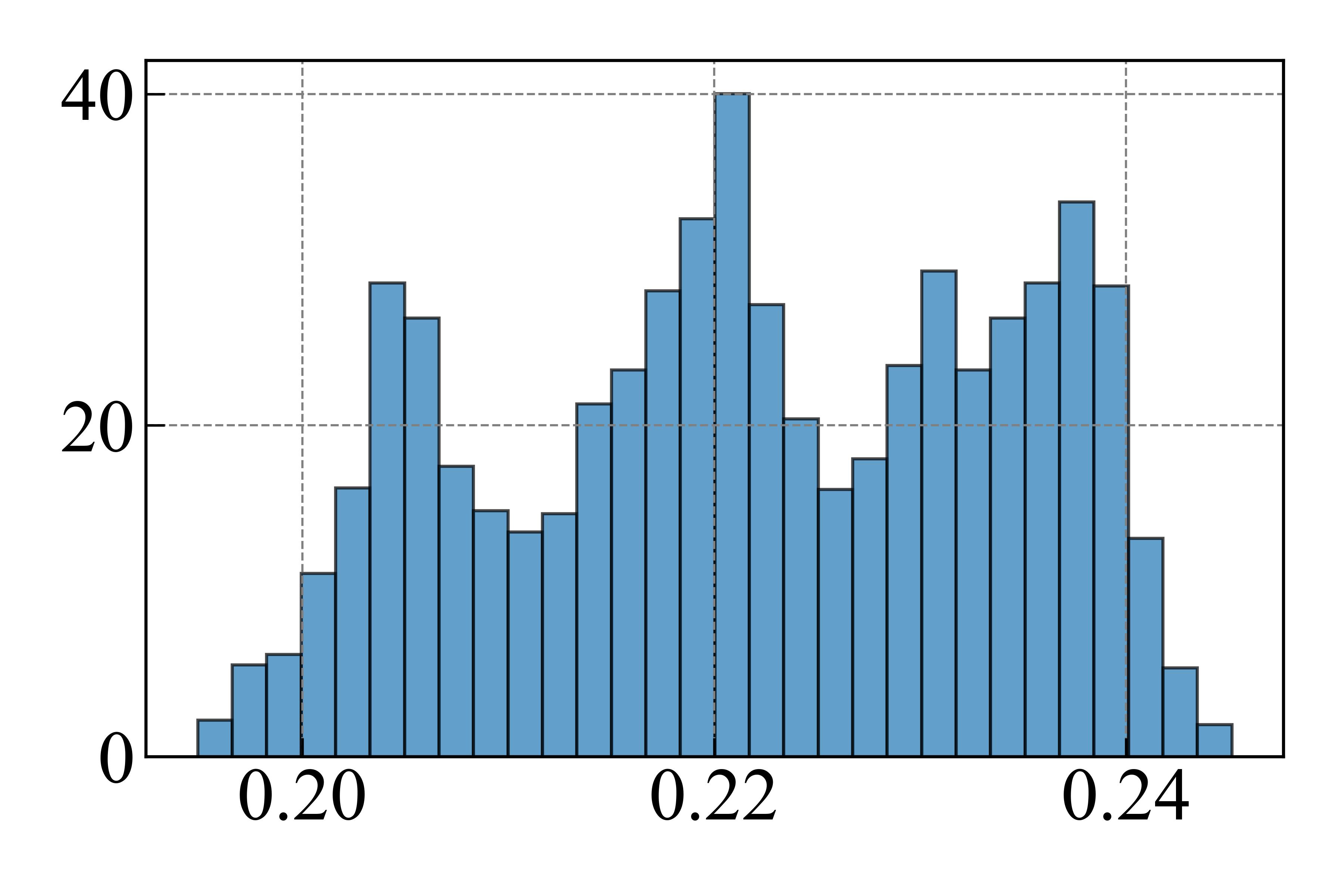}
\end{subfigure} 
\caption{\sf PDFs of $\,\xbar{\langle\rho\rangle}_5$, $\xbar{\langle S\rangle}_5$, and ${\rm tr}(\,\xbar{\mathfrak R}_5)$ approximated in
$D_1$ (top row) and $D_2$ (bottom row) using \texttt{numpy.histogram} function in \texttt{Python} with the \texttt{auto} binning strategy.
\label{fig47}}
\end{figure}

\paragraph{Statistical analysis.}
We first construct the PDFs of $\,\xbar{\langle\rho\rangle}_5$, $\xbar{\langle S\rangle}_5$, and ${\rm tr}(\,\xbar{\mathfrak R}_5)$. To
this end, we select small spatial areas, on which these PDFs are approximated using a histogram approach. Specifically, we choose 
reasonably small spatial windows $D_1$ and $D_2$ outlined by yellow and green squares in \fref{fig41}. The corresponding PDFs, computed 
using \texttt{numpy.histogram} function in \texttt{Python} with the \texttt{auto} binning strategy, are depicted in \fref{fig47}. As one 
can see, the PDFs exhibit nontrivial spreads, that is, non-Dirac-type measure, emphasizing the persistent variability characteristic of 
turbulence \cite{MoninYaglom1971,Pope2000}.
 
\begin{figure}[ht!]
\centering
\begin{subfigure}[b]{0.305\textwidth}
\centering
\text{\small $\rho$}\\
\includegraphics[trim=0.2cm 0.8cm 0.6cm 0.6cm,clip,width=\textwidth]{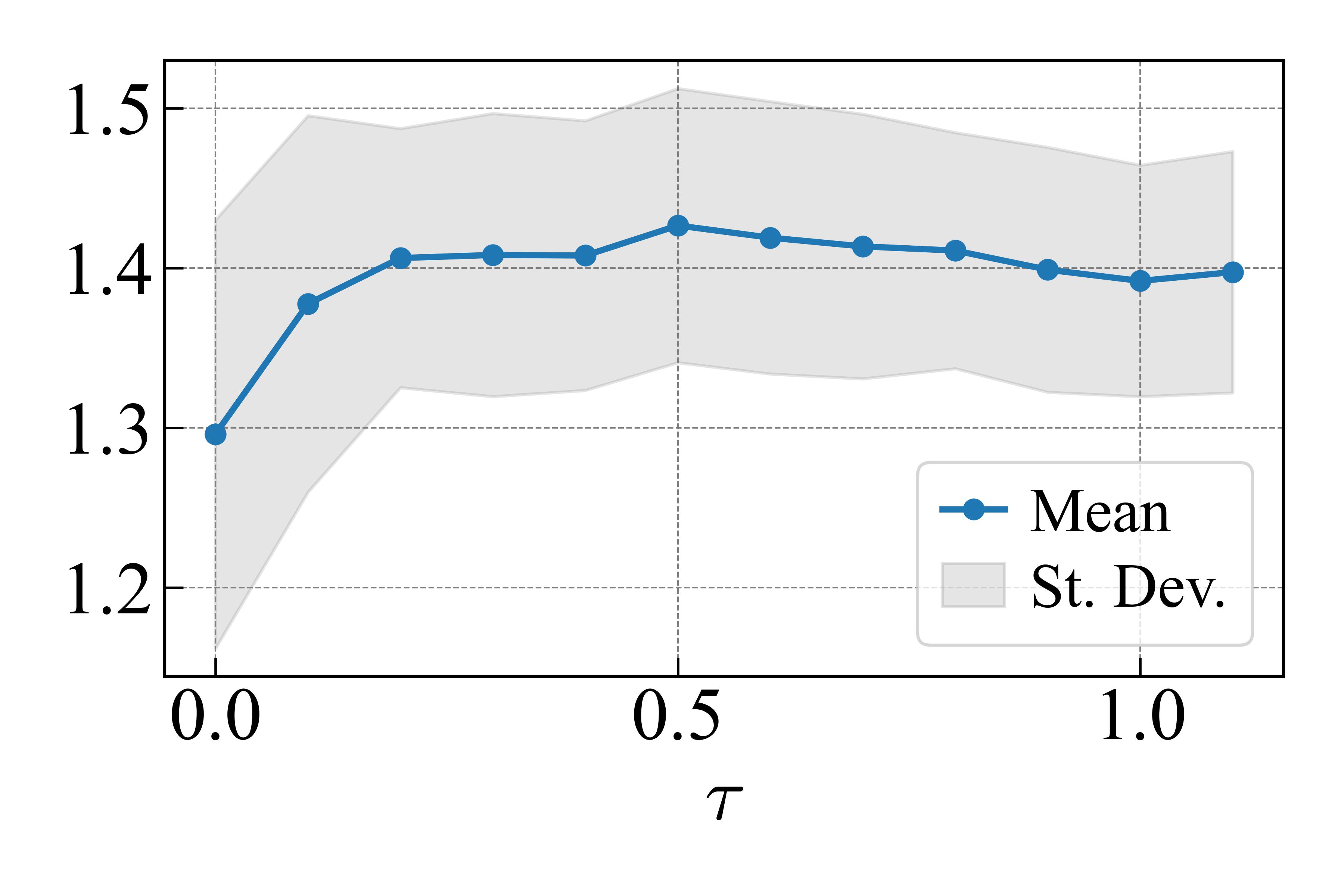}
\end{subfigure}
\hfill
\begin{subfigure}[b]{0.305\textwidth}
\centering
\text{\small $S$}\\
\includegraphics[trim=0.2cm 0.8cm 0.6cm 0.6cm,clip,width=\textwidth]{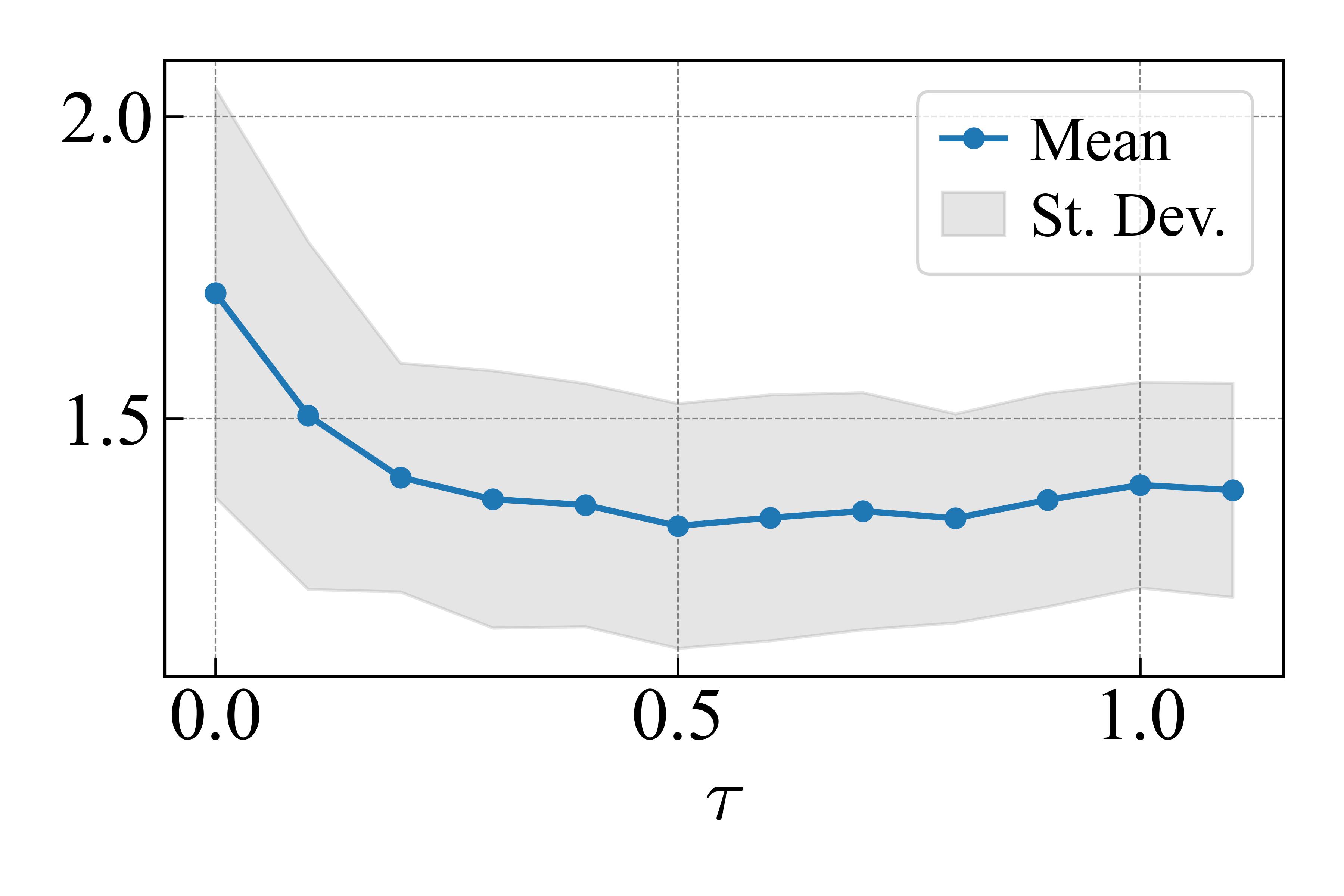}
\end{subfigure} 
\hfill
\begin{subfigure}[b]{0.305\textwidth}
\centering
\text{\small ${\rm tr}\,\mathfrak R$}\\
\includegraphics[trim=0.2cm 0.8cm 0.6cm 0.6cm,clip,width=\textwidth]{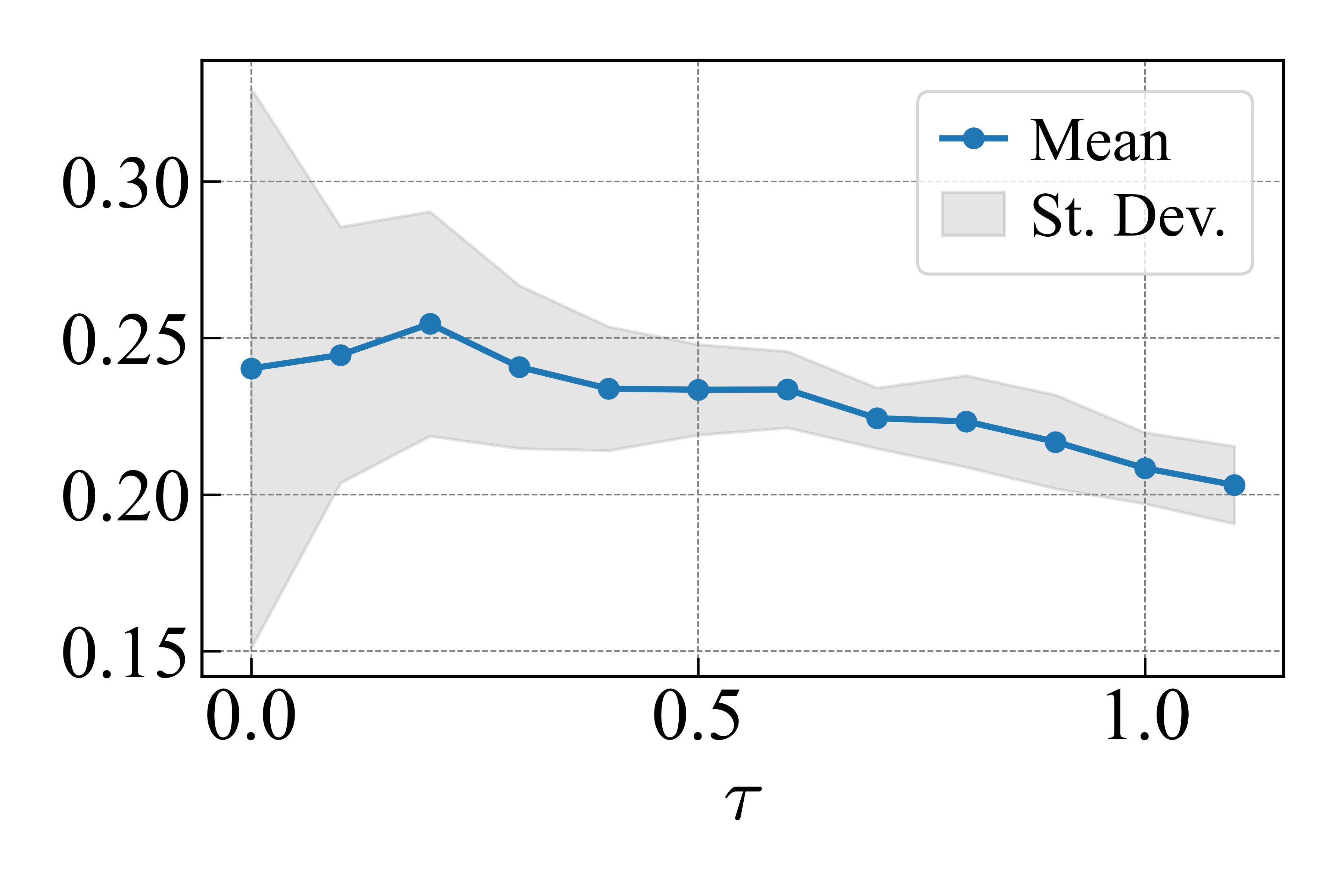}
\end{subfigure} 
\hfill
\vskip8pt
\centering
\begin{subfigure}[b]{0.305\textwidth}
\centering
\text{\small $\rho$}\\
\includegraphics[trim=0.2cm 0.8cm 0.6cm 0.6cm,clip,width=\textwidth]{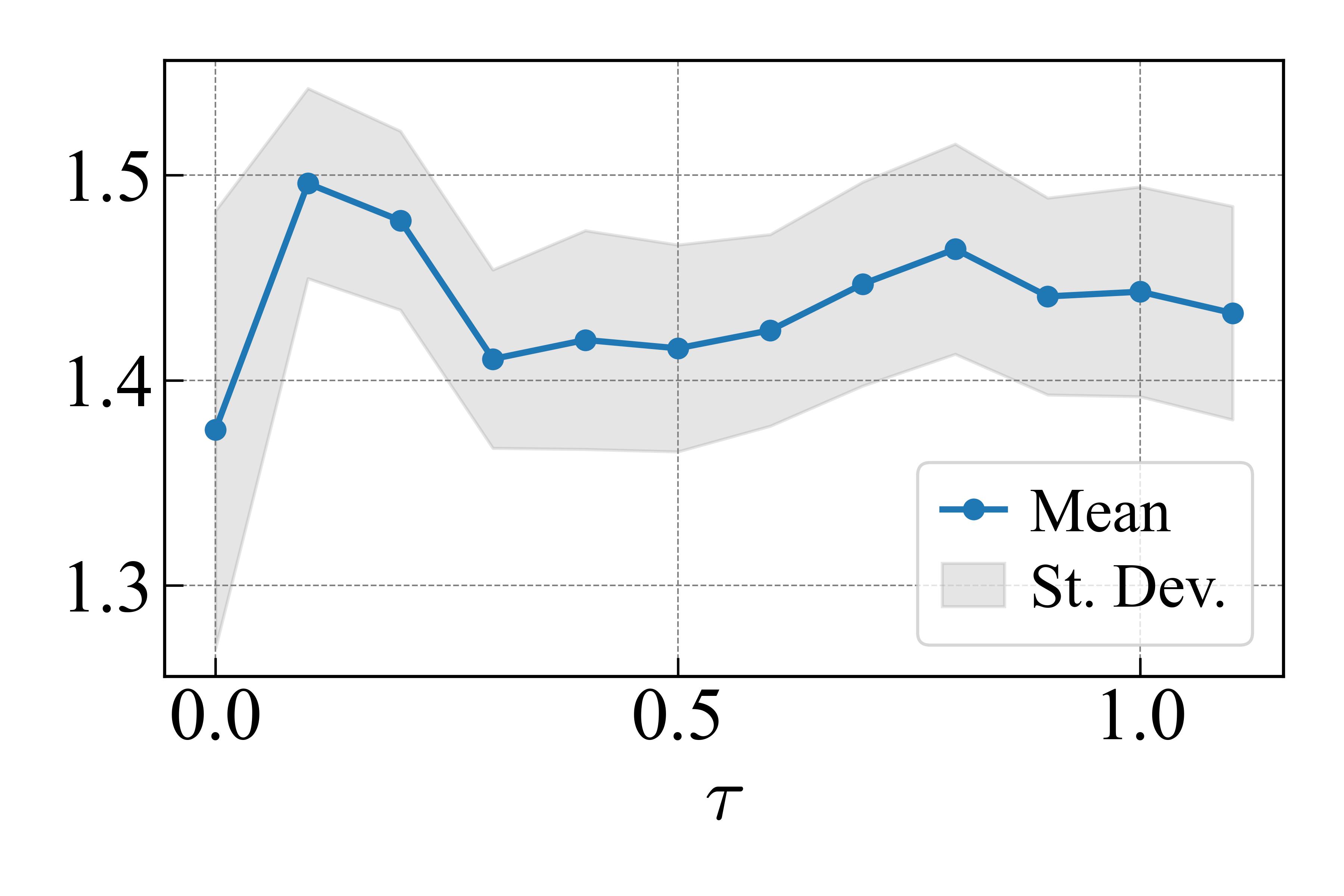}
\end{subfigure}
\hfill
\begin{subfigure}[b]{0.305\textwidth}
\centering
\text{\small $S$}\\
\includegraphics[trim=0.2cm 0.8cm 0.6cm 0.6cm,clip,width=\textwidth]{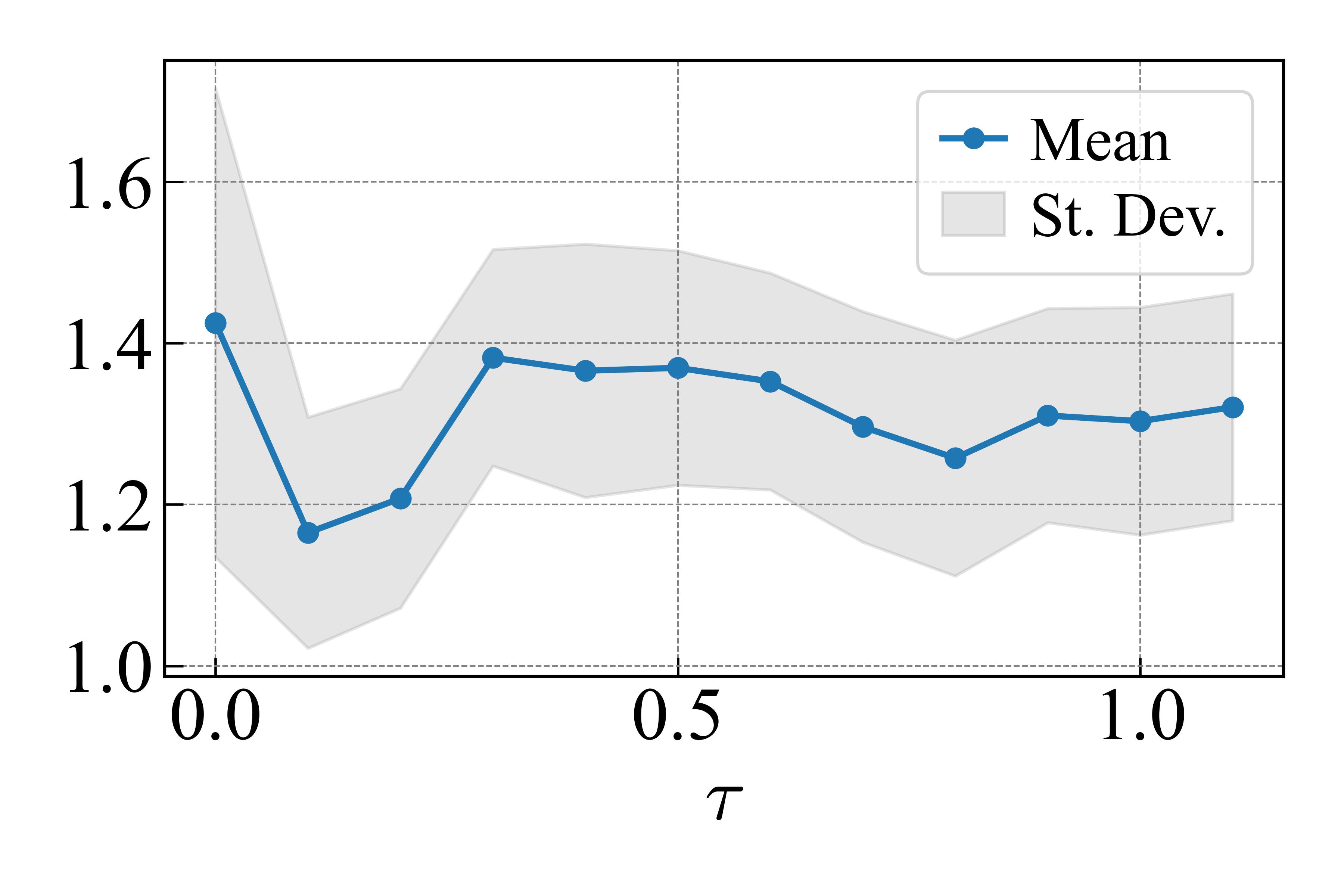}
\end{subfigure} 
\hfill
\begin{subfigure}[b]{0.305\textwidth}
\centering
\text{\small ${\rm tr}\,\mathfrak R$}\\
\includegraphics[trim=0.2cm 0.8cm 0.6cm 0.6cm,clip,width=\textwidth]{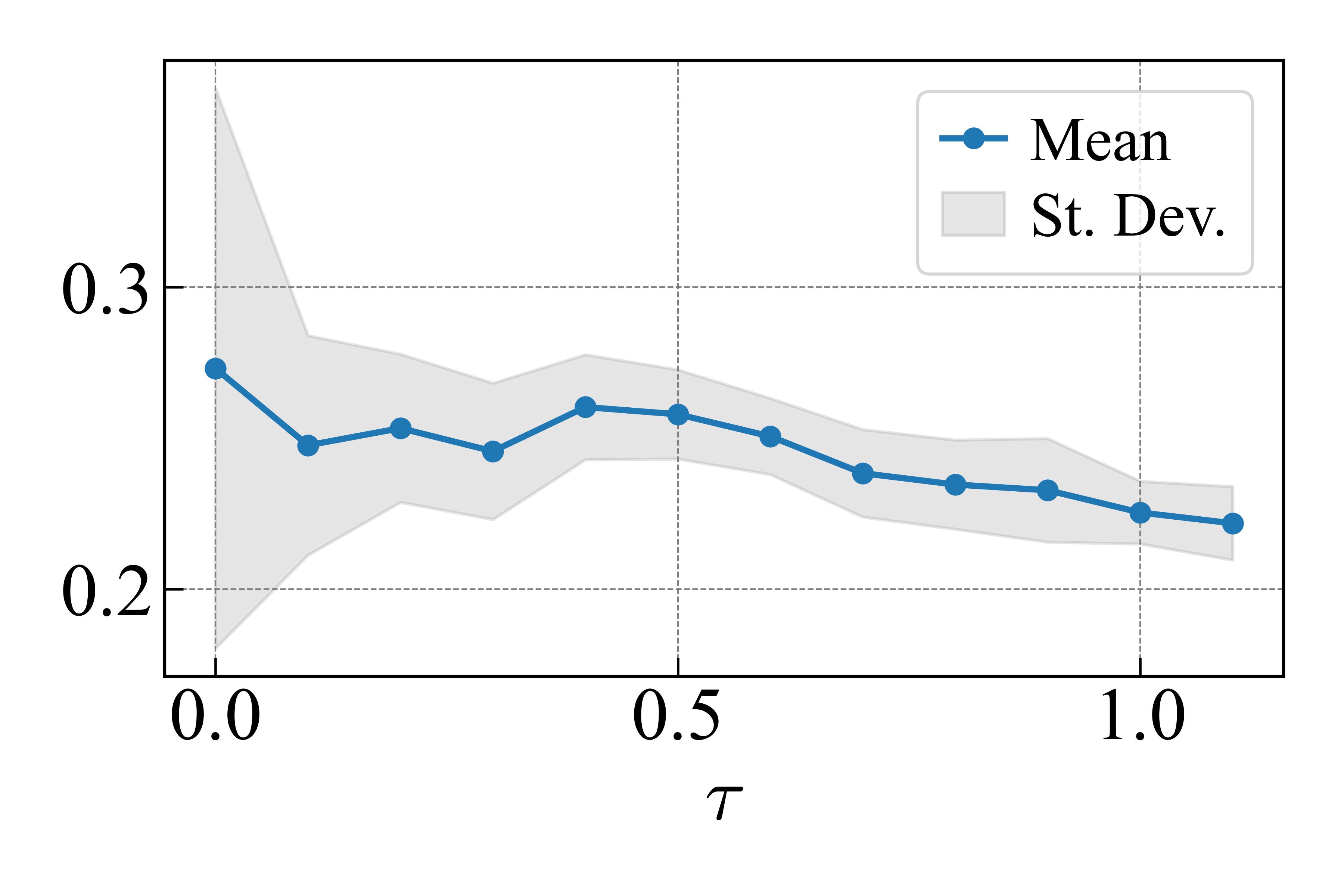}
\end{subfigure} 
\caption{\sf Mean and standard deviation of $\rho$, $S$, and ${\rm tr}\,\mathfrak R$ in $D_1$ (top row) and $D_2$ (bottom row) with 
respect to PDFs of $\,\xbar{\langle\rho\rangle}_5$, $\xbar{\langle S\rangle}_5$, and ${\rm tr}(\,\xbar{\mathfrak R}_5)$, respectively, for 
different $\tau\in[0,1.1]$.\label{fig48}}
\end{figure}

In each spatial window $D_1$ and $D_2,$  we consider statistical properties represented by the PDFs in \fref{fig47} for 
$\,\xbar{\langle\rho\rangle}_5$, $\xbar{\langle S\rangle}_5$, and ${\rm tr}(\,\xbar{\mathfrak R}_5)$. Within each histogram, 
the mean and the standard deviation are computed. The procedure is repeated for several different values of $\tau$ in the range of 
$[0,1.1]$ and reported in \fref{fig48}, showing the influence of the initial interface perturbation $\tau$ on the statistical properties 
of the dissipative Young measure solutions. As one can see, the standard deviation of the PDFs does not seem to depend on $\tau$, except 
for ${\rm tr}(\,\xbar{\mathfrak R}_5)$, for which larger values of $\tau$ yield smaller standard deviation. Understanding the existence of 
the Young measures with non-singular support as a possible evidence of turbulence, we infer that turbulent aspects remain present even 
under strong mixing of the initial flow corresponding to large values of $\tau$.   

It should be observed that the size of the selected spatial windows chosen for these numerical simulations may, in principle, vary, but our numerical experiments (not reported here for the sake of brevity) provided consistent evidence demonstrating unbiased results.

\subsection{Proper Orthogonal Decomposition (POD)}
In this section, we conduct a POD analysis for the computed solutions at the final time $T=2$. We first (\S\ref{sec411}) consider a POD 
for the variable $\rho$ with similar results observed for the other variables. Secondly, we  (\S\ref{sec412}) perform a POD study for the 
Ces\`aro averages of $\rho$, $\mathfrak E$, and ${\rm tr}\,\mathfrak R$.

\subsubsection{POD for $\rho$}\label{sec411}
For each mesh resolution $m$, we  obtain the following center data:
\begin{equation*}
\widetilde\rho_{m,\ell}:=\rho_{m,\ell}(2)-\frac{1}{L}\sum_{i=1}^L\rho_{m,i}(2),\quad\ell=1,\dots,L,
\end{equation*}
where $\rho_{m,\ell}$ are given by \eref{3.1f}. We then reshape each $\widetilde\rho_{m,\ell}$ into a column vector of length $N_m^2$,
construct the data matrices
\begin{equation*}
R_m=\bigl[\widetilde\rho_{m,1}\;\big|\;\widetilde\rho_{m,2}\;\big|\;\cdots\;\big|\;\tilde\rho_{m,L}\bigr]\;\in\;\mathbb R^{N_m^2\times L},
\quad m=1,\dots,M,
\end{equation*}
and perform the singular‐value decomposition $R_m=W_m\Sigma_mV_m^\top$, which yields the singular values $(s_m)_j:=(\Sigma_m)_{jj}$,
$j=1,\dots,L$ (assuming $L<N_m^2$) and orthonormal spatial modes (columns of $W_m$). The modal energies are given by $(s_m)_j^2$, and the
cumulative energy fraction 
\begin{equation*}
{\rm CEF}_m(k):=\frac{\sum_{j=1}^k(s_m)^2_j}{\sum_{j=1}^L(s_m)^2_j}
\end{equation*}
measures the fraction of total variance captured by the first $k$ modes. We then determine the minimal $K_{0.95}$ such that
${\rm CEF}_m(K_{0.95})\ge0.95$, that is, the number of POD modes required to reconstruct any snapshot with at least 95\% of its variance.
Repeating this procedure over all grid resolutions and different values of the interface-perturbation parameter $\tau$, produces the
corresponding values of $K_{0.95}$ presented in \tref{tab41}.
\begin{table}[ht!]
\begin{center}
\begin{tabular}{|c|*{11}{c|}}
\hline
\diagbox{$m$}{$\tau$} & 0.1 & 0.2 & 0.3 & 0.4 & 0.5 & 0.6 & 0.7 & 0.8 & 0.9 & 1.0 & 1.1 \\
\hline
1 & 5  & 9  & 10 & 11 & 14 & 17 & 20 & 22 & 22 & 23 & 24 \\
\hline
2 & 12 & 20 & 27 & 29 & 34 & 36 & 41 & 45 & 45 & 48 & 49 \\
\hline
3 & 26 & 36 & 47 & 47 & 51 & 55 & 59 & 59 & 60 & 61 & 62 \\
\hline
4 & 54 & 63 & 65 & 67 & 69 & 69 & 70 & 70 & 69 & 70 & 70 \\
\hline
5 & 76 & 77 & 78 & 78 & 77 & 77 & 77 & 77 & 76 & 76 & 76 \\
\hline
\end{tabular}
\end{center}
\caption{\sf Number of POD modes $K_{0.95}$ required to capture 95\% of the variance of $\rho$ for each $m$ and $\tau$.\label{tab41}}
\end{table}

This table demonstrates that $K_{0.95}$ depends on $m$ and $\tau$ in different ways. For small values $m=1$ or $2$, the required number of 
modes grows steadily as $\tau$ increases, ranging from only a few modes at small $\tau$ to nearly $50$ modes at large $\tau$. As $m$ 
increases, the number of modes increases sharply: for instance, at $\tau=0.5$, the required number of modes grows from $14$ for $m=1$ to 
$77$ for $m=5$. However, this growth does not continue indefinitely. At larger $m$, the values quickly level off, stabilizing around 
$76$--$78$ modes, almost independently of $\tau$. This indicates a nonlinear, saturating behavior: while both $m$ and $\tau$ contribute to 
the increase in mode count, $m$ has a stronger effect, and once the system reaches a certain complexity, further increases in $\tau$ no 
longer change the dimensional requirements of the reduced-order model.

From a reduced-order modeling perspective, ensuring 95\% reconstruction accuracy across all meshes up to $m=5$ and for all 
$\tau\in[0.1,1.1]$ would require a basis of roughly $80$ modes. If one is interested in coarser simulations, say, for $m=2$, the 
requirement drops to about $20$--$50$ modes even at the largest $\tau$, offering significant savings in basis size. We note that the 
observation that about $70$--$80$ modes are needed to capture 95\% of the variance is consistent with the broadband spectral content 
typical of turbulent flows; see, e.g., \cite{Berkooz1993,Holmes1996}.

\begin{figure}[ht!]
\centering
\begin{subfigure}[b]{0.315\textwidth}
\centering
\hspace*{0.8cm}\text{\small $m=1$}\\
\includegraphics[width=\textwidth]{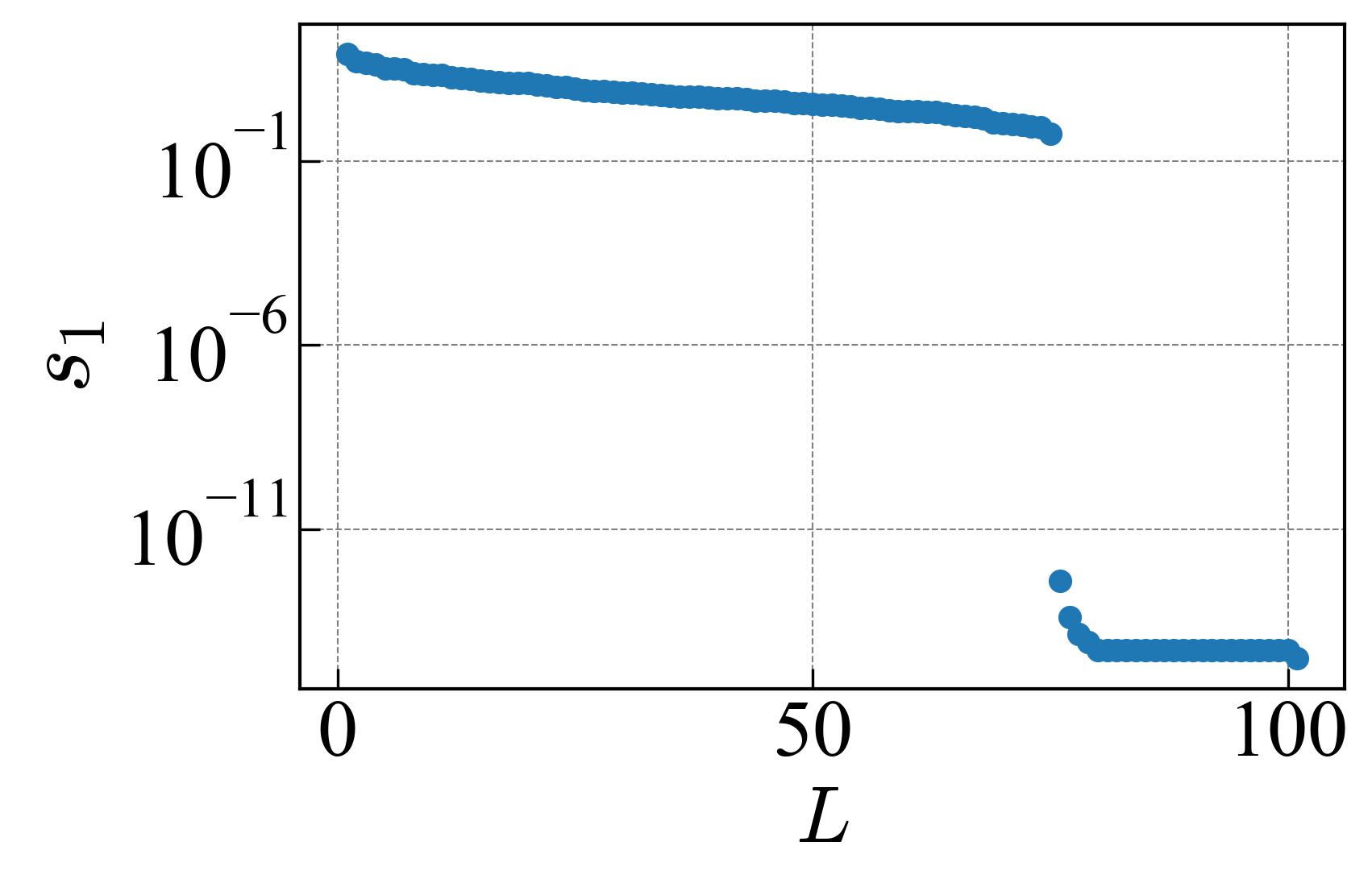}
\end{subfigure}\hspace*{1cm}
\begin{subfigure}[b]{0.315\textwidth}
\centering
\hspace*{0.7cm}\text{\small $m=5$}\\
\includegraphics[width=\textwidth]{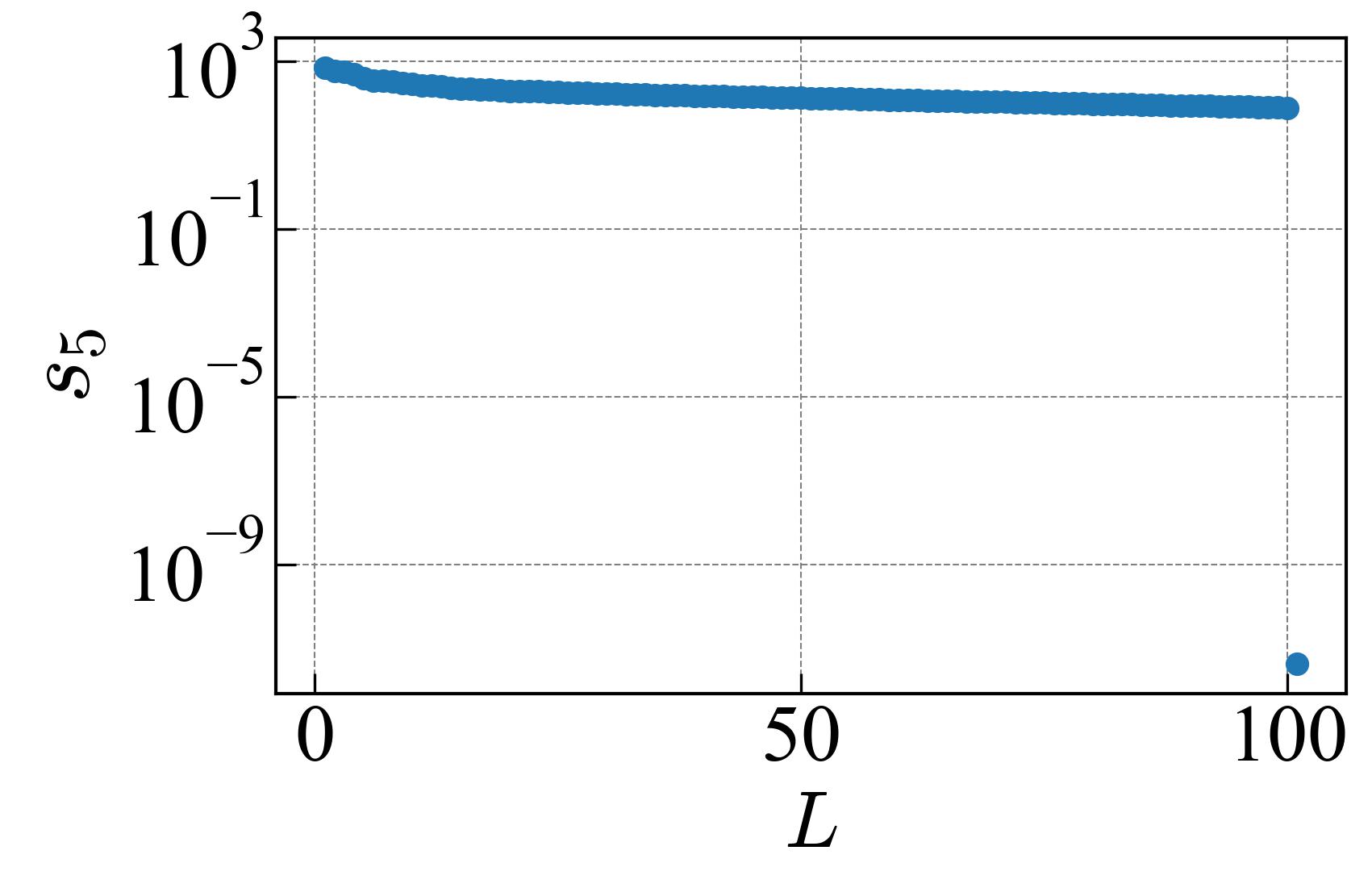}
\end{subfigure} 
\caption{\sf Logarithmic scale decay of the POD singular values for $\tau=1.1$ and $m=1$ (left) and $m=5$ (right).\label{fig49}}
\end{figure}

\fref{fig49} shows the decay of the singular values $(s_m)_j$ for $\tau=1.1$, comparing the coarse ($m=1$) and fine ($m=5$) resolutions. 
In both cases, the spectrum spans many orders of magnitude. On the coarse mesh, the singular values plunge down over the first $75$ modes. 
On the much finer mesh, the singular values decay more gradually, only reaching values around $10^{-9}$ by the hundredth mode, reflecting 
the fact that fine resolution can capture more small-scale features. The slow decay of singular values on finer meshes highlights the wide 
range of active scales---a hallmark of turbulence.

\subsubsection{POD for Ces\`aro Averages}\label{sec412}
We now perform a similar analysis on Ces\`aro averages of $\rho$, $\mathfrak E$, and ${\rm tr}\,\mathfrak R$, and, as before, compute the
minimal $K_{0.95}$ for which ${\rm CEF}_m(K_{0.95})\ge0.95$. The results are reported in Tables \ref{tab42}--\ref{tab44} for
$\,\xbar{\langle\rho\rangle}_5$, $\xbar{\mathfrak E}_5$, and ${\rm tr}(\,\xbar{\mathfrak R}_5)$, respectively. One can observe from Table
\ref{tab42} that significantly fewer POD modes are needed to capture 95\% of the variance of $\,\xbar{\langle\rho\rangle}_5$ compared to the
corresponding data for $\rho$ (Table \ref{tab41}), reflecting a reduction in the effective degrees of freedom. The required number of modes
increases systematically with mesh resolution $m$ and with the perturbation parameter $\tau$, though it tends to saturate once
$\tau\gtrsim0.8$. Among the three considered quantities, density requires the fewest modes, while the energy and (especially) Reynolds
stress defects demand larger modal representations, indicating that defect measures remain more sensitive to small-scale fluctuations even
after averaging. This trend is consistent with turbulence modeling principles, where averaging smooths fine structures reduces effective
complexity while preserving the dominant coherent features of the flow.
\begin{table}[ht!]
\begin{center}
\begin{tabular}{|c|c|c|c|c|c|c|c|c|c|c|c|}
\hline
\diagbox{$M$}{$\tau$} & 0.1 & 0.2 & 0.3 & 0.4 & 0.5 & 0.6 & 0.7 & 0.8 & 0.9 & 1.0 & 1.1 \\
\hline
2 & 11 & 17 & 22 & 24 & 28 & 31 & 36 & 38 & 38 & 41 & 43 \\
\hline
3 & 20 & 27 & 35 & 37 & 40 & 45 & 49 & 50 & 50 & 52 & 52 \\
\hline
4 & 30 & 36 & 44 & 47 & 50 & 54 & 56 & 57 & 57 & 58 & 58 \\
\hline
5 & 45 & 49 & 54 & 57 & 58 & 61 & 62 & 62 & 61 & 62 & 61 \\
\hline
\end{tabular}
\end{center}
\caption{\sf Number of POD modes $K_{0.95}$ required to capture 95\% of the variance of $\,\xbar{\langle\rho\rangle}_5$ for each $m$ and 
$\tau$.\label{tab42}}
\end{table}
\begin{table}[ht!]
\begin{center}
\begin{tabular}{|c|c|c|c|c|c|c|c|c|c|c|c|}
\hline
\diagbox{$M$}{$\tau$} & 0.1 & 0.2 & 0.3 & 0.4 & 0.5 & 0.6 & 0.7 & 0.8 & 0.9 & 1.0 & 1.1 \\
\hline
2 & 12 & 23 & 31 & 33 & 39 & 42 & 46 & 50 & 50 & 53 & 54 \\
\hline
3 & 25 & 36 & 47 & 48 & 53 & 57 & 62 & 63 & 63 & 63 & 63 \\
\hline
4 & 40 & 50 & 57 & 57 & 60 & 65 & 66 & 67 & 66 & 66 & 65 \\
\hline
5 & 54 & 60 & 65 & 64 & 65 & 68 & 68 & 68 & 66 & 66 & 65 \\
\hline
\end{tabular}
\end{center}\caption{\sf Number of POD modes $K_{0.95}$ required to capture 95\% of the variance of $\,\xbar{\mathfrak E}_5$ for each $m$
and $\tau$.\label{tab43}}
\end{table}
\begin{table}[ht!]
\begin{center}
\begin{tabular}{|c|c|c|c|c|c|c|c|c|c|c|c|}
\hline
\diagbox{$M$}{$\tau$} & 0.1 & 0.2 & 0.3 & 0.4 & 0.5 & 0.6 & 0.7 & 0.8 & 0.9 & 1.0 & 1.1 \\
\hline
2 & 13 & 23 & 32 & 34 & 41 & 44 & 48 & 52 & 51 & 55 & 55 \\
\hline
3 & 25 & 38 & 49 & 50 & 55 & 59 & 64 & 65 & 65 & 66 & 65 \\
\hline
4 & 41 & 52 & 59 & 59 & 63 & 67 & 68 & 69 & 68 & 68 & 67 \\
\hline
5 & 56 & 62 & 66 & 66 & 68 & 70 & 71 & 71 & 69 & 69 & 68 \\
\hline
\end{tabular}
\end{center}
\caption{\sf Number of POD modes $K_{0.95}$ required to capture 95\% of the variance of $\,{\rm tr}(\,\xbar{\mathfrak R}_5)$ for each $m$
and $\tau$.\label{tab44}}
\end{table}

\section{Conclusions}\label{sec5}
In this paper, we have investigated the Kelvin-Helmholtz (KH) instability problem for the random compressible Euler equations. Motivated by
the fact that physically reasonable solutions may be seen as inviscid limits of the Navier-Stokes flows, we have investigated random
dissipative weak (DW) solutions of the compressible Euler equations that arise as weak limits of weak solutions of the compressible
Navier-Stokes equations. Motivated by the non-uniqueness of DW solutions, we have adopted a statistical perspective inspired by the
turbulence theory. Our goal was to identify robust features of the solution space that persist across different realizations and mesh
refinements.

To this end, we have computed stable numerical solutions using a stochastic collocation method implemented with the help of a fifth-order
alternative weighted essentially non-oscillatory (A-WENO) scheme and seventh-order central weighted essentially non-oscillatory (CWENO)
interpolation in the random space. We have averaged the computed solutions over several embedded uniform grids and obtained Ces\`aro
averages, whose properties have been investigated using stochastic tools. In particular, we have analyzed Reynolds stresses and energy
defects, probability density functions of averaged quantities, and reduced-order representations via the proper orthogonal decompositions
(POD) analysis.

The numerical experiments have revealed several turbulence-like features. The KH instability produces coherent vortices that roll up and 
subsequently develop fine-scale structures under stronger perturbations, reflecting the transition toward turbulent mixing. Reynolds stress
and energy defect quantify the transport due to unresolved fluctuations and stabilize under mesh refinement, highlighting their role as
turbulence-style observables. Statistical analysis of density and entropy have showed nontrivial probability density functions, emphasizing
the persistent variability characteristic of turbulence. The POD analysis applied to individual solution components has demonstrated a slow
spectral decay and the need for a large number of modes to capture the solution variance, further underlining the broad range of active
scales typical of turbulence. At the same time, POD of Ces\`aro-averaged data requires fewer modes, illustrating how averaging reduces
effective complexity while retaining dominant flow structures.

Overall, our study demonstrates that random KH instabilities can be meaningfully characterized using a turbulence-inspired statistical 
framework. The combination of DW solutions, averaging procedures, and reduced-order models provides a novel approach for describing the
complex and chaotic behavior of inviscid compressible flows. Future research will focus on extending this methodology to more general flow
configurations, exploring long-time dynamics, and further clarifying the role of DW solutions as a statistical model for turbulent
compressible fluid flows.

\begin{acknowledgment}
The work of A. Chertock was supported in part by NSF grant DMS-2208438. The work of A. Kurganov was supported in part by NSFC grants
12171226 and W2431004. M. Herty acknowledges funding by the Deutsche Forschungsgemeinschaft (DFG, German Research Foundation) - SPP 2410 
Hyperbolic Balance Laws in Fluid Mechanics: Complexity, Scales, Randomness (CoScaRa) within Project(s) HE5386/26-1 (Numerische 
Verfahren für gekoppelte Mehrskalenprobleme, 525842915) and (Zufällige kompressible Euler Gleichungen: Numerik und ihre Analysis, 
525853336) HE5386/27-1 as well as though support received funding from the European Union’s Horizon Europe research and innovation 
programme under the Marie Sklodowska-Curie Doctoral Network Datahyking (Grant No. 101072546). M.~Luk\'a\v{c}ov\'a-Medvi{\softd}ov\'a was supported by the
 Deutsche Forschungsgemeinschaft (DFG, German Research Foundation) - SPP 2410  LU1470/10-1  (Random Euler Equations: Numerics and its Analysis, 525853336)  and LU1470/9-1 (An Active Flux method for the Euler equations, 525800857). She was also partially supported by the Gutenberg Research College and
 Mainz Institute of Multiscale Modeling (M3odel).
\end{acknowledgment}

\appendix
\section{Appendix A}\label{appA}
We first introduce the following notation for functional spaces:
$$
\begin{aligned}
C_{\rm weak,loc}\big(L^r(\Td);[0,T)\big):=\Big\{f\in C\big(L^r(\Td);K\big)\mbox{ for any compact set }K\subset[0,T)\mbox{ and }&\\
\int\limits_{\Td}f(\bm x,t)\varphi(\bm x)\,{\rm d}\bm x\in C([0,T))\mbox{ for any }\varphi\in L^{r'}(\Td)\Big\},&
\end{aligned}
$$
where $r'>1$ and $\frac{1}{r}+\frac{1}{r'}=1$. Analogously,
$$
\begin{aligned}
BV_{\rm weak,loc}\big(L^r(\Td);[0,T)\big):=\Big\{f\in C\big(L^r(\Td);K\big)\mbox{ for any compact set }K\subset[0,T)\mbox{ and }&\\
\int\limits_{\Td}f(\bm x,t)\varphi(\bm x)\,{\rm d}\bm x\in BV([0,T))\mbox{ for any }\varphi\in L^{r'}(\Td)\Big\}.&
\end{aligned}
$$
Further, the space ${\cal M}^+(\Td;\R^{d\times d}_{\rm sym})$ is the space of Radon measures ranging in the set of symmetric positive
semi-definite matrices, that is,
$$
\begin{aligned}
&{\cal M}^+(\Td;\R^{d\times d}_{\rm sym}):=\Big\{\mu\in{\cal M}(\Td; \R^{d \times d}_{\rm sym}),\,
\int\limits_{\Td}\phi(\bm\zeta\otimes\bm\zeta):{\rm d}\mu\ge0~~\forall\bm{\zeta}\in\Bbb R^d,~\forall\phi\in C(\Td),~\phi\geq 0\Big\}.
\end{aligned}
$$

A generalized DW solution of the compressible Euler equations is defined in the following way.
\begin{definition}[{\bf DW solution}]\label{def21}
Let the initial data satisfy
$$
\begin{aligned}
&\rho_0\in L^\gamma(\Td),\quad\bm m_0\in L^\frac{2\gamma}{\gamma+1}(\Td;\R^d),\quad S_0\in L^\gamma(\Td)\\
&E_0=E(\rho_0,\bm m_0,S_0),\quad\mbox{and}\quad\int\limits_{\Td}E(\rho_0,\bm m_0,S_0)\,{\rm d}\bm x<\infty,
\end{aligned}
$$
where $S_0(\bm x)$ is the initial total entropy. We say that $(\rho,\bm m,S)$ is a DW solution of the compressible Euler equations in
$\Td\times[0,T)$, $0<T\le\infty$, if the following holds:

\smallskip
\noindent
$\bullet$ {\bf Regularity:}
$$
\begin{aligned}
&\rho\in C_{\rm weak,loc}\big(L^\gamma(\Td);[0,T)\big),\quad\bm m\in C_{\rm weak,loc}\big(L^{\frac{2\gamma}{\gamma+1}}(\Td;\R^d);[0,T)\big),
\\
&S\in L^\infty\big(L^\gamma(\Td);[0,T)\big)\cap BV_{\rm weak,loc}\big(L^\gamma(\Td);[0,T)\big),\\
&\int\limits_{\Td}E(\rho,\bm m, S)(\bm x,t)\,{\rm d}\bm x\le\int\limits_{\Td}E(\rho_0,\bm m_0,S_0)\,{\rm d}\bm x,~\forall t\in[0,T);
\end{aligned}
$$

\smallskip
\noindent
$\bullet$ {\bf Equation of continuity.} The integral identity
\begin{equation}
\int\limits_0^T\int\limits_{\Td}\Big[\rho\varphi_t+\bm m\!\cdot\!\Grad\varphi\Big]{\rm d}\bm x\,{\rm d}t=
-\int\limits_{\Td}\rho_0\varphi(\bm x,0)\,{\rm d}\bm x
\label{2.2}
\end{equation}
for any $\varphi\in C^1_c\big(\Td\times[0,T)\big)$;

\smallskip
\noindent
$\bullet$ {\bf Momentum equation.}  The integral identity
\begin{equation}
\begin{aligned}
&\int\limits_0^T\int\limits_{\Td}\Big[\bm m\!\cdot\!\bfphi+\mathds{1}_{\rho>0}\frac{\bm m\otimes\bm m}{\rho}:\bnabla\bfphi+
p(\rho,S)\bnabla\!\cdot\!\bfphi\Big]{\rm d}\bm x\,{\rm d}t\\
&=\int\limits_0^T\int\limits_{\Td}\Grad\bfphi:{\rm d}\mathfrak R(t)-\int\limits_{\Td}\vm_0\!\cdot\!\bfphi(\bm x,t)\,{\rm d}\bm x
\end{aligned}
\label{2.3}
\end{equation}
for any $\bfphi\in C^1_c\big(\Td\times[0,T);\R^d\big)$, where the Reynolds defect stress reads as
\begin{equation}
\mathfrak R\in L^\infty\big({\cal M}^+\big(\Td;\R^{d\times d}_{\rm sym}\big);[0,T)\big);
\label{2.4}
\end{equation}

\smallskip
\noindent
$\bullet$ {\bf Entropy inequality:}
\begin{equation}
\begin{aligned}
&\int\limits_{\Td}\big[S(\bm x,t_2+)\varphi(\bm x,t_2+)-S(\bm x,t_1-)\varphi(\bm x,t_1-)\big]\,{\rm d}\bm x\\
&\qquad\ge\int\limits_{t_1}^{t_2}\int\limits_{\Td}\Big[S\varphi_t+\Big<{\cal V}_{\bm x,t};\mathds{1}_{\tilde\rho>0}
\Big(\tilde S\,\frac{\tilde{\bm m}}{\tilde\rho}\Big)\Big>\cdot\Grad\varphi\Big]{\rm d}\bm x\,{\rm d}t,\quad S(\bm x,0-)=S_0(\bm x),
\end{aligned}
\label{2.5}
\end{equation}
for any $0\le t_1\le t_2<T$ and any $\varphi\in C^1_c\big(\Td\times[0,T)\big)$, $\varphi\ge0$, where
$\{{\cal V}_{\bm x,t}\}_{(\bm x,t)\in\Td\times(0,T)}$ is a parametrized probability (Young) measure:
\begin{equation}
\begin{aligned}
{\cal V}_{\bm x,t}\in L^\infty\big(\Td\times(0,T);{\cal P}\big(\mathbb R^{d+2}\big)\big),\quad
(\tilde\rho,\tilde{\bm m},\tilde S)^\top\in\mathbb R^{d+2},\\
\left<{\cal V}_{\bm x,t};\tilde\rho\right>=\rho,\quad\left<{\cal V}_{\bm x,t};\tilde{\bm m}\right>=\bm m,\quad
\big<{\cal V}_{\bm x,t};\tilde S\big>=S;
\end{aligned}
\label{2.6}
\end{equation}

\smallskip
\noindent
$\bullet$ {\bf Compatibility of the energy and Reynolds stress defects:}
\begin{equation}
\begin{aligned}
&\int\limits_{\Td}E(\rho_0,\bm m_0,S_0)\,{\rm d}\bm x\ge\int\limits_{\Td}E(\rho,\bm m,S)\,{\rm d}\bm x+
r(d,\gamma)\int\limits_{\Td}{\rm d}\big({\rm tr}\,\mathfrak R(\bm x,t)\big),\\
&r(d,\gamma)=\min\left\{\hf,\frac{1}{d(\gamma-1)}\right\},
\end{aligned}
\label{2.7}
\end{equation}	
for a.a. $t\in(0,T)$.
\end{definition}

\bibliographystyle{siamnodash}
\bibliography{ref}

\end{document}